\documentclass[sn-mathphys-num]{sn-jnl}

\usepackage[T1]{fontenc}
\usepackage[utf8]{inputenc}
\usepackage{lmodern}
\usepackage{amsmath,amssymb,amsthm,mathtools}
\usepackage{mathrsfs}
\usepackage{booktabs}
\usepackage{graphicx}
\usepackage{microtype}
\usepackage{xcolor}
\usepackage{subcaption}
\usepackage{float}

\newcommand{\dd}{\mathrm{d}}
\newcommand{\E}{\mathbb{E}}
\newcommand{\R}{\mathbb{R}}
\newcommand{\orth}{\operatorname{orth}}
\newcommand{\rank}{\operatorname{rank}}

\newcommand{\range}{\operatorname{range}}

\newcommand{\be}{\boldsymbol{e}}

\newcommand{\bz}{\boldsymbol{z}}
\newcommand{\bA}{\mathbf{A}}
\newcommand{\bB}{\mathbf{B}}
\newcommand{\bC}{\mathbf{C}}
\newcommand{\bD}{\mathbf{D}}
\newcommand{\bE}{\mathbf{E}}
\newcommand{\bF}{\mathbf{F}}
\newcommand{\bG}{\mathbf{G}}
\newcommand{\bH}{\mathbf{H}}
\newcommand{\bI}{\mathbf{I}}
\newcommand{\bJ}{\mathbf{J}}
\newcommand{\bK}{\mathbf{K}}
\newcommand{\bL}{\mathbf{L}}
\newcommand{\bM}{\mathbf{M}}
\newcommand{\bN}{\mathbf{N}}
\newcommand{\bP}{\mathbf{P}}
\newcommand{\bQ}{\mathbf{Q}}
\newcommand{\bR}{\mathbf{R}}
\newcommand{\bS}{\mathbf{S}}
\newcommand{\bT}{\mathbf{T}}
\newcommand{\bU}{\mathbf{U}}
\newcommand{\bV}{\mathbf{V}}
\newcommand{\bW}{\mathbf{W}}
\newcommand{\bX}{\mathbf{X}}
\newcommand{\bY}{\mathbf{Y}}
\newcommand{\bZ}{\mathbf{Z}}
\newcommand{\bGamma}{\boldsymbol{\Gamma}}
\newcommand{\bPhi}{\boldsymbol{\Phi}}
\newcommand{\bPi}{\boldsymbol{\Pi}}
\newcommand{\bPsi}{\boldsymbol{\Psi}}
\newcommand{\bSigma}{\boldsymbol{\Sigma}}
\newcommand{\bLambda}{\boldsymbol{\Lambda}}
\newcommand{\bTheta}{\boldsymbol{\Theta}}
\newcommand{\bXi}{\boldsymbol{\Xi}}

\newcommand{\intervalle}[4]{\mathopen{#1}#2\mathclose{},#3\mathclose{#4}}
\newcommand{\ff}[2]{\intervalle{[}{#1}{#2}{]}}

\theoremstyle{thmstyleone}%
\newtheorem{theorem}{Theorem}
\newtheorem{proposition}{Proposition}
\newtheorem{assumption}{Assumption}
\newtheorem{lemma}{Lemma}
\newtheorem{corollary}{Corollary}

\theoremstyle{thmstyletwo}%
\newtheorem{remark}{Remark}%

\theoremstyle{thmstylethree}%

\begin{document}

\title[Article Title]{Robust BUG dynamical low-rank integrators based on the Euler--Maruyama method for stochastic differential equations. Part I: a closed-moment approach}


\author*[1,2]{\fnm{S\'ebastien} \sur{Riffaud}}\email{sebastien.riffaud@epfl.ch}

\affil[1]{\orgdiv{CSQI Chair}, \orgname{\'Ecole Polytechnique F\'ed\'erale de Lausanne}, \city{1015 Lausanne}, \country{Switzerland}}
\affil[2]{\orgdiv{Laboratory of Applied Mechanical Design}, \orgname{\'Ecole Polytechnique F\'ed\'erale de Lausanne}, \city{1015 Lausanne}, \country{Switzerland}}


\abstract{In this work, we introduce Euler--Maruyama-based Basis-Update \& Galerkin (EM--BUG) dynamical low-rank integrators for stochastic differential equations (SDEs). These integrators extend the standard BUG integrator to the Euler--Maruyama scheme by incorporating diffusion directions into the augmented physical basis and exploiting the second moment to update the physical basis independently of the new Brownian increments. This construction avoids the inversion of the reduced second-moment matrix and allows us to derive convergence estimates that are robust with respect to small eigenvalues. We establish preservation of positive semidefiniteness, exactness with respect to Euler--Maruyama at full rank, and second-moment stability. Strong and weak convergence estimates show that the time-discretization accuracy of the underlying Euler--Maruyama scheme is preserved up to the low-rank approximation error. We also derive Monte Carlo error estimates recovering the standard sampling rate up to the time-discretization and low-rank approximation errors. Numerical experiments on linear SDEs with additive and multiplicative noise support these theoretical results and demonstrate accuracy comparable to the standard BUG integrator.}




\maketitle

\section{Introduction}
\label{sec:introduction}

Stochastic differential equations (SDEs) arise in many applications involving uncertainty. Their numerical simulation can be computationally expensive, since the estimation of statistical quantities often relies on Monte Carlo methods, which require computing a large number of trajectories driven by independent realizations of the noise. In addition, if each trajectory is also high-dimensional, the overall computational cost can become prohibitive.

Dynamical low-rank approximation (DLRA) has been successfully applied in various contexts to reduce the computational cost of high-dimensional problems. Examples include kinetic equations \cite{einkemmer2018vlasov,einkemmer2021vlasov,einkemmer2025review}, random PDEs \cite{babaee2017robust,musharbash2018dual,feppon2018dynamically,musharbash2020symplectic,patil2020real,nobile2025petrov}, and data assimilation \cite{nobile2025kalman,trindade2027enkf,kazashi2026filters,marzouk2026smoothing}. The main idea is to approximate the solution by a low-rank representation. Consider an SDE of the form
\begin{equation}
\label{eq:intro-sde}
\dd X(t)=a(t,X(t))\,\dd t+b(t,X(t))\,\dd W(t),
\end{equation}
where $X(t)\in\mathbb{R}^d$, $W(t)$ is a Wiener process, and $a$ and $b$ denote the drift and diffusion coefficients, respectively. Denoting by $\mathbf{X}(t)\in\mathbb{R}^{d\times N}$ the matrix collecting $N$ Monte Carlo realizations $\{X^{(j)}(t)\}_{j=1}^N$, we seek an approximation of the form
\begin{equation}
\mathbf{X}(t)\approx\mathbf{U}(t)\mathbf{Y}(t)^\top,
\end{equation}
where $\mathbf{U}(t)\in\mathbb{R}^{d\times r}$ has orthonormal columns, $\mathbf{Y}(t)\in\mathbb{R}^{N\times r}$, and $r\ll\min\{d,N\}$. Originally introduced in \cite{koch2007dynamical} for matrix differential equations, the evolution of the low-rank approximation is obtained by projecting the dynamics onto the tangent space of the rank-$r$ manifold. DLRA methods have also been investigated for stochastic differential equations \cite{cao2018stochastic,kazashi2025dynamical}. In particular, the time discretization and Monte Carlo approximation of these methods have been analyzed in \cite{kazashi2026time,kazashi2026sampling}, while alternative continuous-time formulations have been investigated in \cite{kazashi2026further}.

In the stochastic setting, the Dynamically Orthogonal (DO) formulation of \cite{sapsis2009dynamically} provides evolution equations for the deterministic basis and the stochastic coefficients. By imposing the dynamical orthogonality condition $\mathbf{U}(t)^\top\dot{\mathbf{U}}(t)=0$, the evolution equations are given by
\begin{equation}
\left\{
\begin{aligned}
\dot{\mathbf{U}}(t)&=\left(\mathbf{I}-\mathbf{U}(t)\mathbf{U}(t)^\top\right)\mathbb{E}\left[a(t,\mathbf{U}(t)Y(t))Y(t)^\top\right]C_Y(t)^{-1},\\
\dd Y(t)&=\mathbf{U}(t)^\top a(t,\mathbf{U}(t)Y(t))\,\dd t+\mathbf{U}(t)^\top b(t,\mathbf{U}(t)Y(t))\,\dd W(t),
\end{aligned}
\right.
\end{equation}
where $Y(t)\in\mathbb{R}^r$ denotes the vector of stochastic coefficients, whose realizations are collected in $\mathbf{Y}(t)$, and $C_Y(t)=\mathbb{E}[Y(t)Y(t)^\top]$ denotes its second-moment matrix.

The numerical integration of the DO equations can be challenging due to the presence of the inverse of the reduced second-moment matrix $C_Y(t)$. Indeed, when some eigenvalues of $C_Y(t)$ become small, the equations can be ill-conditioned, leading to severe restrictions on the time step, since the stability conditions of standard time-integration schemes may depend on the smallest eigenvalue of $C_Y(t)$.

To overcome this difficulty, several DLRA integrators that are robust with respect to small singular values have been developed in the deterministic setting, such as the projector-splitting integrator \cite{lubich2014projector,kieri2016discretized} and the unconventional integrator \cite{ceruti2022unconventional}. These methods avoid the inversion of small singular values and provide robust error bounds independent of their magnitude. Among them, the Basis-Update \& Galerkin (BUG) integrator \cite{ceruti2022rank,ceruti2024parallel} is particularly attractive in practice, since it does not involve backward time-integration step, allows for rank adaptivity, and admits
high-order extensions \cite{ceruti2022rank,ceruti2024parallel,ceruti2024robust,nobile2026high}.

However, extending the BUG integrator to SDEs while preserving its robustness and accuracy is not straightforward. First, unlike in the deterministic setting, the trajectories evolve through stochastic increments, which must be taken into account in the basis updates and the Galerkin step. Second, the direct application of the BUG integrator to a Monte Carlo ensemble would make the updated physical basis depend on the new Brownian increments. This dependence complicates the derivation of strong and weak error estimates, since the projections and rank truncations interact with the stochastic increments.

To address these difficulties, we propose BUG integrators for SDEs that extend the standard BUG formulation to the Euler--Maruyama method. The key idea of these Euler--Maruyama-based BUG (EM--BUG) integrators is to exploit the second moment. First, the physical basis is augmented with additional directions associated with the diffusion terms of the Euler--Maruyama prediction. Second, for linear SDEs, the second-moment equation is closed and deterministic, allowing us to perform the Galerkin projection and rank truncation independently of the new Brownian increments. The stochastic coefficients are then updated by projecting the Euler--Maruyama prediction onto the resulting physical basis. For nonlinear SDEs, where the second-moment equation is generally not closed, we also propose a natural extension based on the empirical conditional second moment.

This construction allows us to establish several theoretical properties for the closed-moment EM--BUG integrator, including preservation of positive semidefiniteness of the second moment, exactness with respect to the Euler--Maruyama scheme at full rank, and second-moment stability. We also derive strong and weak convergence estimates showing that the time-discretization accuracy of the underlying Euler--Maruyama scheme is preserved up to the low-rank approximation error. In particular, these estimates are robust with respect to small eigenvalues of the reduced second-moment matrix. Additionally, we establish Monte Carlo error estimates that recover the standard sampling rate of order $N^{-1/2}$, up to the time-discretization and low-rank approximation errors.

In this paper, we focus on the closed-moment EM--BUG integrator for linear SDEs, for which we provide a detailed theoretical analysis. The empirical conditional-moment EM--BUG integrator for general nonlinear SDEs will be analyzed separately in a second part.

The remainder of the paper is organized as follows. Section~\ref{sec:setting} introduces the problem setting and notation. Then, Section~\ref{sec:method} presents the closed-moment EM--BUG integrator. Section~\ref{sec:analysis} establishes its theoretical properties and convergence estimates. In Section~\ref{sec:numerical-experiments}, these theoretical results, especially regarding the accuracy and robustness, are supported by several numerical experiments. Finally, Section~\ref{sec:conclusion} draws some conclusions and perspectives.

\section{Problem setting and notation}
\label{sec:setting}

\subsection{Stochastic differential equation and Monte Carlo ensemble}
\label{sec:setting-sde}

We first introduce the setting and notation used throughout the paper. Although the analysis focuses on linear stochastic differential equations, we present the general formulation, which will also be used to introduce the extension to nonlinear problems.

Let $(\Omega,\mathcal{F},(\mathcal{F}_t)_{t\geq 0},\mathbb{P})$ be a filtered probability space and let $W(t)$ be an $m$-dimensional standard $(\mathcal{F}_t)$-Brownian motion. We consider the stochastic differential equation (SDE)
\begin{equation}
\label{eq:sde}
\dd X(t) = a(t,X(t))\,\dd t + b(t,X(t))\,\dd W(t), \qquad X(t)\in L^2(\Omega;\R^d),
\end{equation}
where
\begin{equation*}
a:\ff{0}{T}\times\R^d\to\R^d, \qquad b:\ff{0}{T}\times\R^d\to\R^{d\times m}.
\end{equation*}
We assume that $X(t)$, $a(t,X(t))$, and $b(t,X(t))$ are square integrable for $t\in[0,T]$, so that the moments and expectations introduced below are well defined.

We consider $N$ Monte Carlo realizations of \eqref{eq:sde}, driven by independent Brownian motions $\{W^{(j)}(t)\}_{j=1}^N$, and collect them in the ensemble matrix
\begin{equation}
\label{eq:ensemble}
\bX(t) = \begin{bmatrix} X^{(1)}(t) & \cdots & X^{(N)}(t) \end{bmatrix} \in \R^{d\times N}.
\end{equation}

The non-centered second-moment matrix is defined by
\begin{equation}
\label{eq:second-moment}
\bC(t) = \E[X(t)X(t)^\top].
\end{equation}
According to It\^o's formula, it satisfies
\begin{equation}
\label{eq:second-moment-equation}
\dot \bC(t) = \E\!\left[a(t,X(t))X(t)^\top + X(t)a(t,X(t))^\top + b(t,X(t))b(t,X(t))^\top\right].
\end{equation}
In general, the second-moment equation \eqref{eq:second-moment-equation} is not closed, since its right-hand side cannot be expressed only in terms of $\bC(t)$. However, for linear SDEs, the second moment satisfies a closed evolution equation. This property will be exploited in Section~\ref{sec:bug-closed} to construct the closed-moment EM--BUG integrator.

For the Monte Carlo ensemble $\bX(t)$, we also define the empirical non-centered second moment by
\begin{equation}
\label{eq:empirical-second-moment}
\bC_N(t) = \frac{1}{N}\bX(t)\bX(t)^\top.
\end{equation}

\subsection{Euler--Maruyama scheme}
\label{sec:setting-em}

Let $0=t_0<t_1<\cdots<t_{N_t}=T$ be a uniform time grid with step size
$\Delta t=T/N_t$. For $n\in\{0,\ldots,N_t-1\}$, we denote the Brownian increments by
\begin{equation}
\label{eq:brownian-increments}
\Delta W_n^{(j)} = W^{(j)}(t_{n+1})-W^{(j)}(t_n), \qquad j=1,\ldots,N.
\end{equation}
For each realization, the Euler--Maruyama discretization of \eqref{eq:sde} is given by
\begin{equation}
\label{eq:euler-maruyama}
X_{n+1}^{\mathrm{EM},(j)}=X_n^{\mathrm{EM},(j)}+\Delta t\,a(t_n,X_n^{\mathrm{EM},(j)})+b(t_n,X_n^{\mathrm{EM},(j)})\Delta W_n^{(j)},
\qquad j=1,\ldots,N.
\end{equation}

Collecting the $N$ realizations gives the Euler--Maruyama ensemble
\begin{equation}
\label{eq:em-ensemble}
\bX_n^{\mathrm{EM}}=
\begin{bmatrix}
X_n^{\mathrm{EM},(1)} & \cdots & X_n^{\mathrm{EM},(N)}
\end{bmatrix}
\in\R^{d\times N}.
\end{equation}
This Euler--Maruyama update can be written in matrix form as
\begin{equation}
\label{eq:em-matrix}
\bX_{n+1}^{\mathrm{EM}}=\bX_n^{\mathrm{EM}}+\Delta t\,\bA_n(\bX_n^{\mathrm{EM}})+\Delta\bB_n(\bX_n^{\mathrm{EM}}),
\end{equation}
where
\begin{equation}
\label{eq:em-drift-ensemble}
\bA_n(\bX_n^{\mathrm{EM}})=\begin{bmatrix}a(t_n,X_n^{\mathrm{EM},(1)}) & \cdots & a(t_n,X_n^{\mathrm{EM},(N)})\end{bmatrix},
\end{equation}
and
\begin{equation}
\label{eq:em-diffusion-increment}
\Delta\bB_n(\bX_n^{\mathrm{EM}})=\begin{bmatrix}b(t_n,X_n^{\mathrm{EM},(1)})\Delta W_n^{(1)} & \cdots & b(t_n,X_n^{\mathrm{EM},(N)})\Delta W_n^{(N)}\end{bmatrix}.
\end{equation}

Let $X_n^{\mathrm{EM}}$ denote a generic Euler--Maruyama particle at time $t_n$. Its non-centered second moment is defined by
\begin{equation}
\bC_n^{\mathrm{EM}}=\mathbb{E}\!\left[X_n^{\mathrm{EM}}(X_n^{\mathrm{EM}})^\top\right].
\end{equation}
The Euler--Maruyama discretization then induces the second-moment update
\begin{equation}
\label{eq:em-second-moment}
\begin{aligned}
\bC_{n+1}^{\mathrm{EM}}&=\E\!\left[X_{n+1}^{\mathrm{EM}}(X_{n+1}^{\mathrm{EM}})^\top\right]\\
&=\bC_n^{\mathrm{EM}}+\Delta t\,\E\!\left[a_n^{\mathrm{EM}}(X_n^{\mathrm{EM}})^\top+X_n^{\mathrm{EM}}(a_n^{\mathrm{EM}})^\top+b_n^{\mathrm{EM}}(b_n^{\mathrm{EM}})^\top\right]\\
&\qquad+\Delta t^2\,\E[a_n^{\mathrm{EM}}(a_n^{\mathrm{EM}})^\top],
\end{aligned}
\end{equation}
where $a_n^{\mathrm{EM}}=a(t_n,X_n^{\mathrm{EM}})$ and $b_n^{\mathrm{EM}}=b(t_n,X_n^{\mathrm{EM}})$. Compared with the forward Euler discretization of \eqref{eq:second-moment-equation}, the Euler--Maruyama second-moment update \eqref{eq:em-second-moment} contains the additional term $\Delta t^2\,\E[a_n^{\mathrm{EM}}(a_n^{\mathrm{EM}})^\top]$, which is essential for preserving positive semidefiniteness (see Proposition~\ref{prop:positivity}).

\subsection{Dynamical low-rank approximation}
\label{sec:setting-low-rank}

At time $t_n$, the Euler--Maruyama ensemble $\bX_n^{\mathrm{EM}}$ is approximated by
\begin{equation}
\label{eq:low-rank-approximation}
\bX_n^{\mathrm{BUG}} = \bU_n\bY_n^\top,
\end{equation}
where $\bU_n\in\R^{d\times r}$ has orthonormal columns, $\bY_n\in\R^{N\times r}$ is not necessarily full rank, and $r\leq\min\{d,N\}$ denotes the prescribed rank.

The empirical second moment of the low-rank ensemble $\bX_n^{\mathrm{BUG}}$ is given by
\begin{equation}
\label{eq:lr-empirical-second-moment}
\bC_{N,n}^{\mathrm{BUG}} = \frac{1}{N}\bX_n^{\mathrm{BUG}}(\bX_n^{\mathrm{BUG}})^\top = \bU_n\left(\frac{1}{N}\bY_n^\top\bY_n\right)\bU_n^\top.
\end{equation}
For the closed-moment EM--BUG integrator, the second moment is approximated separately using the low-rank representation
\begin{equation}
\label{eq:lr-second-moment}
\bC_n^{\mathrm{BUG}}=\bU_n\bD_n\bU_n^\top,
\end{equation}
where $\bD_n\in\R^{r\times r}$ is symmetric positive semidefinite. The matrix $\bD_n$ is evolved independently using the closed second-moment equation and, in general, $\bD_n\neq\bY_n^\top\bY_n/N$. In contrast, the empirical conditional-moment extension described in Section~\ref{sec:bug-empirical} uses $\bD_n=\bY_n^\top\bY_n/N$.

\section{BUG integrators based on the Euler--Maruyama method}
\label{sec:method}

We first introduce an Euler–Maruyama-based BUG (EM--BUG) integrator that exploits the closed evolution equation of the second moment for linear SDEs. Then, we present an extension to general nonlinear SDEs based on the empirical conditional second moment.

\subsection{Closed-moment EM--BUG integrator for linear SDEs}
\label{sec:bug-closed}

We first consider the linear SDE
\begin{equation}
\label{eq:linear-sde}
\dd X(t)=\bA X(t)\,\dd t+\bB_0\,\dd W_0(t)+\sum_{\ell=1}^q \bB_\ell X(t)\,\dd W_\ell(t),
\end{equation}
where $W_0(t)$ is a $m_0$-dimensional Brownian motion and $\{W_\ell(t)\}_{\ell=1}^q$ are independent scalar Brownian motions, also independent of $W_0(t)$. The Euler--Maruyama discretization of \eqref{eq:linear-sde} is
\begin{equation}
\label{eq:closed-em}
\begin{aligned}
\bX_{n+1}^{\mathrm{EM}}&=\bX_n^{\mathrm{EM}}+\Delta t\,\bA\bX_n^{\mathrm{EM}}+\bB_0\Delta\bW_{0,n}+\sum_{\ell=1}^q\bB_\ell\bX_n^{\mathrm{EM}}\Delta\bW_{\ell,n},\\
&=\bM\bX_n^{\mathrm{EM}}+\bB_0\Delta\bW_{0,n}+\sum_{\ell=1}^q\bB_\ell\bX_n^{\mathrm{EM}}\Delta\bW_{\ell,n},
\end{aligned}
\end{equation}
where $\bM=\bI+\Delta t\,\bA\in\R^{d\times d}$, $\Delta\bW_{0,n}=[\Delta W_{0,n}^{(1)},\ldots,\Delta W_{0,n}^{(N)}]\in\R^{m_0\times N}$, and $\Delta\bW_{\ell,n}=\operatorname{diag}(\Delta W_{\ell,n}^{(1)},\ldots,\Delta W_{\ell,n}^{(N)})\in\R^{N\times N}$. The second moment satisfies the closed equation
\begin{equation}
\label{eq:linear-second-moment}
\dot \bC(t)=\bA\bC+\bC\bA^\top+\bB_0\bB_0^\top+\sum_{\ell=1}^q \bB_\ell \bC\bB_\ell^\top.
\end{equation}
According to the Euler--Maruyama discretization, the closed second-moment update is
\begin{equation}
\label{eq:closed-moment-map}
\begin{aligned}
\bC_{n+1}^{\mathrm{EM}}&=\bC_n^{\mathrm{EM}}+\Delta t\left(\bA\bC_n^{\mathrm{EM}}+\bC_n^{\mathrm{EM}}\bA^\top+\bB_0\bB_0^\top+\sum_{\ell=1}^q\bB_\ell\bC_n^{\mathrm{EM}}\bB_\ell^\top\right)+(\Delta t)^2\bA\bC_n^{\mathrm{EM}}\bA^\top\\
&=\bM\bC_n^{\mathrm{EM}}\bM^\top+\Delta t\,\bB_0\bB_0^\top+\Delta t\sum_{\ell=1}^q\bB_\ell\bC_n^{\mathrm{EM}}\bB_\ell^\top.
\end{aligned}
\end{equation}

In the BUG integrator \cite{ceruti2022rank,ceruti2024parallel}, originally developed for deterministic matrix differential equations, the solution is approximated by the low-rank factorization $\bU_n\bS_n\bV_n^\top$. In the time-explicit case, the augmented bases $\widehat{\bU}_n$ and $\widehat{\bV}_n$ are first constructed to represent the tangent-space projection of the forward Euler prediction. The augmented core matrix $\widehat{\bS}_n$ is then updated by Galerkin projection of the forward Euler prediction onto the augmented bases. Finally, the resulting approximation is truncated back to rank $r$. Several modifications are required to adapt this construction to SDEs:
\begin{enumerate}
\item \textbf{Stochastic representation and Galerkin projection.} Unlike deterministic matrix differential equations, the ensemble $\bX(t)$ of SDEs is generally not differentiable in time. To avoid deriving continuous-time evolution equations for the stochastic factors $\bS_n$ and $\bV_n$, which would require accounting for their It\^o dynamics, we use the low-rank approximation $\bU_n\bY_n^\top$, where $\bY_n$ combines these two factors. The stochastic coefficients $\bY_n$ are then updated by Galerkin projection of the discrete Euler--Maruyama prediction onto the physical basis.
\item \textbf{Second-moment formulation.} Once the update for $\bY_n$ is defined, it remains to determine how to update $\bU_n$. In the Dynamically Orthogonal (DO) approach, the evolution equation for the physical basis is given by
\begin{equation}
\label{eq:do-basis}
\dot{\bU}=(\bI-\bU\bU^\top)\E[a(t,\bU Y)Y^\top]\left(\E[YY^\top]\right)^{-1},
\end{equation}
where $Y(t)\in\R^r$ denotes the stochastic coefficients of a generic realization. However, this equation has two drawbacks. First, it involves the inverse of the reduced second moment, which may become ill-conditioned when some eigenvalues are small. Second, it does not account for the directions generated by the noise. To overcome these difficulties, we exploit the fact that $\bU$ is also the physical basis of the second moment. Consequently, instead of using the DO equation, we consider the second-moment equation \eqref{eq:linear-second-moment}, which includes the contribution of the diffusion. Moreover, since this equation is deterministic, we can apply the standard BUG approach to avoid inverting the reduced second moment. More precisely, we construct the augmented physical basis $\widehat{\bU}_{n+1}$ to approximate the Euler--Maruyama second-moment prediction, and we compute the augmented reduced second moment $\widehat{\bD}_{n+1}$ by Galerkin projection of the Euler--Maruyama second-moment prediction onto this basis.
\item \textbf{Physical basis augmentation.} Due to the second-moment formulation, the augmented physical basis is constructed to approximate the Euler--Maruyama second-moment prediction. However, since our primary objective is to approximate the ensemble $\bX_n^{\mathrm{EM}}$, the augmented physical basis must also approximate its Euler--Maruyama prediction. Moreover, compared to the forward Euler prediction used in the standard BUG integrator, the Euler--Maruyama ensemble prediction contains additional contributions due to the noise. For these reasons, we construct the augmented physical basis $\widehat{\bU}_{n+1}$ to approximate the Euler--Maruyama ensemble and second-moment predictions, while including diffusion directions to account for the noise contributions.
\item \textbf{Moment-based truncation.} Finally, instead of computing the augmented stochastic coefficients $\widehat{\bY}_n$ by Galerkin projection and then truncating the resulting approximation, we perform the truncation on the augmented reduced second moment $\widehat{\bD}_n$. This allows us to update the physical basis independently of the new Brownian increments. Moreover, since the empirical second moment converges to the second moment as the number of particles increases, truncating the eigenvalues of $\widehat{\bD}_n$ is consistent with truncating the squared singular values of $\widehat{\bY}_n$ in the large-sample limit.
\end{enumerate}

In the standard BUG integrator, the initial factors $(\bU_0,\bS_0,\bV_0)$ are obtained from the rank-$r$ truncated SVD of the initial ensemble $\bX_0$. For the closed-moment BUG integrator, we keep $\bU_0$ and set $\bY_0=\bV_0\bS_0^\top$. In addition, we need to initialize the reduced second moment $\bD_0$. Two choices are possible. The first is the consistent initialization
\begin{equation}
\label{eq:closed-initialization}
\bD_0=\frac1N\bY_0^\top\bY_0,
\end{equation}
which ensures that the initial second moment and the empirical second moment coincide, i.e., $\bC_0^{\mathrm{BUG}}=\bC_{N,0}^{\mathrm{BUG}}$. Alternatively, when the exact initial second moment $\bC(0)$ is available, we can set $\bD_0=\bU_0^\top\bC(0)\bU_0$. This initialization exploits the exact second-moment information, but does not ensure in general that $\bC_0^{\mathrm{BUG}}=\bC_{N,0}^{\mathrm{BUG}}$.

Then, given $(\bU_n,\bD_n,\bY_n)$, one time step of the closed-moment EM--BUG integrator reads as follows.

\paragraph{Step 1: physical basis augmentation.}
Let us define the Euler--Maruyama ensemble prediction
\begin{equation}
\label{eq:closed-stochastic-prediction}
\bZ_{n+1}=\bM\bX_n^{\mathrm{BUG}}+\bB_0\Delta\bW_{0,n}+\sum_{\ell=1}^q\bB_\ell\bX_n^{\mathrm{BUG}}\Delta\bW_{\ell,n},
\end{equation}
and the Euler--Maruyama second-moment prediction
\begin{equation}
\label{eq:closed-moment-prediction}
\bS_{n+1}=\bM\bC_n^{\mathrm{BUG}}\bM^\top+\Delta t\,\bB_0\bB_0^\top+\Delta t\sum_{\ell=1}^q\bB_\ell\bC_n^{\mathrm{BUG}}\bB_\ell^\top.
\end{equation}
The augmented physical basis is constructed to approximate the Euler--Maruyama predictions $\bZ_{n+1}$ and $\bS_{n+1}$:
\begin{equation}
\label{eq:closed-basis-augmentation}
\widehat{\bU}_{n+1}=\orth
\begin{bmatrix}
\bM\bU_n & \bQ_0 & \bB_1\bU_n & \cdots & \bB_q\bU_n
\end{bmatrix}\in\R^{d\times\hat r},
\end{equation}
where $\bQ_0\in\R^{d\times q_0}$ contains additive-noise directions. It can be obtained from a basis of $\range(\bB_0)$ when $q_0=\rank(\bB_0)$, or from the dominant eigenvectors of $\bB_0\bB_0^\top$ when fewer directions are retained. When all additive-noise directions are included, the augmented basis represents exactly the Euler--Maruyama predictions $\bZ_{n+1}$ and $\bS_{n+1}$. Otherwise, the partial augmentation may introduce an additional approximation error associated with the omitted additive-noise directions.

\paragraph{Step 2: second-moment update.}
The augmented reduced second moment is computed by Galerkin projection of the Euler--Maruyama second-moment prediction \eqref{eq:closed-moment-prediction} onto the augmented basis:
\begin{equation}
\label{eq:closed-moment-update}
\widehat{\bD}_{n+1}=\widehat{\bU}_{n+1}^\top\bS_{n+1}\widehat{\bU}_{n+1}\in\R^{\hat r\times\hat r}.
\end{equation}

\paragraph{Step 3: truncation step.}
Given the rank-$r$ truncated eigendecomposition
\begin{equation}
\label{eq:closed-spectral-decomposition}
\widehat{\bD}_{n+1}\approx\bPhi_{n+1}\bLambda_{n+1}\bPhi_{n+1}^\top,
\end{equation}
the physical basis and reduced second moment are given by
\begin{equation}
\label{eq:closed-moment-truncation}
\bU_{n+1}=\widehat{\bU}_{n+1}\bPhi_{n+1},
\qquad
\bD_{n+1}=\bLambda_{n+1}.
\end{equation}

\paragraph{Step 4: stochastic update.}
The stochastic coefficients are updated by Galerkin projection of the Euler--Maruyama ensemble prediction \eqref{eq:closed-stochastic-prediction} onto the updated physical basis:
\begin{equation}
\label{eq:closed-stochastic-update}
\bY_{n+1}=\bZ_{n+1}^\top\bU_{n+1}.
\end{equation}

\subsection{Extension to general nonlinear SDEs: an empirical conditional-moment approach}
\label{sec:bug-empirical}

The closed-moment EM--BUG construction suggests a natural extension to general SDEs, for which the second-moment equation \eqref{eq:second-moment-equation} is not closed. The idea is to replace the closed second-moment prediction by the empirical conditional second moment of the Euler--Maruyama ensemble prediction
\begin{equation}
\label{eq:empirical-stochastic-prediction}
\bZ_{n+1}=\bX_n^{\mathrm{BUG}}+\Delta t\,\bA_n(\bX_n^{\mathrm{BUG}})+\Delta\bB_n(\bX_n^{\mathrm{BUG}})
\end{equation}
given by
\begin{equation}
\label{eq:empirical-moment-prediction}
\bS_{n+1}=\E\left[\frac{1}{N}\bZ_{n+1}\bZ_{n+1}^\top\mid\mathcal{F}_{t_n}\right]=\frac{1}{N}\bF_n\bF_n^\top+\Delta t\,\bG_n,
\end{equation}
where
\begin{equation*}
\bF_n=\bX_n^{\mathrm{BUG}}+\Delta t\,\bA_n(\bX_n^{\mathrm{BUG}}),\quad
\bG_n=\frac{1}{N}\sum_{j=1}^N b(t_n,X_n^{\mathrm{BUG},(j)})b(t_n,X_n^{\mathrm{BUG},(j)})^\top.
\end{equation*}

Given $(\bU_n,\bY_n)$, one time step of the empirical conditional-moment EM--BUG integrator reads as follows.

\paragraph{Step 1: physical basis augmentation.}

Let $\bV_n$ be the stochastic basis, given by the completed thin QR factorization
\begin{equation}
\label{eq:completed-qr}
\bY_n=\bV_n\bR_n,
\end{equation}
where $\bV_n\in\R^{N\times r}$ has orthonormal columns and $\bR_n\in\R^{r\times r}$. If $\bY_n$ is rank deficient, the orthonormal basis of $\range(\bY_n)$ is completed to $r$ columns.

The augmented physical basis is then constructed to represent exactly the tangent-space projection of the deterministic contribution to the Euler--Maruyama ensemble prediction $\bZ_{n+1}$, while allowing the inclusion of noise directions:
\begin{equation}
\label{eq:empirical-basis-augmentation}
\widehat{\bU}_{n+1}=\orth
\begin{bmatrix}
\bU_n & \bA_n(\bX_n^{\mathrm{BUG}})\bV_n & \bQ_n
\end{bmatrix},
\end{equation}
where $\bQ_n$ contains relevant noise directions. It can be obtained from a basis of $\range(\bG_n)$, or from the dominant eigenvectors of $\bG_n$ when fewer directions are retained.

\paragraph{Step 2: second-moment update.}
The augmented reduced second moment is computed by Galerkin projection of the empirical conditional second-moment prediction \eqref{eq:empirical-moment-prediction} onto the augmented basis:
\begin{equation}
\label{eq:empirical-moment-update}
\widehat{\bD}_{n+1}=\widehat{\bU}_{n+1}^\top\bS_{n+1}\widehat{\bU}_{n+1}\in\R^{\hat r\times\hat r}.
\end{equation}

\paragraph{Step 3: truncation step.}
Given the rank-$r$ truncated eigendecomposition
\begin{equation}
\label{eq:empirical-spectral-decomposition}
\widehat{\bD}_{n+1}\approx\bPhi_{n+1}\bLambda_{n+1}\bPhi_{n+1}^\top,
\end{equation}
the physical basis is given by
\begin{equation}
\label{eq:empirical-moment-truncation}
\bU_{n+1}=\widehat{\bU}_{n+1}\bPhi_{n+1}.
\end{equation}

\paragraph{Step 4: stochastic update.}
The stochastic coefficients are updated by Galerkin projection of the Euler--Maruyama ensemble prediction \eqref{eq:empirical-stochastic-prediction} onto the updated physical basis:
\begin{equation}
\label{eq:empirical-stochastic-update}
\bY_{n+1}=\bZ_{n+1}^\top\bU_{n+1}.
\end{equation}

Note that the rank $r$ is fixed here for simplicity, but an adaptive rank can also be used by choosing the number of eigenvalues retained in the truncation step according to a prescribed tolerance.

\section{Analysis}
\label{sec:analysis}

In this section, we establish the main properties of the closed-moment EM--BUG integrator: preservation of positive semidefiniteness, exactness with respect to the Euler--Maruyama scheme at full rank, second-moment stability, strong and weak convergence, and Monte Carlo error estimates. The analysis is conducted for a fixed rank $r$.

Let $\bP_{n+1}=\bU_{n+1}\bU_{n+1}^\top$ denote the orthogonal projector associated with the updated physical basis. The low-rank ensemble is obtained by projecting the Euler--Maruyama ensemble prediction onto this basis:
\begin{equation}
\label{eq:analysis-trajectory-structure}
\bX_{n+1}^{\mathrm{BUG}}=\bP_{n+1}\bZ_{n+1}.
\end{equation}
The resulting second moment satisfies
\begin{equation}
\label{eq:analysis-moment-structure}
\bC_{n+1}^{\mathrm{BUG}}=\bP_{n+1}\bS_{n+1}\bP_{n+1}=\bU_{n+1}\bD_{n+1}\bU_{n+1}^\top.
\end{equation}

For an ensemble matrix $\bX$, we also introduce the normalized Frobenius norm
\begin{equation*}
\|\bX\|_N = \frac{1}{\sqrt{N}}\|\bX\|_{\mathrm F},
\end{equation*}
so that estimates expressed in this norm do not introduce an artificial factor depending on $N$.

\subsection{Positivity of the second moment}
\label{sec:analysis-positivity}

\begin{proposition}[Positive semidefiniteness]
\label{prop:positivity}
Assume that $\bD_0\succeq 0$. Then, the closed-moment EM--BUG integrator preserves positive semidefiniteness of the reduced second moment:
\begin{equation}
\label{eq:positivity-linear}
\bD_n\succeq 0,\qquad \bC_n^{\mathrm{BUG}}=\bU_n\bD_n\bU_n^\top\succeq 0,
\end{equation}
for all $n\geq0$. In particular, this property does not require any restriction on $\Delta t>0$.
\end{proposition}

\begin{proof}
We proceed by induction. At $n=0$, the assumption $\bD_0\succeq 0$ implies $\bC_0^{\mathrm{BUG}}=\bU_0\bD_0\bU_0^\top\succeq 0$. Suppose that $\bD_n\succeq 0$, and hence $\bC_n^{\mathrm{BUG}}\succeq 0$. Since the Euler--Maruyama second-moment prediction \eqref{eq:closed-moment-prediction} is a sum of positive semidefinite matrices, we have $\bS_{n+1}\succeq 0$. Consequently, its Galerkin projection $\widehat{\bD}_{n+1}=\widehat{\bU}_{n+1}^\top\bS_{n+1}\widehat{\bU}_{n+1}$ is also positive semidefinite. Since the spectral truncation retains only nonnegative eigenvalues, we obtain $\bD_{n+1}=\bLambda_{n+1}\succeq 0$, and hence $\bC_{n+1}^{\mathrm{BUG}}=\bU_{n+1}\bD_{n+1}\bU_{n+1}^\top\succeq 0$, which concludes the induction.
\end{proof}

\begin{remark}
This property holds because the additional term $(\Delta t)^2\bA\bC_n^{\mathrm{BUG}}\bA^\top$ completes the factorization $\bM\bC_n^{\mathrm{BUG}}\bM^\top$ in \eqref{eq:closed-moment-prediction}. Without this term, the standard explicit Euler discretization of the second-moment equation remains symmetric but is not necessarily positive semidefinite for all $\Delta t>0$.
\end{remark}

\begin{remark}
This property ensures that all eigenvalues of $\widehat{\bD}_{n+1}$ are nonnegative. Consequently, the spectral truncation can be used to obtain the best rank-$r$ approximation of $\widehat{\bD}_{n+1}$ and update the physical basis.
\end{remark}

\subsection{Exactness with respect to Euler--Maruyama}
\label{sec:analysis-exactness}

\begin{proposition}[Euler--Maruyama exactness at full rank]
Consider the closed-moment EM--BUG integrator with complete augmentation, where the augmented physical basis contains all additive-noise directions. Assume that $r=d$ and $N\geq d$. If the closed-moment EM--BUG integrator and the Euler--Maruyama scheme are initialized with the same ensemble and use the same Wiener increments, then
\begin{equation}
\bX_n^{\mathrm{BUG}}=\bX_n^{\mathrm{EM}},
\qquad
\bC_{N,n}^{\mathrm{BUG}}=\bC_{N,n}^{\mathrm{EM}},
\end{equation}
almost surely for all $n\geq0$. If, additionally, the closed-moment EM--BUG integrator and the Euler--Maruyama scheme are initialized with the same second moment, then
\begin{equation}
\bC_n^{\mathrm{BUG}}=\bC_n^{\mathrm{EM}},
\end{equation}
for all $n\geq0$.
\end{proposition}

\begin{proof}
We proceed by induction. The assertions hold at $n=0$ by assumption. First, suppose that $\bX_n^{\mathrm{BUG}}=\bX_n^{\mathrm{EM}}$. Since both schemes use the same Wiener increments, they yield the same Euler--Maruyama prediction $\bZ_{n+1}=\bX_{n+1}^{\mathrm{EM}}$. Since all additive and multiplicative noise directions are included in the augmented basis, we have $\range(\bZ_{n+1})\subseteq\range(\widehat{\bU}_{n+1})$. Moreover, since $r=d$, the rank-$d$ truncation does not discard any direction of the augmented basis. Consequently,
\begin{equation*}
\bX_{n+1}^{\mathrm{BUG}}=\bP_{n+1}\bZ_{n+1}=\bZ_{n+1}=\bX_{n+1}^{\mathrm{EM}}.
\end{equation*}
It follows that
\begin{equation*}
\bC_{N,n+1}^{\mathrm{BUG}}=\frac{1}{N}\bX_{n+1}^{\mathrm{BUG}}(\bX_{n+1}^{\mathrm{BUG}})^\top=\frac{1}{N}\bX_{n+1}^{\mathrm{EM}}(\bX_{n+1}^{\mathrm{EM}})^\top=\bC_{N,n+1}^{\mathrm{EM}}.
\end{equation*}
Then, suppose additionally that $\bC_n^{\mathrm{BUG}}=\bC_n^{\mathrm{EM}}$. The closed second-moment recursion gives $\bS_{n+1}=\bC_{n+1}^{\mathrm{EM}}$. Moreover, since the augmentation is complete and no direction is discarded, we obtain
\begin{equation*}
\bC_{n+1}^{\mathrm{BUG}}=\bP_{n+1}\bS_{n+1}\bP_{n+1}=\bS_{n+1}=\bC_{n+1}^{\mathrm{EM}},
\end{equation*}
which concludes the induction.
\end{proof}

\begin{proposition}[Euler--Maruyama exactness at reduced rank]
\label{prop:exactness}
Consider the closed-moment EM--BUG integrator with complete augmentation. Assume that $r<d$ and the Euler--Maruyama particles are identically distributed. Assume additionally that the closed-moment EM--BUG integrator and the Euler--Maruyama scheme are initialized with the same ensemble and second moment, and use the same Wiener increments. If, at each time step, the Euler--Maruyama second moment has rank at most $r$:
\begin{equation}
\label{eq:exactness-condition-moment}
\operatorname{rank}(\bC_n^{\mathrm{EM}})\le r,
\end{equation}
then the closed-moment EM--BUG integrator reproduces exactly the Euler--Maruyama ensemble, its empirical second moments, and the Euler--Maruyama second moment:
\begin{equation}
\label{eq:exactness-solution}
\bX_n^{\mathrm{BUG}}=\bX_n^{\mathrm{EM}}, \qquad
\bC_{N,n}^{\mathrm{BUG}}=\bC_{N,n}^{\mathrm{EM}}, \qquad
\bC_n^{\mathrm{BUG}}=\bC_n^{\mathrm{EM}},
\end{equation}
almost surely for all $n\geq0$.
\end{proposition}

\begin{proof}
We proceed by induction. The assertions hold at $n=0$ by assumption. Suppose that $\bX_n^{\mathrm{BUG}}=\bX_n^{\mathrm{EM}}$ and $\bC_n^{\mathrm{BUG}}=\bC_n^{\mathrm{EM}}$. Since both schemes use the same Wiener increments, they yield the same Euler--Maruyama prediction $\bZ_{n+1}=\bX_{n+1}^{\mathrm{EM}}$. Moreover, the closed second-moment recursion gives $\bS_{n+1}=\bC_{n+1}^{\mathrm{EM}}$. Since the augmentation is complete, we have $\range(\bS_{n+1})\subseteq\range(\widehat{\bU}_{n+1})$. Under the rank condition \eqref{eq:exactness-condition-moment}, the spectral truncation discards no nonzero eigenvalues. Consequently, $\range(\bS_{n+1})\subseteq\range(\bU_{n+1})$, and hence
\begin{equation*}
\bC_{n+1}^{\mathrm{BUG}}=\bP_{n+1}\bS_{n+1}\bP_{n+1}=\bS_{n+1}=\bC_{n+1}^{\mathrm{EM}}.
\end{equation*}
Since the Euler--Maruyama particles are identically distributed, each particle $Z_{n+1}^{(j)}$ has second moment $\bS_{n+1}$, and we have $\range(\bZ_{n+1})\subseteq\range(\bS_{n+1})$ almost surely. Consequently, $\range(\bZ_{n+1})\subseteq\range(\bU_{n+1})$ almost surely, and hence
\begin{equation*}
\bX_{n+1}^{\mathrm{BUG}}=\bP_{n+1}\bZ_{n+1}=\bZ_{n+1}=\bX_{n+1}^{\mathrm{EM}}.
\end{equation*}
It follows that
\begin{equation*}
\bC_{N,n+1}^{\mathrm{BUG}}=\frac{1}{N}\bX_{n+1}^{\mathrm{BUG}}(\bX_{n+1}^{\mathrm{BUG}})^\top=\frac{1}{N}\bX_{n+1}^{\mathrm{EM}}(\bX_{n+1}^{\mathrm{EM}})^\top=\bC_{N,n+1}^{\mathrm{EM}},
\end{equation*}
which concludes the induction.
\end{proof}

\begin{corollary}[Sufficient conditions for Euler--Maruyama exactness at reduced rank]
Consider the closed-moment EM--BUG integrator with complete augmentation. Assume that $r<d\leq N$ and the Euler--Maruyama particles are independent and identically distributed. Assume additionally that the closed-moment EM--BUG integrator and the Euler--Maruyama scheme are initialized with the same ensemble and second moment, and use the same Wiener increments. If, at each time step, the Euler--Maruyama ensemble has rank at most $r$:
\begin{equation}
\label{eq:exactness-condition-ensemble}
\operatorname{rank}(\bX_n^{\mathrm{EM}})\leq r,\qquad\text{almost surely},
\end{equation}
then the closed-moment EM--BUG integrator reproduces exactly the Euler--Maruyama ensemble, its empirical second moments, and the Euler--Maruyama second moment
\begin{equation}
\bX_n^{\mathrm{BUG}}=\bX_n^{\mathrm{EM}}, \qquad
\bC_{N,n}^{\mathrm{BUG}}=\bC_{N,n}^{\mathrm{EM}}, \qquad
\bC_n^{\mathrm{BUG}}=\bC_n^{\mathrm{EM}},
\end{equation}
almost surely for all $n\geq0$.
\end{corollary}

\begin{proof}
Since $N\geq d>r$, we have $N\geq r+1$. Moreover, since the Euler--Maruyama particles are independent and identically distributed, the condition $\rank(\bX_n^{\mathrm{EM}})\leq r$ almost surely implies that their common distribution is supported on a subspace of dimension at most $r$. Indeed, otherwise, $r+1$ independent particles would be linearly independent with positive probability. Consequently, $\rank(\bC_n^{\mathrm{EM}})\leq r$, and the result follows from Proposition~\ref{prop:exactness}.
\end{proof}

\subsection{Second-moment stability}
\label{sec:analysis-stability}

The closed-moment EM--BUG integrator inherits the second-moment stability bound of the Euler--Maruyama scheme, without any additional time-step restriction. This estimate follows from the fact that orthogonal projection does not increase the trace of the positive semidefinite second-moment prediction.

\begin{proposition}[Second-moment stability]
Assume that $\bC_0^{\mathrm{BUG}}\succeq 0$. Let $\bR=\bM^\top\bM+\Delta t\sum_{\ell=1}^{q}\bB_\ell^\top\bB_\ell$, and assume that $\bR\preceq\rho\bI$ for some $0\leq\rho<1$. Then, the closed-moment EM--BUG integrator satisfies
\begin{equation}
\label{eq:stability}
\operatorname{tr}(\bC_n^{\mathrm{BUG}})\leq \rho^n\operatorname{tr}(\bC_0^{\mathrm{BUG}})+\Delta t\,\frac{1-\rho^n}{1-\rho}\|\bB_0\|_{\mathrm F}^2,
\end{equation}
for all $0\leq n\leq N_t$. In particular, apart from $\bC_0^{\mathrm{BUG}}$, this bound is independent of the rank $r$, the number of particles $N$, the smallest eigenvalue of the reduced second moments, and spectral gaps.
\end{proposition}

\begin{proof}
Since $\bS_{n+1}\succeq 0$ and $\bP_{n+1}$ is an orthogonal projector, we have
\begin{equation*}
\label{eq:stability-projection}
\operatorname{tr}(\bP_{n+1}\bS_{n+1}\bP_{n+1})\leq\operatorname{tr}(\bS_{n+1}).
\end{equation*}
Moreover, since $\bC_n^{\mathrm{BUG}}\succeq0$ and the Löwner order is preserved under congruence transformations, the assumption $\bR\preceq\rho\bI$ implies
\begin{equation*}
(\bC_n^{\mathrm{BUG}})^{1/2}\bR(\bC_n^{\mathrm{BUG}})^{1/2}\preceq\rho\bC_n^{\mathrm{BUG}}.
\end{equation*}
Consequently, taking the trace and using its cyclicity yields
\begin{equation*}
\operatorname{tr}(\bR\bC_n^{\mathrm{BUG}})\leq\rho\operatorname{tr}(\bC_n^{\mathrm{BUG}}).
\end{equation*}
Applying these inequalities to the Euler--Maruyama second-moment prediction, we obtain
\begin{equation*}
\begin{aligned}
\operatorname{tr}(\bC_{n+1}^{\mathrm{BUG}})
&=\operatorname{tr}(\bP_{n+1}\bS_{n+1}\bP_{n+1})\leq\operatorname{tr}(\bS_{n+1})\\
&=\operatorname{tr}(\bM\bC_n^{\mathrm{BUG}}\bM^\top)+\Delta t\,\|\bB_0\|_{\mathrm F}^2+\Delta t\sum_{\ell=1}^{q}\operatorname{tr}(\bB_\ell\bC_n^{\mathrm{BUG}}\bB_\ell^\top)\\
&=\operatorname{tr}(\bR\bC_n^{\mathrm{BUG}})+\Delta t\,\|\bB_0\|_{\mathrm F}^2\\
&\leq\rho\operatorname{tr}(\bC_n^{\mathrm{BUG}})+\Delta t\,\|\bB_0\|_{\mathrm F}^2.
\end{aligned}
\end{equation*}
Finally, iterating this recursion yields \eqref{eq:stability}.
\end{proof}

\begin{remark}
The stability condition $\bR\preceq\rho\bI$ with $0\leq\rho<1$ can be made explicit for the two linear problems considered in Section \ref{sec:numerical-experiments}. For the Ornstein--Uhlenbeck process with additive noise, if $\bA=\bA^\top\prec0$, then
\begin{equation*}
\rho=\|\bI+\Delta t\,\bA\|_2^2=\max_{1\leq i\leq d}|1+\Delta t\,\lambda_i(\bA)|^2<1
\qquad\text{for}\qquad
0<\Delta t<\frac{2}{|\lambda_{\min}(\bA)|}.
\end{equation*}
For commuting multiplicative noise, if $\bA=\bPsi\operatorname{diag}(\lambda_1,\ldots,\lambda_d)\bPsi^\top$ and $\bB_\ell=\bPsi\operatorname{diag}(\mu_{1\ell},\ldots,\mu_{d\ell})\bPsi^\top$ with $\bPsi$ orthogonal, then
\begin{equation*}
\rho=\max_{1\leq i\leq d}\left((1+\Delta t\,\lambda_i)^2+\Delta t\sum_{\ell=1}^{q}\mu_{i\ell}^2\right),
\end{equation*}
and the stability condition holds if
\begin{equation*}
2\lambda_i+\sum_{\ell=1}^{q}\mu_{i\ell}^2+\Delta t\,\lambda_i^2<0,
\qquad
i=1,\ldots,d.
\end{equation*}
More generally, a sufficient condition for $\bR\preceq\rho\bI$ with $0\leq\rho<1$ is the existence of constants $\kappa>0$ and $L_A>0$ such that
\begin{equation*}
\bA+\bA^\top+\sum_{\ell=1}^{q}\bB_\ell^\top\bB_\ell\preceq-\kappa\bI,
\qquad
\bA^\top\bA\preceq L_A^2\bI.
\end{equation*}
In this case,
\begin{equation*}
\rho=1-\kappa\Delta t+L_A^2(\Delta t)^2<1
\qquad\text{for}\qquad
0<\Delta t<\frac{\kappa}{L_A^2}.
\end{equation*}
This general sufficient condition may be more restrictive than the exact conditions obtained above.
\end{remark}

\subsection{Strong convergence}
\label{sec:analysis-strong}

We now establish the strong convergence estimate of the closed-moment EM--BUG integrator. We recall that, under suitable assumptions on the coefficients, including temporal regularity (see \cite{kloeden1992numerical}), the Euler--Maruyama scheme satisfies
\begin{equation}
\label{eq:analysis-strong-em}
\left(\E\left[\max_{0\le n\le N_t}\|\bX_n^{\mathrm{EM}}-\bX(t_n)\|_N^2\right]\right)^{1/2}\leq C_T^{\mathrm{EM}}(\Delta t)^\alpha,
\end{equation}
where $C_T^{\mathrm{EM}}$ is independent of $N$ and $\Delta t$. The convergence order is $\alpha=1/2$ in general, and $\alpha=1$ for linear SDEs with constant coefficients and additive noise. We also define the initial low-rank approximation error
\begin{equation}
\label{eq:analysis-strong-initial}
\delta_{r,N}:=\left(\E\left[\|\bX_0^{\mathrm{BUG}}-\bX_0^{\mathrm{EM}}\|_N^2\right]\right)^{1/2},
\end{equation}
and the cumulative discarded energy
\begin{equation}
\label{eq:analysis-strong-error}
\gamma_r^2:=\sum_{n=0}^{N_t-1}\E\left[\operatorname{tr}\left(\bS_{n+1}-\bP_{n+1}\bS_{n+1}\bP_{n+1}\right)\right].
\end{equation}

\begin{proposition}[Uniform bound on $\gamma_r^2$]
Let $\bN_n=\bB_0\bB_0^\top+\sum_{\ell=1}^q\bB_\ell\bC_n^{\mathrm{BUG}}\bB_\ell^\top$, and let $\bP_n^{\mathrm{M}}$ denote the orthogonal projector onto $\range(\bM\bU_n)$. Assume that $\bD_0\succeq0$. Then,
\begin{equation}
\operatorname{tr}\left(\bS_{n+1}-\bP_{n+1}\bS_{n+1}\bP_{n+1}\right)\leq \Delta t\,\operatorname{tr}\left((\bI-\bP_n^{\mathrm{M}})\bN_n\right)\leq \Delta t\,\operatorname{tr}(\bN_n),
\end{equation}
for all $0\leq n<N_t$. Consequently, the cumulative discarded energy satisfies
\begin{equation}
\gamma_r^2\leq C_{T,r,N}^{\mathrm{bound}},
\end{equation}
where $C_{T,r,N}^{\mathrm{bound}}$ is independent of the time-step size $\Delta t$, the smallest eigenvalue of the reduced second moment, and spectral gaps. In addition, if the initial ensemble $\bX_0^{\mathrm{BUG}}$ is obtained from the rank-$r$ truncated SVD of $\bX_0^{\mathrm{EM}}$, the reduced second moment is initialized as $\bD_0=\bY_0^\top\bY_0/N$, and $\E[\bC_{N,0}^{\mathrm{EM}}]=\bC_0^{\mathrm{EM}}$, then
\begin{equation}
\gamma_r^2\leq C_T^{\mathrm{bound}},
\end{equation}
where $C_T^{\mathrm{bound}}$ is independent of the time-step size $\Delta t$, the rank $r$, the number of particles $N$, the smallest eigenvalue of the reduced second moment, and spectral gaps.
\end{proposition}

\begin{proof}
Since $\bD_0\succeq0$, the closed-moment EM--BUG integrator preserves positive semidefiniteness according to Proposition~\ref{prop:positivity}, and hence $\bC_n^{\mathrm{BUG}}\succeq0$, $\bN_n\succeq0$, and $\bS_{n+1}\succeq0$ for all $n$.

Since $\range(\bM\bU_n)\subseteq\range(\widehat{\bU}_{n+1})$ and $\rank(\bP_n^{\mathrm{M}})\leq r$, the optimality of the rank-$r$ spectral truncation gives
\begin{equation*}
\operatorname{tr}(\bP_{n+1}\bS_{n+1}\bP_{n+1})
\geq\operatorname{tr}(\bP_n^{\mathrm{M}}\bS_{n+1}\bP_n^{\mathrm{M}}).
\end{equation*}
Consequently, since $\range(\bM\bC_n^{\mathrm{BUG}}\bM^\top)\subseteq\range(\bM\bU_n)$, we obtain
\begin{align*}
\operatorname{tr}(\bS_{n+1}-\bP_{n+1}\bS_{n+1}\bP_{n+1})
&\leq\operatorname{tr}(\bS_{n+1}-\bP_n^{\mathrm M}\bS_{n+1}\bP_n^{\mathrm M})\\
&=\Delta t\,\operatorname{tr}((\bI-\bP_n^{\mathrm M})\bN_n)\\
&\leq\Delta t\,\operatorname{tr}(\bN_n),
\end{align*}
where the last inequality follows from $\bN_n\succeq0$.

Since $\bC_n^{\mathrm{BUG}}\succeq0$ and $\bB_\ell^\top\bB_\ell\preceq\|\bB_\ell\|_2^2\bI$, the preservation of the Löwner order under congruence transformations gives
\begin{equation*}
\operatorname{tr}(\bN_n)\leq\|\bB_0\|_{\mathrm F}^2+\sum_{\ell=1}^q\|\bB_\ell\|_2^2\operatorname{tr}(\bC_n^{\mathrm{BUG}}).
\end{equation*}
Moreover, using the second-moment update, $\bM=\bI+\Delta t\bA$, and the same congruence argument, we obtain
\begin{equation*}
\operatorname{tr}(\bC_{n+1}^{\mathrm{BUG}})\leq(1+c\Delta t)\operatorname{tr}(\bC_n^{\mathrm{BUG}})+\Delta t\|\bB_0\|_{\mathrm F}^2,
\end{equation*}
where
\begin{equation*}
c=2\|\bA\|_2+T\|\bA\|_2^2+\sum_{\ell=1}^q\|\bB_\ell\|_2^2.
\end{equation*}
By the discrete Grönwall inequality, it follows that
\begin{equation*}
\operatorname{tr}(\bC_n^{\mathrm{BUG}})\leq e^{cT}\left(\operatorname{tr}(\bC_0^{\mathrm{BUG}})+T\|\bB_0\|_{\mathrm F}^2\right).
\end{equation*}
Summing over $n$ and taking expectations yield
\begin{equation*}
\gamma_r^2\leq C_{T,r,N}^{\mathrm{bound}},
\end{equation*}
where
\begin{equation*}
C_{T,r,N}^{\mathrm{bound}}=T\left(\|\bB_0\|_{\mathrm F}^2+\sum_{\ell=1}^q\|\bB_\ell\|_2^2e^{cT}\left(\E\left[\operatorname{tr}(\bC_0^{\mathrm{BUG}})\right]+T\|\bB_0\|_{\mathrm F}^2\right)\right).
\end{equation*}

Finally, if $\bX_0^{\mathrm{BUG}}$ is obtained from the rank-$r$ truncated SVD of $\bX_0^{\mathrm{EM}}$ and $\bD_0=\bY_0^\top\bY_0/N$, then
\begin{equation*}
\operatorname{tr}(\bC_0^{\mathrm{BUG}})=\|\bX_0^{\mathrm{BUG}}\|_N^2\leq\|\bX_0^{\mathrm{EM}}\|_N^2.
\end{equation*}
Moreover, if $\E[\bC_{N,0}^{\mathrm{EM}}]=\bC_0^{\mathrm{EM}}$, then taking expectations gives
\begin{equation*}
\begin{aligned}
\E\left[\operatorname{tr}(\bC_0^{\mathrm{BUG}})\right]
&\leq\E\left[\|\bX_0^{\mathrm{EM}}\|_N^2\right]\\
&=\operatorname{tr}\left(\E\left[\bC_{N,0}^{\mathrm{EM}}\right]\right)\\
&=\operatorname{tr}(\bC_0^{\mathrm{EM}})\\
&=\E\left[\|X(0)\|_2^2\right].
\end{aligned}
\end{equation*}
Consequently, we obtain
\begin{equation*}
\gamma_r^2\leq C_T^{\mathrm{bound}},
\end{equation*}
where
\begin{equation*}
C_T^{\mathrm{bound}}=T\left(\|\bB_0\|_{\mathrm F}^2+\sum_{\ell=1}^q\|\bB_\ell\|_2^2e^{cT}\left(\E\left[\|X(0)\|_2^2\right]+T\|\bB_0\|_{\mathrm F}^2\right)\right).
\end{equation*}
\end{proof}

\begin{proposition}[Strong comparison]
\label{prop:analysis-strong-closed-em}
Consider the closed-moment EM--BUG integrator with the reduced second moment initialized as $\bD_0=\bY_0^\top\bY_0/N$. Assume that the closed-moment EM--BUG integrator and the Euler--Maruyama scheme use the same Wiener increments. Then, for additive noise ($q=0$), the closed-moment EM--BUG integrator satisfies
\begin{equation}
\label{eq:analysis-strong-additive}
\left(\E\left[\max_{0\le n\le N_t}\|\bX_n^{\mathrm{BUG}}-\bX_n^{\mathrm{EM}}\|_N^2\right]\right)^{1/2}\leq C^{\mathrm{comp,add}}_T(\delta_{r,N}+\gamma_r),
\end{equation}
while for general additive and multiplicative noise, it satisfies
\begin{equation}
\label{eq:analysis-strong-multiplicative}
\max_{0\le n\le N_t}\left(\E\left[\|\bX_n^{\mathrm{BUG}}-\bX_n^{\mathrm{EM}}\|_N^2\right]\right)^{1/2}\leq C^{\mathrm{comp}}_T(\delta_{r,N}+\gamma_r),
\end{equation}
where the constants $C^{\mathrm{comp,add}}_T$ and $C^{\mathrm{comp}}_T$ are independent of the rank $r$, the number of particles $N$, the time-step size $\Delta t$, the smallest eigenvalue of the reduced second moment, and spectral gaps.
\end{proposition}

\begin{proof}
Let $\bX_n^{\mathrm{aux}}$ be the auxiliary ensemble introduced in Lemma~\ref{lem:analysis-strong-closed-aux}, and let $\widehat{\bP}_{n+1}=\widehat{\bU}_{n+1}\widehat{\bU}_{n+1}^\top$ denote the orthogonal projector onto the augmented physical basis. Since this basis contains the ranges of $\bM\bU_n$ and $\{\bB_\ell\bU_n\}_{\ell=1}^q$, we have
\begin{equation*}
\begin{aligned}
\operatorname{tr}\left(\bS_{n+1}-\widehat{\bP}_{n+1}\bS_{n+1}\widehat{\bP}_{n+1}\right)
&=\Delta t\,\operatorname{tr}\left(\bB_0\bB_0^\top-\widehat{\bP}_{n+1}\bB_0\bB_0^\top\widehat{\bP}_{n+1}\right)\\
&=\Delta t\,\|(\bI-\widehat{\bP}_{n+1})\bB_0\|_{\mathrm F}^2.
\end{aligned}
\end{equation*}
Consequently, using the spectral truncation defect $\tau_{n+1}$ defined in \eqref{eq:analysis-strong-closed-tau}, we obtain
\begin{equation}
\label{eq:analysis-strong-loss-decomposition}
\gamma_r^2=\sum_{n=0}^{N_t-1}\E\left[\tau_{n+1}+\Delta t\,\|(\bI-\widehat{\bP}_{n+1})\bB_0\|_{\mathrm F}^2\right].
\end{equation}
Let $c_1=12\exp(3T^2\|\bA\|_2^2)$ and $c_2=32\exp(2T^2\|\bA\|_2^2)$. For additive noise ($q=0$), the estimates \eqref{eq:analysis-strong-closed-max}, \eqref{eq:analysis-strong-closed-aux-em}, and \eqref{eq:analysis-strong-loss-decomposition} give
\begin{equation*}
\begin{aligned}
&\E\left[\max_{0\le n\le N_t}\|\bX_n^{\mathrm{BUG}}-\bX_n^{\mathrm{EM}}\|_N^2\right]\\
&\qquad\leq 2\E\left[\max_{0\le n\le N_t}\|\bX_n^{\mathrm{BUG}}-\bX_n^{\mathrm{aux}}\|_N^2\right]+2\E\left[\max_{0\le n\le N_t}\|\bX_n^{\mathrm{aux}}-\bX_n^{\mathrm{EM}}\|_N^2\right]\\
&\qquad\leq 2c_2\sum_{n=0}^{N_t-1}\E[\tau_{n+1}]+2c_1\left(\delta_{r,N}^2+\Delta t\sum_{n=0}^{N_t-1}\E\left[\|(\bI-\widehat{\bP}_{n+1})\bB_0\|_{\mathrm F}^2\right]\right)\\
&\qquad\leq 2\max\{c_1,c_2\}(\delta_{r,N}^2+\gamma_r^2).
\end{aligned}
\end{equation*}
Taking the square root and using $(a^2+b^2)^{1/2}\leq a+b$ for $a,b\geq0$, we obtain
\begin{equation*}
\begin{aligned}
\left(\E\left[\max_{0\le n\le N_t}\|\bX_n^{\mathrm{BUG}}-\bX_n^{\mathrm{EM}}\|_N^2\right]\right)^{1/2}
&\leq\sqrt{2\max\{c_1,c_2\}}\,(\delta_{r,N}^2+\gamma_r^2)^{1/2}\\
&\leq C_T^{\mathrm{comp,add}}(\delta_{r,N}+\gamma_r),
\end{aligned}
\end{equation*}
where
\begin{equation*}
C_T^{\mathrm{comp,add}}=\sqrt{2\max\{c_1,c_2\}}.
\end{equation*}

For general additive and multiplicative noise, the triangle inequality and $\max_n\E[\|\cdot\|_N^2]\leq\E[\max_n\|\cdot\|_N^2]$ give
\begin{equation*}
\begin{aligned}
&\max_{0\le n\le N_t}\E\left[\|\bX_n^{\mathrm{BUG}}-\bX_n^{\mathrm{EM}}\|_N^2\right]\\
&\qquad\leq 2\max_{0\le n\le N_t}\E\left[\|\bX_n^{\mathrm{BUG}}-\bX_n^{\mathrm{aux}}\|_N^2\right]+2\max_{0\le n\le N_t}\E\left[\|\bX_n^{\mathrm{aux}}-\bX_n^{\mathrm{EM}}\|_N^2\right]\\
&\qquad\leq 2\max_{0\le n\le N_t}\E\left[\|\bX_n^{\mathrm{BUG}}-\bX_n^{\mathrm{aux}}\|_N^2\right]+2\E\left[\max_{0\le n\le N_t}\|\bX_n^{\mathrm{aux}}-\bX_n^{\mathrm{EM}}\|_N^2\right].
\end{aligned}
\end{equation*}
Let $c_3=12\exp(T(3T\|\bA\|_2^2+12\sum_{\ell=1}^q\|\bB_\ell\|_2^2))$ and $c_4=\exp(T(2\|\bA\|_2+T\|\bA\|_2^2+\sum_{\ell=1}^q\|\bB_\ell\|_2^2))$. The estimates \eqref{eq:analysis-strong-closed-aux-lr}, \eqref{eq:analysis-strong-closed-aux-em}, and \eqref{eq:analysis-strong-loss-decomposition} then yield
\begin{equation*}
\begin{aligned}
&\max_{0\le n\le N_t}\E\left[\|\bX_n^{\mathrm{BUG}}-\bX_n^{\mathrm{EM}}\|_N^2\right]\\
&\qquad\leq 2c_4\sum_{n=0}^{N_t-1}\E[\tau_{n+1}]+2c_3\left(\delta_{r,N}^2+\Delta t\sum_{n=0}^{N_t-1}\E\left[\|(\bI-\widehat{\bP}_{n+1})\bB_0\|_{\mathrm F}^2\right]\right)\\
&\qquad\leq 2\max\{c_3,c_4\}(\delta_{r,N}^2+\gamma_r^2).
\end{aligned}
\end{equation*}
Taking the square root, we obtain
\begin{equation*}
\begin{aligned}
\max_{0\le n\le N_t}\left(\E\left[\|\bX_n^{\mathrm{BUG}}-\bX_n^{\mathrm{EM}}\|_N^2\right]\right)^{1/2}
&\leq\sqrt{2\max\{c_3,c_4\}}\,(\delta_{r,N}^2+\gamma_r^2)^{1/2}\\
&\leq C_T^{\mathrm{comp}}(\delta_{r,N}+\gamma_r),
\end{aligned}
\end{equation*}
where
\begin{equation*}
C_T^{\mathrm{comp}}=\sqrt{2\max\{c_3,c_4\}}.
\end{equation*}
\end{proof}

\begin{theorem}[Strong convergence]
\label{thm:analysis-strong-closed}
Consider the closed-moment EM--BUG integrator with the reduced second moment initialized as $\bD_0=\bY_0^\top\bY_0/N$. Assume that the hypotheses of the Euler--Maruyama strong convergence estimate \eqref{eq:analysis-strong-em} hold, and that the Wiener increments used by the closed-moment EM--BUG integrator and the Euler--Maruyama scheme are those of the Wiener processes driving the exact solution. Then, for additive noise ($q=0$), the closed-moment EM--BUG integrator satisfies
\begin{equation}
\label{eq:analysis-strong-closed-additive-exact}
\left(\E\left[\max_{0\le n\le N_t}\|\bX_n^{\mathrm{BUG}}-\bX(t_n)\|_N^2\right]\right)^{1/2}\leq C^{\mathrm{strong,add}}_T(\delta_{r,N}+\gamma_r+\Delta t),
\end{equation}
while for general additive and multiplicative noise, it satisfies
\begin{equation}
\label{eq:analysis-strong-closed-multiplicative-exact}
\max_{0\le n\le N_t}\left(\E\left[\|\bX_n^{\mathrm{BUG}}-\bX(t_n)\|_N^2\right]\right)^{1/2}\leq C^{\mathrm{strong}}_T(\delta_{r,N}+\gamma_r+(\Delta t)^{1/2}),
\end{equation}
where the constants $C^{\mathrm{strong,add}}_T$ and $C^{\mathrm{strong}}_T$ are independent of the rank $r$, the number of particles $N$, the time-step size $\Delta t$, the smallest eigenvalue of the reduced second moment, and spectral gaps.
\end{theorem}

\begin{proof}
For additive noise, Proposition~\ref{prop:analysis-strong-closed-em} and the strong Euler--Maruyama estimate \eqref{eq:analysis-strong-em} with $\alpha=1$ yield
\begin{equation*}
\begin{aligned}
&\left(\E\left[\max_{0\le n\le N_t}\|\bX_n^{\mathrm{BUG}}-\bX(t_n)\|_N^2\right]\right)^{1/2}\\
&\qquad\leq\left(\E\left[\max_{0\le n\le N_t}\|\bX_n^{\mathrm{BUG}}-\bX_n^{\mathrm{EM}}\|_N^2\right]\right)^{1/2}+\left(\E\left[\max_{0\le n\le N_t}\|\bX_n^{\mathrm{EM}}-\bX(t_n)\|_N^2\right]\right)^{1/2}\\
&\qquad\leq C^{\mathrm{comp,add}}_T(\delta_{r,N}+\gamma_r)+C^{\mathrm{EM}}_T\Delta t\\
&\qquad\leq C^{\mathrm{strong,add}}_T(\delta_{r,N}+\gamma_r+\Delta t).
\end{aligned}
\end{equation*}
where
\begin{equation*}
C^{\mathrm{strong,add}}_T=\max\{C^{\mathrm{comp,add}}_T,C^{\mathrm{EM}}_T\}.
\end{equation*}

For general additive and multiplicative noise, the triangle inequality and $\max_n\E[\|\cdot\|_N^2]\leq\E[\max_n\|\cdot\|_N^2]$ give
\begin{equation*}
\begin{aligned}
&\max_{0\le n\le N_t}\left(\E\left[\|\bX_n^{\mathrm{BUG}}-\bX(t_n)\|_N^2\right]\right)^{1/2}\\
&\qquad\leq \max_{0\le n\le N_t}\left(\E\left[\|\bX_n^{\mathrm{BUG}}-\bX_n^{\mathrm{EM}}\|_N^2\right]\right)^{1/2}
+\max_{0\le n\le N_t}\left(\E\left[\|\bX_n^{\mathrm{EM}}-\bX(t_n)\|_N^2\right]\right)^{1/2}\\
&\qquad\leq \max_{0\le n\le N_t}\left(\E\left[\|\bX_n^{\mathrm{BUG}}-\bX_n^{\mathrm{EM}}\|_N^2\right]\right)^{1/2}
+\left(\E\left[\max_{0\le n\le N_t}\|\bX_n^{\mathrm{EM}}-\bX(t_n)\|_N^2\right]\right)^{1/2}.
\end{aligned}
\end{equation*}
Proposition~\ref{prop:analysis-strong-closed-em} and the strong Euler--Maruyama estimate \eqref{eq:analysis-strong-em} with $\alpha=1/2$ then yield
\begin{equation*}
\begin{aligned}
\max_{0\le n\le N_t}\left(\E\left[\|\bX_n^{\mathrm{BUG}}-\bX(t_n)\|_N^2\right]\right)^{1/2}&\leq C_T^{\mathrm{comp}}(\delta_{r,N}+\gamma_r)+C_T^{\mathrm{EM}}(\Delta t)^{1/2}\\
&\leq C_T^{\mathrm{strong}}(\delta_{r,N}+\gamma_r+(\Delta t)^{1/2}),
\end{aligned}
\end{equation*}
where
\begin{equation*}
C_T^{\mathrm{strong}}=\max\{C_T^{\mathrm{comp}},C_T^{\mathrm{EM}}\}.
\end{equation*}
\end{proof}

\subsection{Weak convergence}
\label{sec:analysis-weak}

We now establish the weak convergence estimate of the closed-moment EM--BUG integrator. For an observable $\varphi:\R^d\to\R$, we introduce
\begin{equation*}
\Phi_N(\bX):=\frac{1}{N}\sum_{j=1}^{N}\varphi(X_j).
\end{equation*}
We recall that, under suitable assumptions on the coefficients and the observable, the Euler--Maruyama scheme satisfies
\begin{equation}
\label{eq:analysis-weak-em}
\max_{0\le n\le N_t}\left|\E[\Phi_N(\bX_n^{\mathrm{EM}})]-\E[\varphi(X(t_n))]\right|\le C_{T,\varphi}^{\mathrm{EM}}(\Delta t)^\beta,
\end{equation}
where $C_{T,\varphi}^{\mathrm{EM}}$ is independent of $N$ and $\Delta t$. For globally Lipschitz observables, $\beta=1/2$ follows from the standard strong convergence assumptions. The order $\beta=1$ can be obtained in specific cases, for example under additional regularity assumptions on the coefficients and the observable, or for linear SDEs with constant coefficients and additive noise.

\begin{theorem}[Weak convergence]
\label{thm:analysis-weak}
Let $\varphi:\R^d\to\R$ be globally Lipschitz continuous, with Lipschitz constant $L_\varphi$. Consider the closed-moment EM--BUG integrator with the reduced second moment initialized as $\bD_0=\bY_0^\top\bY_0/N$. Assume that the hypotheses of Theorem~\ref{thm:analysis-strong-closed} and of the Euler--Maruyama weak convergence estimate \eqref{eq:analysis-weak-em} hold. Then, the closed-moment EM--BUG integrator satisfies
\begin{equation}
\label{eq:analysis-weak-closed}
\max_{0\le n\le N_t}\left|\E[\Phi_N(\bX_n^{\mathrm{BUG}})]-\E[\varphi(X(t_n))]\right|\le C^{\mathrm{weak}}_{T,\varphi}\left(\delta_{r,N}+\gamma_r+(\Delta t)^\beta\right),
\end{equation}
where the constant  $C^{\mathrm{weak}}_{T,\varphi}$ is independent of the rank $r$, the number of particles $N$, the time-step size $\Delta t$, the smallest eigenvalue of the reduced second moment, and spectral gaps.
\end{theorem}

\begin{proof}
The Lipschitz continuity of $\varphi$ and the Cauchy--Schwarz inequality give
\begin{equation}
\label{eq:analysis-weak-lipschitz}
\begin{aligned}
\left|\Phi_N(\bX)-\Phi_N(\bY)\right|&\leq\frac{L_\varphi}{N}\sum_{j=1}^N\|X_j-Y_j\|\\
&\leq L_\varphi\|\bX-\bY\|_N.
\end{aligned}
\end{equation}
Consequently, using the triangle inequality and the Cauchy--Schwarz inequality, we obtain
\begin{equation*}
\begin{aligned}
&\left|\E[\Phi_N(\bX_n^{\mathrm{BUG}})]-\E[\varphi(X(t_n))]\right|\\
&\qquad\leq\left|\E[\Phi_N(\bX_n^{\mathrm{BUG}})]-\E[\Phi_N(\bX_n^{\mathrm{EM}})]\right|+\left|\E[\Phi_N(\bX_n^{\mathrm{EM}})]-\E[\varphi(X(t_n))]\right|\\
&\qquad\leq\E\left[\left|\Phi_N(\bX_n^{\mathrm{BUG}})-\Phi_N(\bX_n^{\mathrm{EM}})\right|\right]+C_{T,\varphi}^{\mathrm{EM}}(\Delta t)^\beta\\
&\qquad\leq L_\varphi\left(\E\left[\|\bX_n^{\mathrm{BUG}}-\bX_n^{\mathrm{EM}}\|_N^2\right]\right)^{1/2}+C_{T,\varphi}^{\mathrm{EM}}(\Delta t)^\beta.
\end{aligned}
\end{equation*}
Taking the maximum over $n$ and using the strong comparison estimate \eqref{eq:analysis-strong-multiplicative} then yield
\begin{equation*}
\begin{aligned}
&\max_{0\le n\le N_t}\left|\E[\Phi_N(\bX_n^{\mathrm{BUG}})]-\E[\varphi(X(t_n))]\right|\\
&\qquad\leq L_\varphi\max_{0\le n\le N_t}\left(\E\left[\|\bX_n^{\mathrm{BUG}}-\bX_n^{\mathrm{EM}}\|_N^2\right]\right)^{1/2}+C_{T,\varphi}^{\mathrm{EM}}(\Delta t)^\beta\\
&\qquad\leq L_\varphi C_T^{\mathrm{comp}}(\delta_{r,N}+\gamma_r)+C_{T,\varphi}^{\mathrm{EM}}(\Delta t)^\beta\\
&\qquad\leq C_{T,\varphi}^{\mathrm{weak}}(\delta_{r,N}+\gamma_r+(\Delta t)^\beta),
\end{aligned}
\end{equation*}
where
\begin{equation*}
C_{T,\varphi}^{\mathrm{weak}}=\max\{L_\varphi C_T^{\mathrm{comp}},C_{T,\varphi}^{\mathrm{EM}}\}.
\end{equation*}
\end{proof}

\begin{proposition}[Weak convergence for linear observables]
\label{prop:analysis-weak-linear}
Consider the linear SDE \eqref{eq:linear-sde} with constant coefficients, and let $\varphi_v(x)=v^\top x$ with $v\in\R^d$. Consider the closed-moment EM--BUG integrator with the reduced second moment initialized as $\bD_0=\bY_0^\top\bY_0/N$. Assume that the hypotheses of Theorem~\ref{thm:analysis-weak} hold. Then, the closed-moment EM--BUG integrator satisfies
\begin{equation}
\label{eq:analysis-weak-linear-closed}
\max_{0\le n\le N_t}\left|\E[\Phi_N(\bX_n^{\mathrm{BUG}})]-\E[\varphi_v(X(t_n))]\right|\leq C^{\mathrm{weak}}_{T,v}\left(\delta_{r,N}+\gamma_r+\Delta t\right),
\end{equation}
where the constant $C^{\mathrm{weak}}_{T,v}$ is independent of the rank $r$, the number of particles $N$, the time-step size $\Delta t$, the smallest eigenvalue of the reduced second moments, and spectral gaps.
\end{proposition}

\begin{proof}
For linear SDEs with constant coefficients, the exact mean $m(t)=\E[X(t)]$ and its Euler--Maruyama discretization satisfy
\begin{equation*}
\dot m(t)=\bA m(t),\qquad m_{n+1}^{\mathrm{EM}}=(\bI+\Delta t\,\bA)m_n^{\mathrm{EM}}.
\end{equation*}
This corresponds to the explicit Euler discretization, which converges with order one:
\begin{equation*}
\max_{0\le n\le N_t}\|m_n^{\mathrm{EM}}-m(t_n)\|\leq C_T^{\mathrm{EM,mean}}\Delta t.
\end{equation*}
Since $\varphi_v(x)=v^\top x$ is linear, we obtain
\begin{equation*}
\begin{aligned}
\max_{0\le n\le N_t}\left|\E[\Phi_N(\bX_n^{\mathrm{EM}})]-\E[\varphi_v(X(t_n))]\right|
&=\max_{0\le n\le N_t}|v^\top(m_n^{\mathrm{EM}}-m(t_n))|\\
&\leq\|v\|\max_{0\le n\le N_t}\|m_n^{\mathrm{EM}}-m(t_n)\|\\
&\leq C_T^{\mathrm{EM,mean}}\|v\|\Delta t.
\end{aligned}
\end{equation*}
The Lipschitz continuity of $\varphi_v$ with $L_{\varphi_v}=\|v\|$ and the strong comparison estimate \eqref{eq:analysis-strong-multiplicative} then yield
\begin{equation*}
\begin{aligned}
&\max_{0\le n\le N_t}\left|\E[\Phi_N(\bX_n^{\mathrm{BUG}})]-\E[\varphi_v(X(t_n))]\right|\\
&\qquad\leq\|v\|\max_{0\le n\le N_t}\left(\E\left[\|\bX_n^{\mathrm{BUG}}-\bX_n^{\mathrm{EM}}\|_N^2\right]\right)^{1/2}+C_T^{\mathrm{EM,mean}}\|v\|\Delta t\\
&\qquad\leq\|v\|C_T^{\mathrm{comp}}(\delta_{r,N}+\gamma_r)+C_T^{\mathrm{EM,mean}}\|v\|\Delta t\\
&\qquad\leq C_{T,v}^{\mathrm{weak}}(\delta_{r,N}+\gamma_r+\Delta t),
\end{aligned}
\end{equation*}
where
\begin{equation*}
C_{T,v}^{\mathrm{weak}}=\|v\|\max\{C_T^{\mathrm{comp}},C_T^{\mathrm{EM,mean}}\}.
\end{equation*}
\end{proof}

\begin{proposition}[Weak convergence for quadratic observables]
\label{prop:analysis-weak-quadratic}
Consider the linear SDE \eqref{eq:linear-sde} with constant coefficients, and let $\varphi_H(x)=x^\top\bH x$, where $\bH\in\R^{d\times d}$ is symmetric. Consider the closed-moment EM--BUG integrator with complete augmentation and with the reduced second moment initialized as $\bD_0=\bY_0^\top\bY_0/N$, where the initial ensemble $\bX_0^{\mathrm{BUG}}$ is obtained from the rank-$r$ truncated SVD of $\bX_0^{\mathrm{EM}}$. Assume that $\E[\bC_{N,0}^{\mathrm{EM}}]=\bC_0^{\mathrm{EM}}$. Then, the closed-moment EM--BUG integrator satisfies
\begin{equation}
\label{eq:analysis-weak-quadratic}
\max_{0\le n\le N_t}\left|\E[\Phi_N(\bX_n^{\mathrm{BUG}})]-\E[\varphi_H(X(t_n))]\right|\leq C^{\mathrm{weak}}_{T,H}\left(\delta_{r,N}^2+\gamma_r^2+\Delta t\right),
\end{equation}
where the constant $C^{\mathrm{weak}}_{T,H}$ is independent of the rank $r$, the number of particles $N$, the time-step size $\Delta t$, the smallest eigenvalue of the reduced second moments, and spectral gaps.
\end{proposition}

\begin{proof}
Let $\mathcal G=\sigma(\bX_0^{\mathrm{EM}})$ denote the sigma-algebra generated by the initial Euler--Maruyama ensemble. Since the initial low-rank ensemble $\bX_0^{\mathrm{BUG}}$ is obtained from $\bX_0^{\mathrm{EM}}$ and the closed second-moment recursion is deterministic, the projector $\bP_{n+1}$ is $\mathcal G$-measurable. Given the initialization $\bD_0=\bY_0^\top\bY_0/N$, the conditional second-moment recursion yields $\E[\bC_{N,n}^{\mathrm{BUG}}\mid\mathcal G]=\bC_n^{\mathrm{BUG}}$ for all $n$. Taking expectations and using \eqref{eq:analysis-moment-structure}, we obtain
\begin{equation*}
\E[\bC_{N,n+1}^{\mathrm{BUG}}]=\E[\bP_{n+1}\bS_{n+1}\bP_{n+1}].
\end{equation*}
Since
\begin{equation*}
\begin{aligned}
\bC_{n+1}^{\mathrm{EM}}&=\bM\bC_n^{\mathrm{EM}}\bM^\top+\Delta t\,\bB_0\bB_0^\top+\Delta t\sum_{\ell=1}^q\bB_\ell\bC_n^{\mathrm{EM}}\bB_\ell^\top,\\
\E[\bS_{n+1}]&=\bM\E[\bC_{N,n}^{\mathrm{BUG}}]\bM^\top+\Delta t\,\bB_0\bB_0^\top+\Delta t\sum_{\ell=1}^q\bB_\ell\E[\bC_{N,n}^{\mathrm{BUG}}]\bB_\ell^\top,
\end{aligned}
\end{equation*}
the second-moment error $\bE_n=\bC_n^{\mathrm{EM}}-\E[\bC_{N,n}^{\mathrm{BUG}}]$ satisfies
\begin{equation}
\label{eq:analysis-weak-quadratic-recursion}
\bE_{n+1}=\bM\bE_n\bM^\top+\Delta t\sum_{\ell=1}^q\bB_\ell\bE_n\bB_\ell^\top+\E[\bL_{n+1}],
\end{equation}
where $\bL_{n+1}=\bS_{n+1}-\bP_{n+1}\bS_{n+1}\bP_{n+1}$.

Since the augmentation is complete, the augmented basis $\widehat{\bU}_{n+1}$ contains the range of $\bS_{n+1}$. Hence, the columns of $\bU_{n+1}$ selected by spectral truncation are eigenvectors of $\bS_{n+1}$, and the corresponding projector $\bP_{n+1}=\bU_{n+1}\bU_{n+1}^\top$ satisfies
\begin{equation*}
\bP_{n+1}\bS_{n+1}=\bP_{n+1}\bS_{n+1}\bP_{n+1}=\bS_{n+1}\bP_{n+1}.
\end{equation*}
Consequently,
\begin{equation*}
\bL_{n+1}=(\bI-\bP_{n+1})\bS_{n+1}(\bI-\bP_{n+1})\succeq0.
\end{equation*}
Moreover, since $\bX_0^{\mathrm{BUG}}$ is obtained by a truncated SVD of $\bX_0^{\mathrm{EM}}$, the retained and discarded components are orthogonal, and we obtain
\begin{equation*}
\bX_0^{\mathrm{EM}}(\bX_0^{\mathrm{EM}})^\top-\bX_0^{\mathrm{BUG}}(\bX_0^{\mathrm{BUG}})^\top=(\bX_0^{\mathrm{EM}}-\bX_0^{\mathrm{BUG}})(\bX_0^{\mathrm{EM}}-\bX_0^{\mathrm{BUG}})^\top.
\end{equation*}
It follows that
\begin{equation*}
\bE_0=\E\left[\frac1N(\bX_0^{\mathrm{EM}}-\bX_0^{\mathrm{BUG}})(\bX_0^{\mathrm{EM}}-\bX_0^{\mathrm{BUG}})^\top\right]\succeq0,
\qquad
\operatorname{tr}(\bE_0)=\delta_{r,N}^2.
\end{equation*}
Since each term on the right-hand side of \eqref{eq:analysis-weak-quadratic-recursion} preserves positive semidefiniteness, we deduce that $\bE_n\succeq0$ for all $n$.

Taking the trace of \eqref{eq:analysis-weak-quadratic-recursion} gives
\begin{equation*}
\operatorname{tr}(\bE_{n+1})=\operatorname{tr}\left[\left(\bM^\top\bM+\Delta t\sum_{\ell=1}^q\bB_\ell^\top\bB_\ell\right)\bE_n\right]+\E[\operatorname{tr}(\bL_{n+1})].
\end{equation*}
Moreover, using $\Delta t\leq T$ yields
\begin{equation*}
\begin{aligned}
\bM^\top\bM+\Delta t\sum_{\ell=1}^q\bB_\ell^\top\bB_\ell
&=\bI+\Delta t\left(\bA+\bA^\top+\sum_{\ell=1}^q\bB_\ell^\top\bB_\ell\right)+(\Delta t)^2\bA^\top\bA\\
&\preceq(1+c\Delta t)\bI,
\end{aligned}
\end{equation*}
where $c=2\|\bA\|_2+\sum_{\ell=1}^q\|\bB_\ell\|_2^2+T\|\bA\|_2^2$. Since $\bE_n\succeq0$, the preservation of the Löwner order under congruence transformations gives
\begin{equation*}
\operatorname{tr}(\bE_{n+1})\leq(1+c\Delta t)\operatorname{tr}(\bE_n)+\E[\operatorname{tr}(\bL_{n+1})].
\end{equation*}
Consequently, by the discrete Gr\"onwall inequality, we obtain
\begin{equation}
\label{eq:analysis-weak-quadratic-moment-error}
\begin{aligned}
\operatorname{tr}(\bE_n)
&\leq(1+c\Delta t)^n\operatorname{tr}(\bE_0)+\sum_{k=0}^{n-1}(1+c\Delta t)^{n-1-k}\E[\operatorname{tr}(\bL_{k+1})]\\
&\leq e^{cT}\left(\operatorname{tr}(\bE_0)+\sum_{k=0}^{n-1}\E[\operatorname{tr}(\bL_{k+1})]\right)\\
&\leq e^{cT}(\delta_{r,N}^2+\gamma_r^2).
\end{aligned}
\end{equation}

The Euler--Maruyama discretization of the closed second-moment equation converges with order one, since it differs from the explicit Euler discretization only by a term of order $(\Delta t)^2$:
\begin{equation}
\label{eq:analysis-weak-quadratic-moment-em}
\max_{0\le n\le N_t}\|\bC_n^{\mathrm{EM}}-\bC(t_n)\|_{\mathrm F}\leq C_T^{\mathrm{EM,mom}}\Delta t.
\end{equation}
Using $\varphi_H(x)=x^\top\bH x$ and $x^\top\bH x=\operatorname{tr}(\bH xx^\top)$ yields
\begin{equation*}
\begin{aligned}
\left|\E[\Phi_N(\bX_n^{\mathrm{BUG}})]-\E[\varphi_H(X(t_n))]\right|
&=\left|\operatorname{tr}\left(\bH(\E[\bC_{N,n}^{\mathrm{BUG}}]-\bC(t_n))\right)\right|\\
&\leq\left|\operatorname{tr}(\bH\bE_n)\right|+\left|\operatorname{tr}\left(\bH(\bC_n^{\mathrm{EM}}-\bC(t_n))\right)\right|\\
&\leq\|\bH\|_2\operatorname{tr}(\bE_n)+\left|\left\langle\bH,\bC_n^{\mathrm{EM}}-\bC(t_n)\right\rangle_{\mathrm F}\right|\\
&\leq\|\bH\|_2\operatorname{tr}(\bE_n)+\|\bH\|_{\mathrm F}\|\bC_n^{\mathrm{EM}}-\bC(t_n)\|_{\mathrm F},
\end{aligned}
\end{equation*}
where $|\operatorname{tr}(\bH\bE_n)|\leq\|\bH\|_2\operatorname{tr}(\bE_n)$ follows from $\bE_n\succeq0$, $-\|\bH\|_2\bI\preceq\bH\preceq\|\bH\|_2\bI$ (since $\bH$ is symmetric), and the preservation of the Löwner order under congruence transformations. Taking the maximum over $n$ and using the estimates \eqref{eq:analysis-weak-quadratic-moment-error} and \eqref{eq:analysis-weak-quadratic-moment-em}, we obtain
\begin{equation*}
\begin{aligned}
\max_{0\le n\le N_t}\left|\E[\Phi_N(\bX_n^{\mathrm{BUG}})]-\E[\varphi_H(X(t_n))]\right|
&\leq e^{cT}\|\bH\|_2(\delta_{r,N}^2+\gamma_r^2)+C_T^{\mathrm{EM,mom}}\|\bH\|_{\mathrm F}\Delta t\\
&\leq C_{T,H}^{\mathrm{weak}}(\delta_{r,N}^2+\gamma_r^2+\Delta t),
\end{aligned}
\end{equation*}
where
\begin{equation*}
C_{T,H}^{\mathrm{weak}}=\max\left\{e^{cT}\|\bH\|_2,C_T^{\mathrm{EM,mom}}\|\bH\|_{\mathrm F}\right\}.
\end{equation*}
\end{proof}

\subsection{Monte Carlo error}
\label{sec:analysis-monte-carlo}

We now establish estimates for the root-mean-square Monte Carlo error of the closed-moment EM--BUG integrator, defined for an observable $\varphi:\R^d\to\R$ by
\begin{equation}
\label{eq:analysis-mc-error}
\max_{0\le n\le N_t}\left(\E\left[\left|\Phi_N(\bX_n^{\mathrm{BUG}})-\E[\varphi(X(t_n))]\right|^2\right]\right)^{1/2}.
\end{equation}
To control the sampling error, we introduce the following assumption.

\begin{assumption}[Independent initial particles]
\label{ass:analysis-mc-iid}
The initial particles are independent and identically distributed with the same law as $X(0)$, and independent of the Brownian motions.
\end{assumption}

Under Assumption~\ref{ass:analysis-mc-iid}, the Euler--Maruyama particles are independent and identically distributed at each time step. For the linear SDE \eqref{eq:linear-sde}, the standard Euler--Maruyama second-moment estimate gives
\begin{equation*}
\max_{0\le n\le N_t}\E[\|X_n^{\mathrm{EM}}\|^2]\leq C_T^{\mathrm{EM,stab}},
\end{equation*}
where $C_T^{\mathrm{EM,stab}}$ is independent of $N$ and $\Delta t$. It follows that, for a globally Lipschitz observable $\varphi$,
\begin{equation}
\label{eq:analysis-mc-variance}
\begin{aligned}
\max_{0\le n\le N_t}\operatorname{Var}(\varphi(X_n^{\mathrm{EM}}))&\leq\max_{0\le n\le N_t}\E\left[|\varphi(X_n^{\mathrm{EM}})-\varphi(0)|^2\right]\\
&\leq L_\varphi^2\max_{0\le n\le N_t}\E[\|X_n^{\mathrm{EM}}\|^2]\\
&\leq L_\varphi^2C_T^{\mathrm{EM,stab}}:=\sigma_\varphi^2.
\end{aligned}
\end{equation}

\begin{theorem}[Monte Carlo error]
\label{thm:analysis-monte-carlo}
Let $\varphi:\R^d\to\R$ be globally Lipschitz continuous, with Lipschitz constant $L_\varphi$. Consider the closed-moment EM--BUG integrator with the reduced second moment initialized as $\bD_0=\bY_0^\top\bY_0/N$. Suppose that Assumption~\ref{ass:analysis-mc-iid} and the hypotheses of Theorem~\ref{thm:analysis-weak} hold. Then, the closed-moment EM--BUG integrator satisfies
\begin{equation}
\label{eq:analysis-mc}
\max_{0\le n\le N_t}\left(\E\left[\left|\Phi_N(\bX_n^{\mathrm{BUG}})-\E[\varphi(X(t_n))]\right|^2\right]\right)^{1/2}\leq C_{T,\varphi}^{\mathrm{MC}}\left(\delta_{r,N}+\gamma_r+(\Delta t)^\beta\right)+\frac{\sigma_\varphi}{\sqrt N},
\end{equation}
where the constant $C_{T,\varphi}^{\mathrm{MC}}$ is independent of the rank $r$, the number of particles $N$, the time-step size $\Delta t$, the smallest eigenvalue of the reduced second moment, and spectral gaps.
\end{theorem}

\begin{proof}
The triangle inequality in $L^2(\Omega)$ gives
\begin{equation*}
\begin{aligned}
\left(\E\left[\left|\Phi_N(\bX_n^{\mathrm{BUG}})-\E[\varphi(X(t_n))]\right|^2\right]\right)^{1/2}
&\leq\left(\E\left[\left|\Phi_N(\bX_n^{\mathrm{BUG}})-\Phi_N(\bX_n^{\mathrm{EM}})\right|^2\right]\right)^{1/2}\\
&\quad+\left(\E\left[\left|\Phi_N(\bX_n^{\mathrm{EM}})-\E[\Phi_N(\bX_n^{\mathrm{EM}})]\right|^2\right]\right)^{1/2}\\
&\quad+\left|\E[\Phi_N(\bX_n^{\mathrm{EM}})]-\E[\varphi(X(t_n))]\right|.
\end{aligned}
\end{equation*}
By the Lipschitz estimate \eqref{eq:analysis-weak-lipschitz}, we have
\begin{equation*}
\left|\Phi_N(\bX_n^{\mathrm{BUG}})-\Phi_N(\bX_n^{\mathrm{EM}})\right|\leq L_\varphi\|\bX_n^{\mathrm{BUG}}-\bX_n^{\mathrm{EM}}\|_N.
\end{equation*}
Moreover, under Assumption~\ref{ass:analysis-mc-iid}, the Euler--Maruyama particles are independent and identically distributed. Consequently, the variance bound \eqref{eq:analysis-mc-variance} yields
\begin{equation*}
\begin{aligned}
\E\left[\left|\Phi_N(\bX_n^{\mathrm{EM}})-\E[\Phi_N(\bX_n^{\mathrm{EM}})]\right|^2\right]
&=\frac{1}{N}\operatorname{Var}(\varphi(X_n^{\mathrm{EM}}))\\
&\leq\frac{\sigma_\varphi^2}{N}.
\end{aligned}
\end{equation*}
Combining these estimates with the Euler--Maruyama weak convergence estimate \eqref{eq:analysis-weak-em}, we obtain
\begin{equation*}
\begin{aligned}
&\max_{0\le n\le N_t}\left(\E\left[\left|\Phi_N(\bX_n^{\mathrm{BUG}})-\E[\varphi(X(t_n))]\right|^2\right]\right)^{1/2}\\
&\qquad\leq L_\varphi\max_{0\le n\le N_t}\left(\E\left[\|\bX_n^{\mathrm{BUG}}-\bX_n^{\mathrm{EM}}\|_N^2\right]\right)^{1/2}+C_{T,\varphi}^{\mathrm{EM}}(\Delta t)^\beta+\frac{\sigma_\varphi}{\sqrt N}.
\end{aligned}
\end{equation*}
The strong comparison estimate \eqref{eq:analysis-strong-multiplicative} then yields
\begin{equation*}
\begin{aligned}
&\max_{0\le n\le N_t}\left(\E\left[\left|\Phi_N(\bX_n^{\mathrm{BUG}})-\E[\varphi(X(t_n))]\right|^2\right]\right)^{1/2}\\
&\qquad\leq L_\varphi C_T^{\mathrm{comp}}(\delta_{r,N}+\gamma_r)+C_{T,\varphi}^{\mathrm{EM}}(\Delta t)^\beta+\frac{\sigma_\varphi}{\sqrt N}\\
&\qquad\leq C_{T,\varphi}^{\mathrm{MC}}\left(\delta_{r,N}+\gamma_r+(\Delta t)^\beta\right)+\frac{\sigma_\varphi}{\sqrt N},
\end{aligned}
\end{equation*}
where
\begin{equation*}
C_{T,\varphi}^{\mathrm{MC}}=\max\{L_\varphi C_T^{\mathrm{comp}},C_{T,\varphi}^{\mathrm{EM}}\}.
\end{equation*}
\end{proof}

\begin{proposition}[Monte Carlo error for linear observables]
\label{pro:analysis-mc-linear}
Consider the linear SDE \eqref{eq:linear-sde} with constant coefficients, and let $\varphi_v(x)=v^\top x$ with $v\in\R^d$. Consider the closed-moment EM--BUG integrator with the reduced second moment initialized as $\bD_0=\bY_0^\top\bY_0/N$. Assume that the hypotheses of Theorem~\ref{thm:analysis-monte-carlo} hold. Then, the closed-moment EM--BUG integrator satisfies
\begin{equation}
\label{eq:analysis-mc-linear-closed}
\max_{0\le n\le N_t}\left(\E\left[\left|\Phi_N(\bX_n^{\mathrm{BUG}})-\E[\varphi_v(X(t_n))]\right|^2\right]\right)^{1/2}\leq C_{T,v}^{\mathrm{MC}}\left(\delta_{r,N}+\gamma_r+\Delta t\right)+\frac{\sigma_{\varphi_v}}{\sqrt N},
\end{equation}
where the constant $C_{T,v}^{\mathrm{MC}}$ is independent of the rank $r$, the number of particles $N$, the time-step size $\Delta t$, the smallest eigenvalue of the reduced second moment, and spectral gaps.
\end{proposition}

\begin{proof}
As established in the proof of Proposition~\ref{prop:analysis-weak-linear}, the first-order accuracy of the Euler--Maruyama approximation of the mean for linear SDEs with constant coefficients implies, for $\varphi_v(x)=v^\top x$,
\begin{equation*}
\max_{0\le n\le N_t}\left|\E[\Phi_N(\bX_n^{\mathrm{EM}})]-\E[\varphi_v(X(t_n))]\right|\leq\|v\|C_T^{\mathrm{EM,mean}}\Delta t.
\end{equation*}
Then, using the triangle inequality, the Lipschitz estimate \eqref{eq:analysis-weak-lipschitz}, the variance bound \eqref{eq:analysis-mc-variance}, and the strong comparison estimate \eqref{eq:analysis-strong-multiplicative}, as in the proof of Theorem~\ref{thm:analysis-monte-carlo}, we obtain
\begin{equation*}
\begin{aligned}
&\max_{0\le n\le N_t}\left(\E\left[\left|\Phi_N(\bX_n^{\mathrm{BUG}})-\E[\varphi_v(X(t_n))]\right|^2\right]\right)^{1/2}\\
&\quad\leq\|v\|\max_{0\le n\le N_t}\left(\E\left[\|\bX_n^{\mathrm{BUG}}-\bX_n^{\mathrm{EM}}\|_N^2\right]\right)^{1/2}
+\|v\|C_T^{\mathrm{EM,mean}}\Delta t+\frac{\sigma_{\varphi_v}}{\sqrt N}\\
&\quad\leq\|v\|C_T^{\mathrm{comp}}(\delta_{r,N}+\gamma_r)+\|v\|C_T^{\mathrm{EM,mean}}\Delta t+\frac{\sigma_{\varphi_v}}{\sqrt N}\\
&\quad\leq C_{T,v}^{\mathrm{MC}}(\delta_{r,N}+\gamma_r+\Delta t)+\frac{\sigma_{\varphi_v}}{\sqrt N},
\end{aligned}
\end{equation*}
where
\begin{equation*}
C_{T,v}^{\mathrm{MC}}=\|v\|\max\{C_T^{\mathrm{comp}},C_T^{\mathrm{EM,mean}}\}.
\end{equation*}
\end{proof}

\begin{proposition}[Monte Carlo error for quadratic observables]
\label{prop:analysis-mc-quadratic}
Consider the linear SDE \eqref{eq:linear-sde} with constant coefficients, and let $\varphi_H(x)=x^\top\bH x$ with $\bH\in\R^{d\times d}$ symmetric. Consider the closed-moment EM--BUG integrator with complete augmentation and with the reduced second moment initialized as $\bD_0=\bY_0^\top\bY_0/N$, where the initial ensemble $\bX_0^{\mathrm{BUG}}$ is obtained from the rank-$r$ truncated SVD of $\bX_0^{\mathrm{EM}}$. Suppose that Assumption~\ref{ass:analysis-mc-iid} holds and that $\E[\|X(0)\|_2^4]$ is finite. Then, the closed-moment EM--BUG integrator satisfies
\begin{equation}
\label{eq:analysis-mc-quadratic}
\max_{0\le n\le N_t}\left(\E\left[\left|\Phi_N(\bX_n^{\mathrm{BUG}})-\E[\varphi_H(X(t_n))]\right|^2\right]\right)^{1/2}\leq C_{T,H}^{\mathrm{MC}}\left(\delta_{r,N}+\gamma_r+\Delta t\right)+\frac{\sigma_{T,\varphi_H}}{\sqrt N},
\end{equation}
where the constants $C_{T,H}^{\mathrm{MC}}$ and $\sigma_{T,\varphi_H}$ are independent of the rank $r$, the number of particles $N$, the time-step size $\Delta t$, the smallest eigenvalue of the reduced second moment, and spectral gaps.
\end{proposition}

\begin{proof}
Let $\bGamma_n$ denote the cumulative loss matrix introduced in Lemma~\ref{lem:analysis-mc-quadratic}. By the recursion \eqref{eq:analysis-mc-quadratic-lemma-loss}, we have
\begin{equation*}
\bGamma_{n+1}=\bM\bGamma_n\bM^\top+\Delta t\sum_{\ell=1}^q\bB_\ell\bGamma_n\bB_\ell^\top+\bL_{n+1},
\end{equation*}
with $\bGamma_n\succeq0$ and $\bL_{n+1}\succeq0$. The same trace argument as in the proof of Proposition~\ref{prop:analysis-weak-quadratic} then gives
\begin{equation*}
\E[\operatorname{tr}(\bGamma_{n+1})]\leq(1+c\Delta t)\E[\operatorname{tr}(\bGamma_n)]+\E[\operatorname{tr}(\bL_{n+1})],
\end{equation*}
where $c=2\|\bA\|_2+\sum_{\ell=1}^q\|\bB_\ell\|_2^2+T\|\bA\|_2^2$. Moreover, since $\bX_0^{\mathrm{BUG}}$ is obtained by a truncated SVD of $\bX_0^{\mathrm{EM}}$, we have
\begin{equation*}
\begin{aligned}
\bGamma_0
&=\frac1N\left(\bX_0^{\mathrm{EM}}(\bX_0^{\mathrm{EM}})^\top
-\bX_0^{\mathrm{BUG}}(\bX_0^{\mathrm{BUG}})^\top\right),\\
\E[\operatorname{tr}(\bGamma_0)]
&=\E[\|\bX_0^{\mathrm{EM}}-\bX_0^{\mathrm{BUG}}\|_N^2]
=\delta_{r,N}^2.
\end{aligned}
\end{equation*}
Consequently, the discrete Gr\"onwall inequality yields
\begin{equation*}
\begin{aligned}
\E[\operatorname{tr}(\bGamma_n)]
&\leq(1+c\Delta t)^n\delta_{r,N}^2+\sum_{k=0}^{n-1}(1+c\Delta t)^{n-1-k}\E[\operatorname{tr}(\bL_{k+1})]\\
&\leq e^{cT}(\delta_{r,N}^2+\gamma_r^2).
\end{aligned}
\end{equation*}

By the decomposition \eqref{eq:analysis-mc-quadratic-decomposition}, we have
\begin{equation*}
\bC_{N,n}^{\mathrm{BUG}}=\bC_n^{\mathrm{EM}}+\bXi_n-\bGamma_n,
\end{equation*}
with $\bC_{N,n}^{\mathrm{BUG}}\succeq0$ and $\bGamma_n\succeq0$. Taking the trace then gives
\begin{equation*}
0\leq\operatorname{tr}(\bGamma_n)\leq\operatorname{tr}(\bC_n^{\mathrm{EM}})+\operatorname{tr}(\bXi_n).
\end{equation*}
Using $\operatorname{tr}(\bC_n^{\mathrm{EM}})=\E[\|X_n^{\mathrm{EM}}\|^2]\leq C_T^{\mathrm{EM,stab}}$, the Cauchy--Schwarz inequality, and the estimate \eqref{eq:analysis-mc-quadratic-fluctuation} with $\bK=\bI$, we obtain
\begin{equation*}
\begin{aligned}
\E\left[(\operatorname{tr}(\bGamma_n))^2\right]&\leq\operatorname{tr}(\bC_n^{\mathrm{EM}})\E[\operatorname{tr}(\bGamma_n)]+\E[\operatorname{tr}(\bGamma_n)\operatorname{tr}(\bXi_n)]\\
&\leq C_T^{\mathrm{EM,stab}}e^{cT}(\delta_{r,N}^2+\gamma_r^2)+\left(\E\left[(\operatorname{tr}(\bGamma_n))^2\right]\right)^{1/2}\left(\E\left[(\operatorname{tr}(\bXi_n))^2\right]\right)^{1/2}\\
&\leq C_T^{\mathrm{EM,stab}}e^{cT}(\delta_{r,N}^2+\gamma_r^2)+\frac{\sigma_{T,I}}{\sqrt N}\left(\E\left[(\operatorname{tr}(\bGamma_n))^2\right]\right)^{1/2}.
\end{aligned}
\end{equation*}
Solving this quadratic inequality yields
\begin{equation}
\label{eq:analysis-mc-quadratic-loss-L2}
\max_{0\le n\le N_t}\left(\E\left[(\operatorname{tr}(\bGamma_n))^2\right]\right)^{1/2}\leq\sqrt{C_T^{\mathrm{EM,stab}}e^{cT}}(\delta_{r,N}+\gamma_r)+\frac{\sigma_{T,I}}{\sqrt N}.
\end{equation}

Finally, since $\varphi_H(x)=x^\top\bH x$ and $x^\top\bH x=\operatorname{tr}(\bH xx^\top)$, the decomposition \eqref{eq:analysis-mc-quadratic-decomposition} gives
\begin{equation*}
\begin{aligned}
&\left(\E\left[\left|\Phi_N(\bX_n^{\mathrm{BUG}})-\E[\varphi_H(X(t_n))]\right|^2\right]\right)^{1/2}\\
&\quad\leq
\left(\E\left[|\operatorname{tr}(\bH\bXi_n)|^2\right]\right)^{1/2}
+\left(\E\left[|\operatorname{tr}(\bH\bGamma_n)|^2\right]\right)^{1/2}
+\left|\operatorname{tr}\left(\bH(\bC_n^{\mathrm{EM}}-\bC(t_n))\right)\right|\\
&\quad\leq
\left(\E\left[|\operatorname{tr}(\bH\bXi_n)|^2\right]\right)^{1/2}
+\|\bH\|_2\left(\E\left[(\operatorname{tr}(\bGamma_n))^2\right]\right)^{1/2}
+\|\bH\|_{\mathrm F}\|\bC_n^{\mathrm{EM}}-\bC(t_n)\|_{\mathrm F},
\end{aligned}
\end{equation*}
where $|\operatorname{tr}(\bH\bGamma_n)|\leq\|\bH\|_2\operatorname{tr}(\bGamma_n)$ follows from $\bGamma_n\succeq0$, $-\|\bH\|_2\bI\preceq\bH\preceq\|\bH\|_2\bI$ (since $\bH$ is symmetric), and the preservation of the Löwner order under congruence transformations. Taking the maximum over $n$ and using the estimates \eqref{eq:analysis-mc-quadratic-fluctuation} with $\bK=\bH$, the estimate \eqref{eq:analysis-mc-quadratic-loss-L2}, and the Euler--Maruyama second-moment estimate \eqref{eq:analysis-weak-quadratic-moment-em}, we obtain
\begin{equation*}
\begin{aligned}
&\max_{0\le n\le N_t}
\left(\E\left[\left|\Phi_N(\bX_n^{\mathrm{BUG}})-\E[\varphi_H(X(t_n))]\right|^2\right]\right)^{1/2}\\
&\quad\leq
\|\bH\|_2\sqrt{C_T^{\mathrm{EM,stab}}e^{cT}}(\delta_{r,N}+\gamma_r)
+C_T^{\mathrm{EM,mom}}\|\bH\|_{\mathrm F}\Delta t
+\frac{\sigma_{T,H}+\|\bH\|_2\sigma_{T,I}}{\sqrt N}\\
&\quad\leq
C_{T,H}^{\mathrm{MC}}(\delta_{r,N}+\gamma_r+\Delta t)
+\frac{\sigma_{T,\varphi_H}}{\sqrt N},
\end{aligned}
\end{equation*}
where
\begin{equation*}
\begin{aligned}
C_{T,H}^{\mathrm{MC}}&=\max\left\{\|\bH\|_2\sqrt{C_T^{\mathrm{EM,stab}}e^{cT}},C_T^{\mathrm{EM,mom}}\|\bH\|_{\mathrm F}\right\},\\
\sigma_{T,\varphi_H}&=\sigma_{T,H}+\|\bH\|_2\sigma_{T,I}.
\end{aligned}
\end{equation*}
\end{proof}

\section{Numerical experiments}
\label{sec:numerical-experiments}

We assess the closed-moment EM--BUG integrator on two linear test problems: an Ornstein--Uhlenbeck process with additive noise and a linear equation with commuting multiplicative noise. These problems admit closed second-moment equations and exact reference solutions. The objectives of these experiments are threefold:
\begin{enumerate}
\item evaluate the accuracy and robustness of the closed-moment EM--BUG integrator with respect to the rank $r$;
\item validate its strong convergence with respect to the time-step size $\Delta t$;
\item assess the convergence of statistical quantities with respect to the number of particles $N$.
\end{enumerate}
The closed-moment EM--BUG integrator is also compared with the standard BUG integrator \cite{ceruti2022rank,ceruti2024parallel}, which serves as the reference method.

\subsection{Setup and error metrics}
\label{sec:numerical-setup}

The initial ensemble $\bX_0$ consists of independent and identically distributed realizations of $X(0)$, independent of the Brownian motions. The Euler--Maruyama scheme is initialized from this ensemble, and the initial low-rank ensemble $\bU_0\bY_0^\top$ is then obtained from a rank-$r$ truncated SVD of $\bX_0$. For the closed-moment EM--BUG integrator, the reduced second moment is initialized as $\bD_0=\bY_0^\top\bY_0/N$, and all additive and multiplicative noise directions are included in the augmented physical basis. These choices satisfy the assumptions on particle independence, reduced second-moment initialization, and complete noise augmentation used in the theoretical analysis.

To ensure a consistent comparison, the exact solution, the Euler--Maruyama scheme, and the BUG integrators are coupled through the same Brownian paths. In time-step convergence studies, Brownian increments are generated on the finest grid and then aggregated on coarser grids.

The strong convergence is assessed using the maximum-in-time error
\begin{equation}
\label{eq:strong-error}
e_{\mathrm{strong}}=\max_{0\leq n\leq N_t}\frac{1}{\sqrt{N}}\left\|\bX_n-\bX(t_n)\right\|_{\mathrm F},
\end{equation}
where $\bX(t_n)$ denotes the exact ensemble and $\bX_n$ denotes the ensemble computed by either the Euler--Maruyama scheme or the BUG integrators. In addition, we assess the approximation of statistical quantities using the errors in the mean and non-centered second moment at the final time $T$:
\begin{align}
\label{eq:mean-error}
e_{\mathrm{mean}} &= \left\|\E[X(T)]-\frac{1}{N}\bX_{N_t}\boldsymbol{1}_N\right\|_2, \\ \label{eq:moment-error}
e_{\mathrm{moment}} &= \left\|\E[X(T)X(T)^\top]-\frac{1}{N}\bX_{N_t}\bX_{N_t}^\top\right\|_{\mathrm F}.
\end{align}
These errors include the Monte Carlo sampling error, the time-discretization error, and the low-rank approximation error. To distinguish these contributions, we also report the corresponding errors obtained from the exact Monte Carlo ensemble, which includes only the sampling error, and the Euler--Maruyama scheme, which includes the sampling and time-discretization errors.

\subsection{Ornstein--Uhlenbeck process with additive noise}
\label{sec:ou-problem}

We first consider the Ornstein--Uhlenbeck process
\begin{equation}
\label{eq:ou-sde}
\dd X(t) = \bA X(t)\,\dd t + \sigma\bPsi_k\,\dd W(t),
\end{equation}
where $W$ is a $k$-dimensional Brownian motion. The drift matrix is given by
\begin{equation*}
\bA = -\alpha\bI+\beta\bL,
\end{equation*}
where $\bL$ is the standard one-dimensional discrete Laplacian. Since $\bA$ is symmetric, it admits the eigendecomposition
\begin{equation*}
\bA = \bPsi\bLambda\bPsi^\top,
\qquad
\bLambda = \operatorname{diag}(\lambda_1,\ldots,\lambda_d).
\end{equation*}
The noise directions are then chosen as the first $k$ eigenvectors of $\bA$:
\begin{equation*}
\bPsi_k=\bPsi(:,1:k).
\end{equation*}
We consider the deterministic initial condition
\begin{equation}
\label{eq:ou-initial-condition}
X(0)=\be_1,
\end{equation}
where $\be_1=[1,0,\ldots,0]^\top\in\R^d$. The associated initial mean and non-centered second moment are
\begin{equation*}
\E[X(0)]=\be_1,
\qquad
\E[X(0)X(0)^\top]=\be_1\be_1^\top.
\end{equation*}
Hence, the initial ensemble and second-moment matrix have rank one, and the second-moment matrix is singular. Moreover, since the range of $\bPsi_k$ is invariant under $\bA$, the stochastic component of the solution remains in this $k$-dimensional subspace, while the mean contributes to at most one additional direction. Consequently,
\begin{equation*}
\rank(\bX(t))\leq k+1.
\end{equation*}
For this problem, the augmented physical basis of the closed-moment EM--BUG integrator is constructed as
\begin{equation*}
\widehat{\bU}_{n+1}=\orth\begin{bmatrix}\bM\bU_n & \bPsi_k\end{bmatrix}.
\end{equation*}
This augmentation includes all additive-noise directions and therefore ensures Euler--Maruyama exactness for $r\geq k+1$. In the following, we use $d=50$, $k=14$, $\alpha=1$, $\beta=0.1$, $\sigma=0.08$, and $T=1$.

\paragraph{Exact solution and moments.}
Using the modal coordinates $Z(t)=\bPsi^\top X(t)$, the equation \eqref{eq:ou-sde} decouples into
\begin{equation*}
\dd Z_i(t) =
\begin{cases}
\lambda_i Z_i(t)\,\dd t+\sigma\,\dd W_i(t),
& i\leq k,\\
\lambda_i Z_i(t)\,\dd t,
& i>k,
\end{cases}
\end{equation*}
for $i=1,\ldots,d$. Over one time step, the exact transition is
\begin{equation*}
Z_i(t_{n+1})=e^{\lambda_i\Delta t}Z_i(t_n)+\eta_{n,i},
\end{equation*}
where
\begin{equation*}
\eta_{n,i}=
\begin{cases}
\displaystyle
\sigma\int_{t_n}^{t_{n+1}}
e^{\lambda_i(t_{n+1}-s)}\,\dd W_i(s),
& i\leq k,\\
0, & i>k.
\end{cases}
\end{equation*}
Let
\begin{equation*}
c_i(h)=\int_0^h e^{\lambda_i s}\,\dd s,
\qquad
v_i(h)=\int_0^h e^{2\lambda_i s}\,\dd s.
\end{equation*}
The increments $\eta_{n,i}$ are centered Gaussian random variables with variance
\begin{equation*}
\operatorname{Var}(\eta_{n,i})=
\begin{cases}
\sigma^2 v_i(\Delta t), & i\leq k,\\
0, & i>k.
\end{cases}
\end{equation*}
To couple the exact solution with the Brownian increments $\Delta W_{n,i}$ used by the integrators, we generate
\begin{equation*}
\eta_{n,i}=
\begin{cases}
\displaystyle
\sigma\left(\frac{c_i(\Delta t)}{\Delta t}\Delta W_{n,i}+\sqrt{v_i(\Delta t)-\frac{c_i(\Delta t)^2}{\Delta t}}\,\zeta_{n,i}\right), & i\leq k,\\
0, & i>k,
\end{cases}
\end{equation*}
where $\zeta_{n,i}$ are independent standard Gaussian random variables, independent of the Brownian increments. This construction preserves both the variance of the exact stochastic increment and its covariance with the Brownian increment. The exact solution is then recovered as $X(t_n)=\bPsi Z(t_n)$. Given the initial condition \eqref{eq:ou-initial-condition}, the exact mean and non-centered second moment are
\begin{align*}
\E[X(t)]&=e^{\bA t}\E[X(0)],\\ \label{eq:ou-exact-moment}
\E[X(t)X(t)^\top]&=e^{\bA t}\E[X(0)X(0)^\top] e^{\bA^\top t}+\sigma^2\bPsi_k\operatorname{diag}(v_1(t),\ldots,v_k(t))\bPsi_k^\top.
\end{align*}
These expressions provide the exact references for the errors \eqref{eq:mean-error}--\eqref{eq:moment-error}.

\paragraph{Numerical results.}

\begin{figure}[H]
\center
\includegraphics[width=0.8\textwidth]{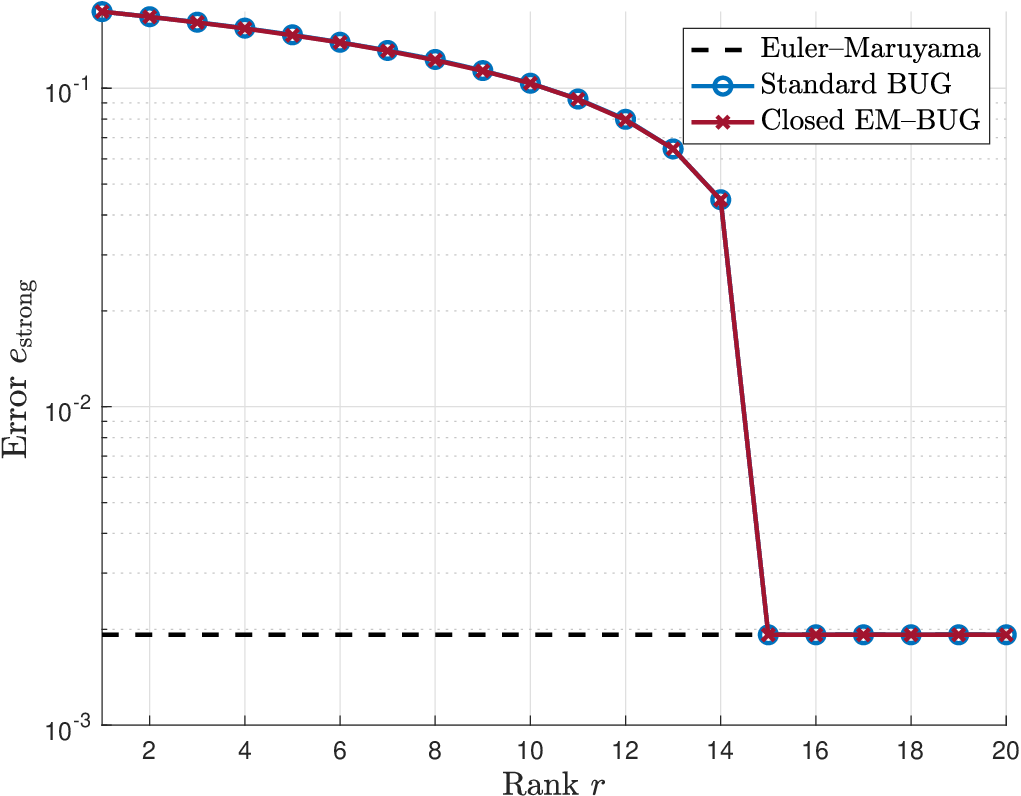}
\caption{\label{fig:ou-rank}Convergence with respect to the rank $r$ for the Ornstein--Uhlenbeck problem, with $N=5000$ and $\Delta t=1/128$.}
\end{figure}

Figure~\ref{fig:ou-rank} shows the error as a function of the rank. The standard BUG and closed-moment EM--BUG integrators yield nearly identical results, with an error decreasing as the rank increases. For $r\geq k+1=15$, the error coincides with that of the full Euler--Maruyama scheme, which is consistent with the exactness property.

\begin{figure}[H]
\begin{subfigure}[c]{0.5\linewidth}
\center
\includegraphics[width=0.9\textwidth]{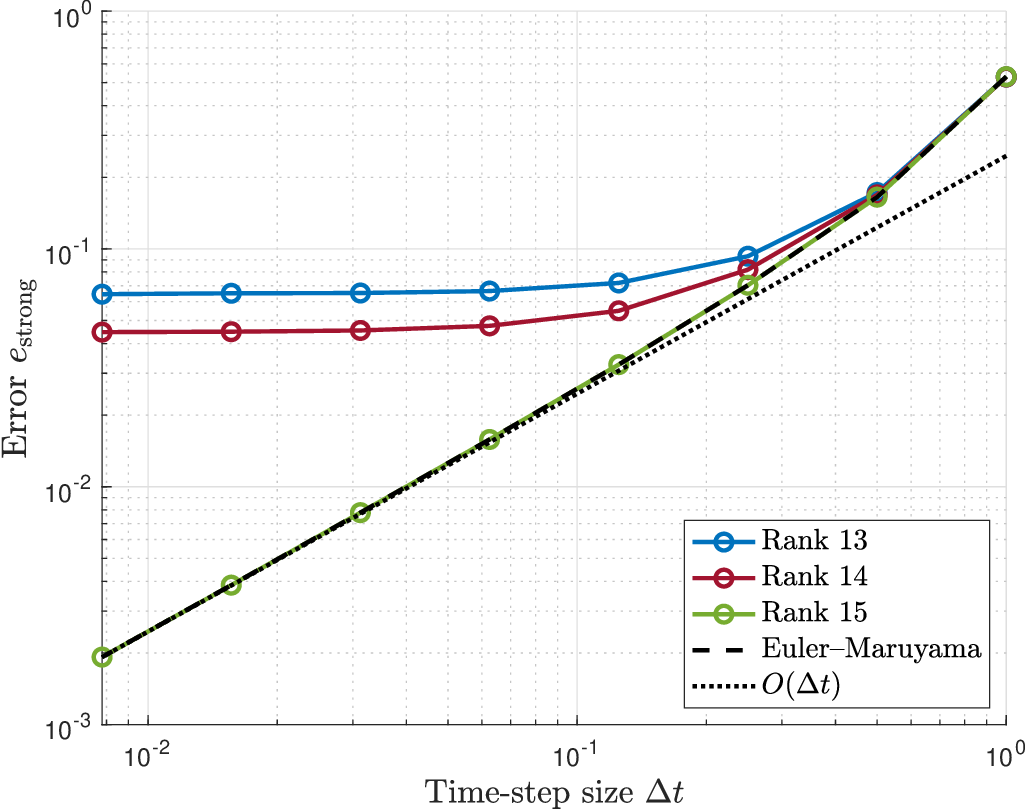}
\caption{Standard BUG}
\end{subfigure}\hfill
\begin{subfigure}[c]{0.5\linewidth}
\center
\includegraphics[width=0.9\textwidth]{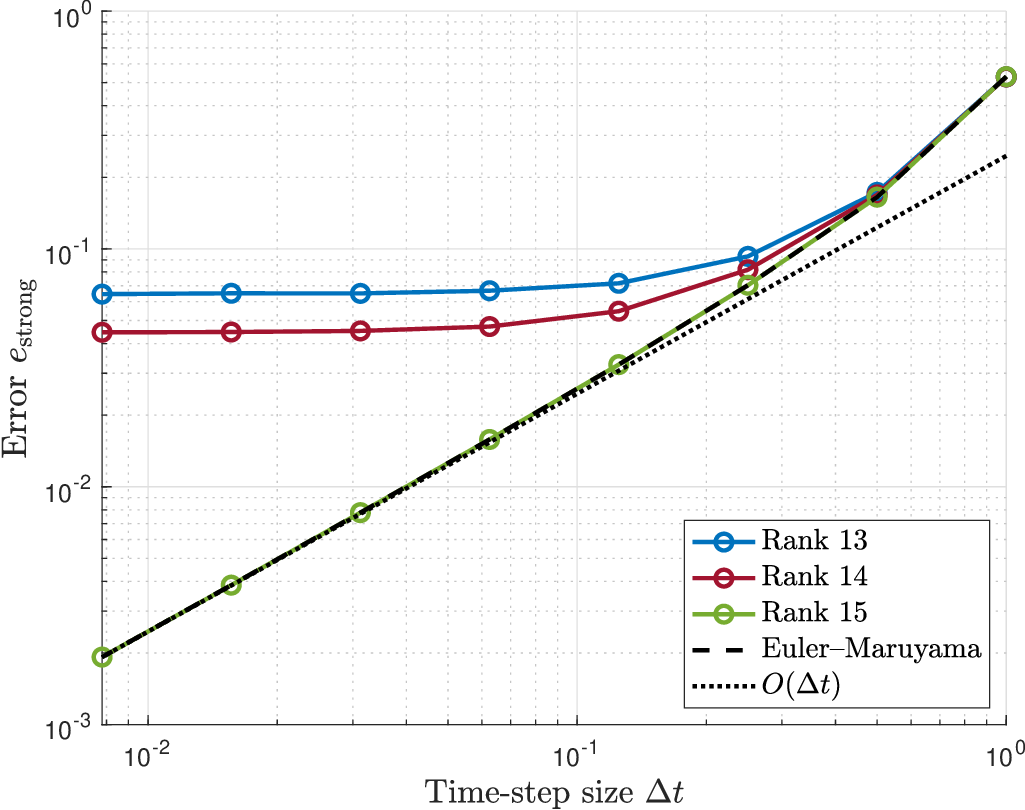}
\caption{Closed-moment EM--BUG}
\end{subfigure}\hfill
\caption{\label{fig:ou-timestep}Convergence with respect to the time-step size $\Delta t$ for the Ornstein--Uhlenbeck problem, with $N=5000$.}
\end{figure}

Figure~\ref{fig:ou-timestep} presents the error as a function of the time-step size for $r\in\{13,14,15\}$. The two BUG integrators yield nearly identical results. For larger time-step sizes, the error exhibits first-order convergence in $\Delta t$, before reaching a plateau due to the low-rank approximation error. This plateau decreases as the rank increases and disappears for $r=15$, as the integrators recover the Euler--Maruyama scheme.

\begin{figure}[H]
\begin{subfigure}[c]{\linewidth}
\center
\includegraphics[width=0.9\textwidth]{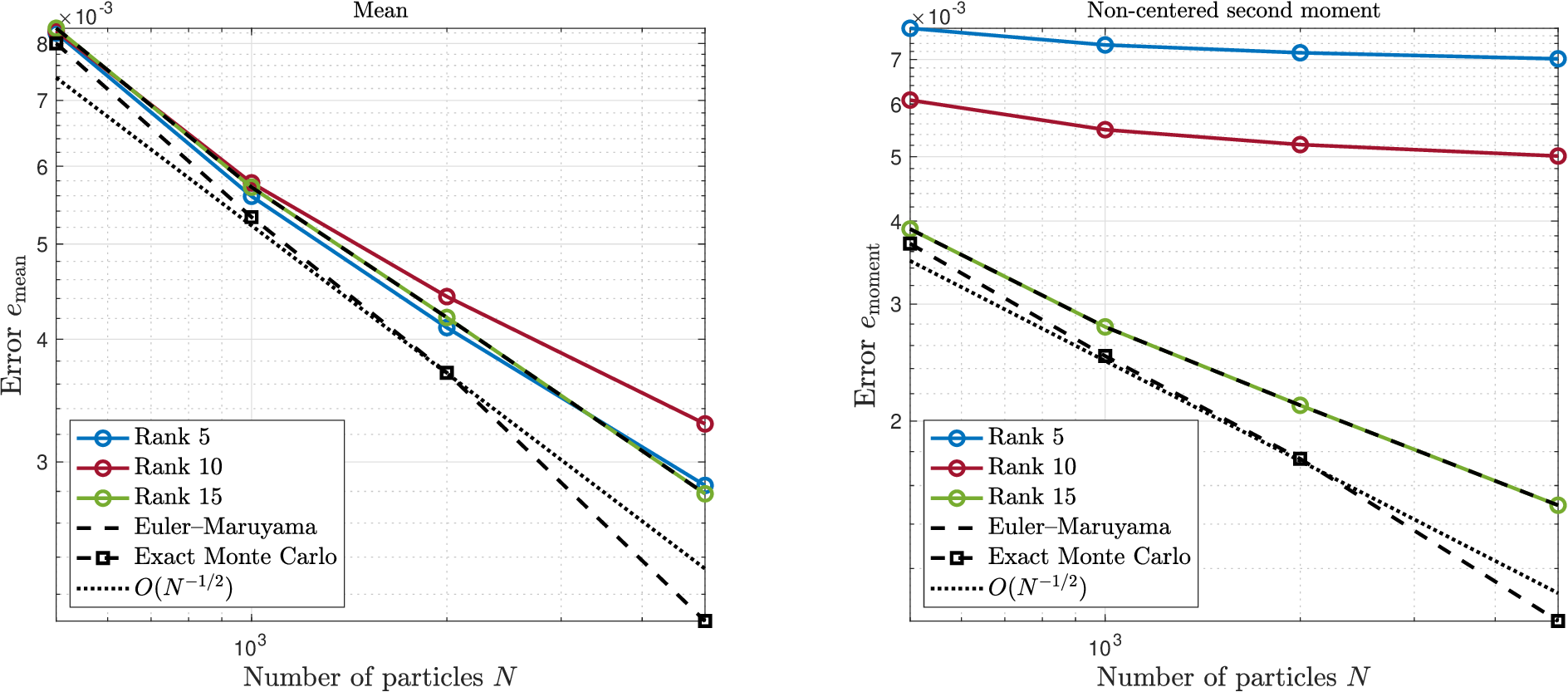}
\caption{Standard BUG}
\end{subfigure}\hfill
\begin{subfigure}[c]{\linewidth}
\center
\includegraphics[width=0.9\textwidth]{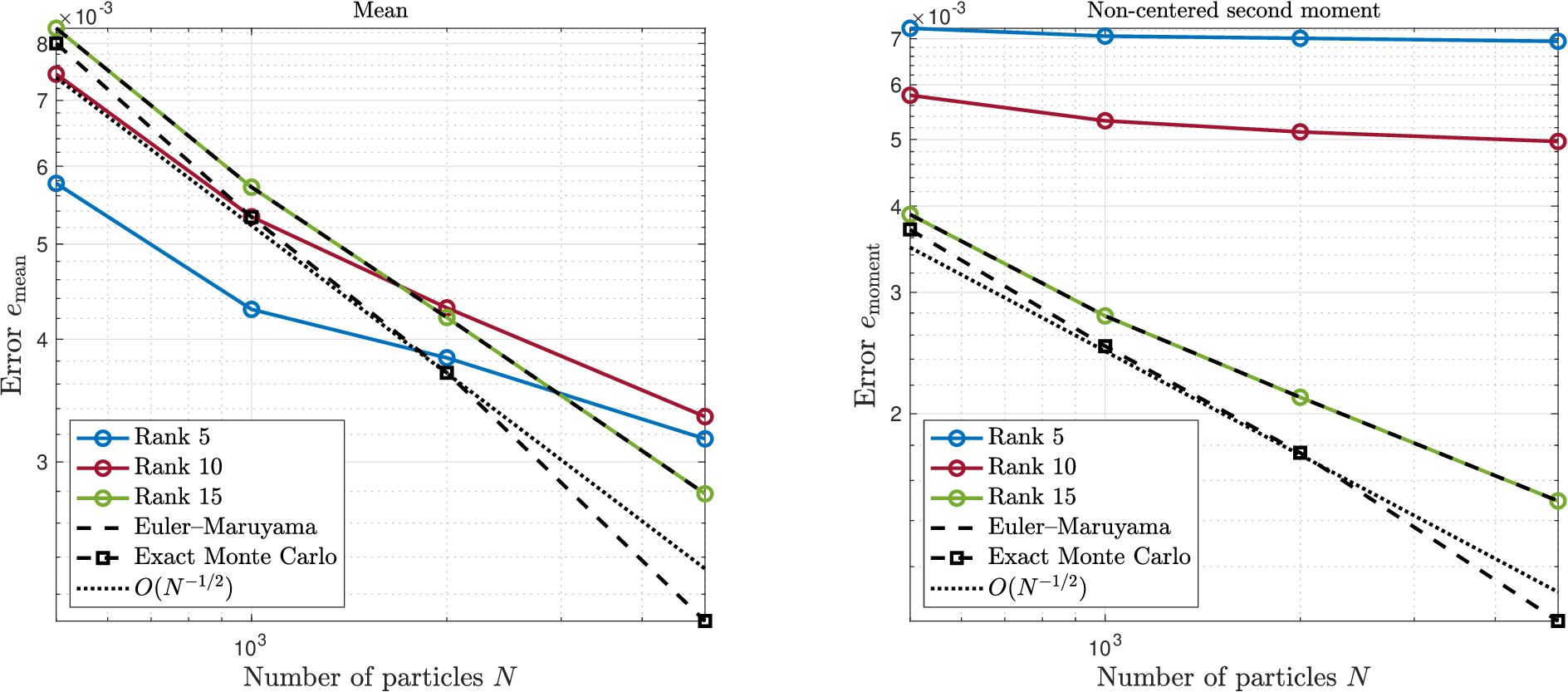}
\caption{Closed-moment EM--BUG}
\end{subfigure}\hfill
\caption{\label{fig:ou-particles}Convergence with respect to the number of particles $N$ for the Ornstein--Uhlenbeck problem, with $\Delta t=1/128$ and averaged over 5 runs.}
\end{figure}

Figure~\ref{fig:ou-particles} reports the error in the mean and non-centered second moment as a function of the number of particles. The two BUG integrators yield similar results. As the rank increases, the error approaches that of the Euler--Maruyama scheme, until it coincides for $r=15$. However, the error does not necessarily decrease monotonically with the rank, particularly for the mean. As $N$ increases, the error tends to decrease, but this reduction may be limited by the low-rank approximation error, especially for the second moment at smaller ranks.

\subsection{Linear equation with commuting multiplicative noise}
\label{sec:multiplicative-problem}

We then consider the linear stochastic differential equation
\begin{equation}
\label{eq:multiplicative-problem}
\dd X(t)=\bA X(t)\,\dd t+\sum_{\ell=1}^{m}\bB_\ell X(t)\,\dd W_\ell(t),
\end{equation}
where $W_\ell$ are independent scalar Brownian motions. The drift and diffusion matrices are
\begin{equation*}
\bA=\bPsi\operatorname{diag}(\lambda_1,\ldots,\lambda_d)\bPsi^\top,
\qquad
\bB_\ell=\bPsi\operatorname{diag}(\mu_{1\ell},\ldots,\mu_{d\ell})\bPsi^\top,
\end{equation*}
where $\bPsi$ is the orthogonal matrix containing the eigenvectors of the standard one-dimensional discrete Laplacian, and
\begin{equation*}
\lambda_i=-0.6-1.2s_i,
\qquad
\mu_{i\ell}=\frac{0.35}{\sqrt m}\cos\!\left(\ell\pi\left(s_i+\frac{1}{2d}\right)\right),
\qquad
s_i=\frac{i-1}{d-1},
\end{equation*}
for $i=1,\ldots,d$ and $\ell=1,\ldots,m$. We consider the deterministic initial condition given in modal coordinates $Z(0)=\bPsi^\top X(0)$ by
\begin{equation}
\label{eq:multiplicative-initial-condition}
Z_i(0)=
\begin{cases}
z_i, & i\leq k,\\
0, & i>k,
\end{cases}
\qquad
z_i=\frac{e^{-4s_i}}{\sqrt{\sum_{j=1}^{k}e^{-8s_j}}}.
\end{equation}
The associated initial mean and non-centered second moment are
\begin{equation*}
\E[X(0)]=\bPsi_k\boldsymbol z,
\qquad
\E[X(0)X(0)^\top]=\bPsi_k\bz\bz^\top\bPsi_k^\top,
\end{equation*}
where $\bPsi_k=\bPsi(:,1:k)$ and $\bz=[z_1,\ldots,z_k]^\top$. Hence, the initial ensemble and second-moment matrix have rank one, and the second-moment matrix is singular. Moreover, since the drift and diffusion matrices are simultaneously diagonalizable and the initial condition belongs to the range of $\bPsi_k$, the solution remains in this $k$-dimensional subspace. Consequently,
\begin{equation*}
\rank(\bX(t))\leq k.
\end{equation*}
For this problem, the augmented physical basis of the closed-moment EM--BUG integrator is constructed as
\begin{equation*}
\widehat{\bU}_{n+1}
=\orth\begin{bmatrix}
\bM\bU_n & \bB_1\bU_n & \cdots & \bB_m\bU_n
\end{bmatrix}.
\end{equation*}
This augmentation includes all multiplicative-noise directions and therefore ensures Euler--Maruyama exactness for $r\geq k$. In the following, we use $d=50$, $m=4$, $k=15$, and $T=1$.

\paragraph{Exact solution and moments.}
Using the modal coordinates $Z(t)=\bPsi^\top X(t)$, the equation \eqref{eq:multiplicative-problem} decouples into
\begin{equation*}
\dd Z_i(t)=\lambda_i Z_i(t)\,\dd t+\sum_{\ell=1}^{m}\mu_{i\ell}Z_i(t)\,\dd W_\ell(t),
\end{equation*}
for $i=1,\ldots,d$. Over one time step, the exact transition is
\begin{equation*}
Z_i(t_{n+1})=Z_i(t_n)\exp\left[\left(\lambda_i-\frac12\sum_{\ell=1}^{m}\mu_{i\ell}^2\right)\Delta t+\sum_{\ell=1}^{m}\mu_{i\ell}\Delta W_{\ell,n}\right].
\end{equation*}
The exact solution is then recovered as $X(t_n)=\bPsi Z(t_n)$. Since the exact transition uses the same Brownian increments as the Euler--Maruyama scheme and BUG integrators, the trajectories are directly coupled. Given the initial condition \eqref{eq:multiplicative-initial-condition}, the exact mean and non-centered second moment are
\begin{align*}
\E[X(t)]&=\bPsi\E[Z(t)],\\
\E[X(t)X(t)^\top]&=\bPsi\E[Z(t)Z(t)^\top]\bPsi^\top,
\end{align*}
where
\begin{align*}
\E[Z_i(t)]&=Z_i(0)e^{\lambda_i t},\\
\E[Z_i(t)Z_j(t)]&=Z_i(0)Z_j(0)\exp\!\left[\left(\lambda_i+\lambda_j+\sum_{\ell=1}^m\mu_{i\ell}\mu_{j\ell}\right)t\right].
\end{align*}
These expressions provide the exact references for the errors \eqref{eq:mean-error}--\eqref{eq:moment-error}.

\paragraph{Numerical results.}

\begin{figure}[H]
\center
\includegraphics[width=0.8\textwidth]{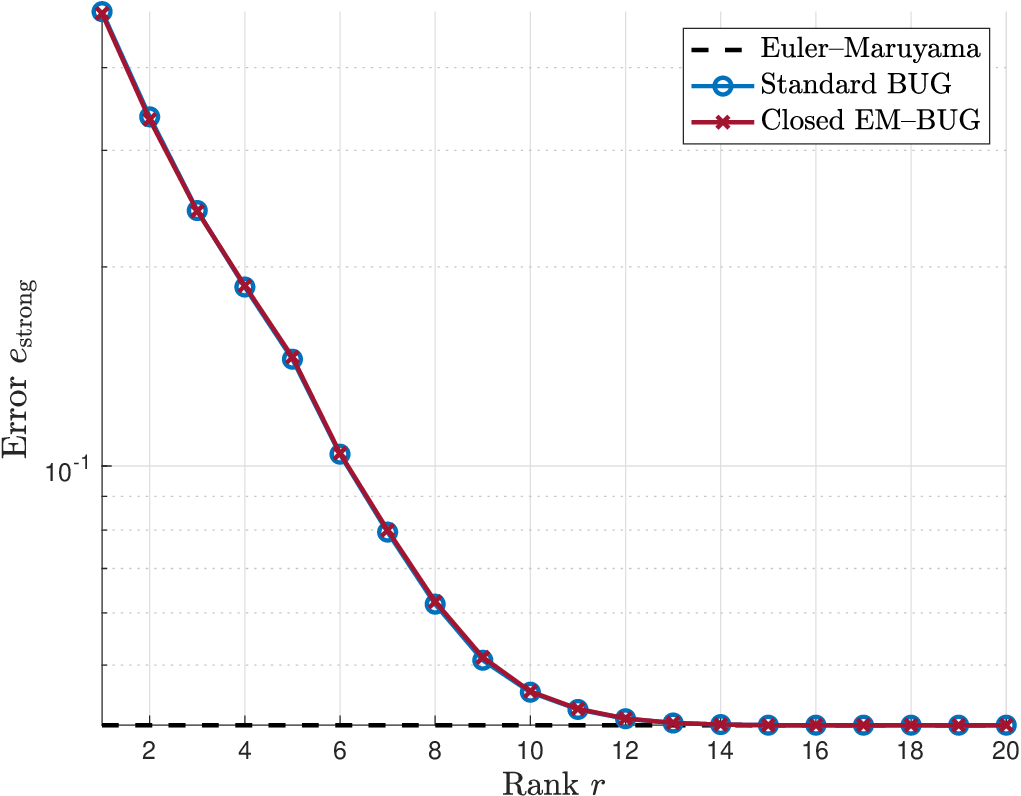}
\caption{\label{fig:multiplicative-rank}Convergence with respect to the rank $r$ for the commuting multiplicative-noise problem, with $N=5000$ and $\Delta t=1/128$.}
\end{figure}

Figure~\ref{fig:multiplicative-rank} shows the error as a function of the rank. The standard BUG and closed-moment EM--BUG integrators yield nearly identical results, with an error decreasing as the rank increases. For $r\geq k=15$, the error coincides with that of the full Euler--Maruyama scheme, which is consistent with the exactness property.

\begin{figure}[H]
\begin{subfigure}[c]{0.5\linewidth}
\center
\includegraphics[width=0.9\textwidth]{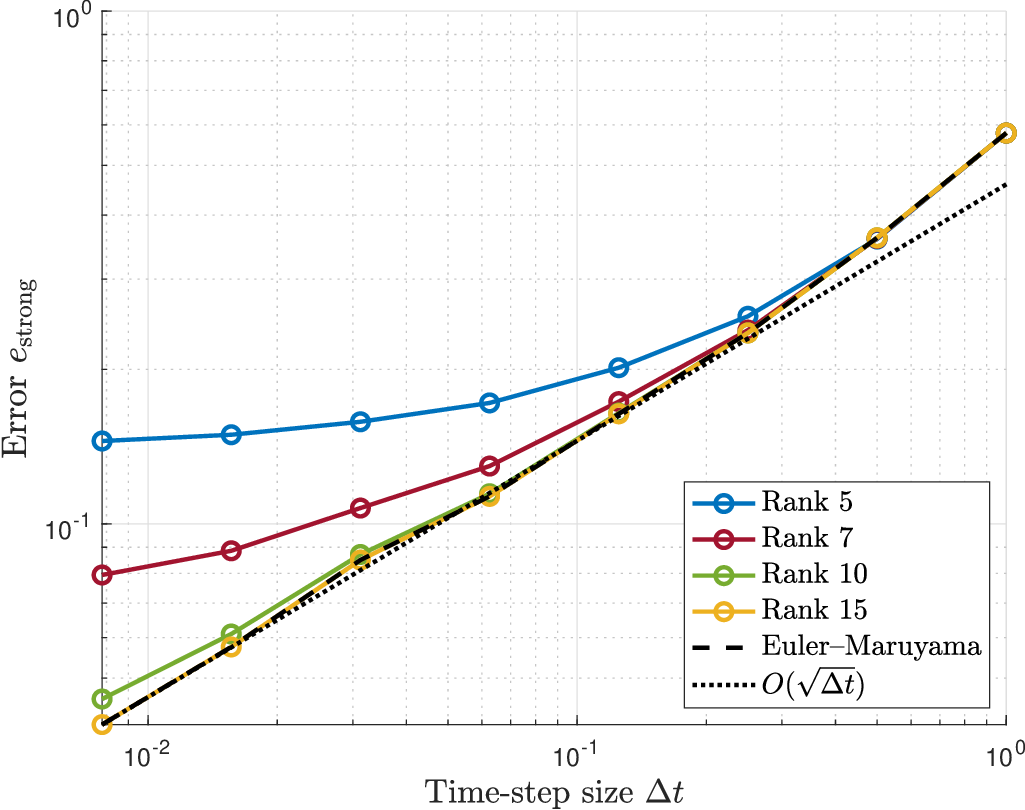}
\caption{Standard BUG}
\end{subfigure}\hfill
\begin{subfigure}[c]{0.5\linewidth}
\center
\includegraphics[width=0.9\textwidth]{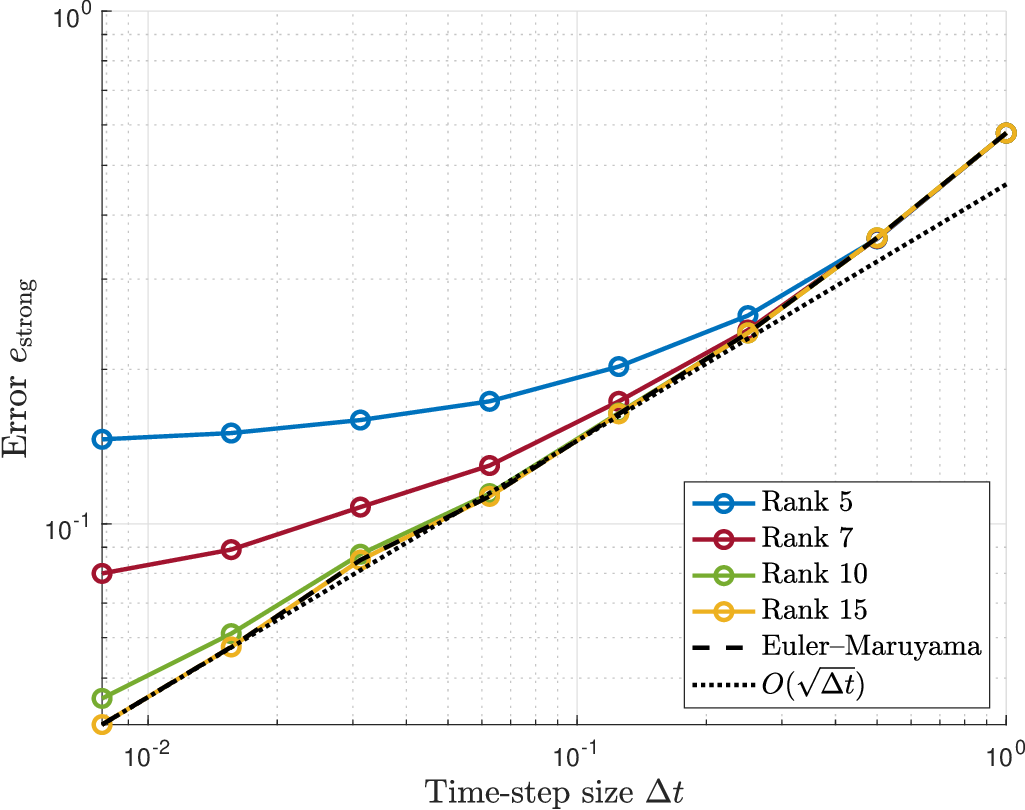}
\caption{Closed-moment EM--BUG}
\end{subfigure}\hfill
\caption{\label{fig:multiplicative-timestep}Convergence with respect to the time-step size $\Delta t$ for the commuting multiplicative-noise problem, with $N=5000$.}
\end{figure}

Figure~\ref{fig:multiplicative-timestep} presents the error as a function of the time-step size for $r\in\{5,7,10,15\}$. The two BUG integrators yield nearly identical results. The error decreases with $\Delta t$, following the convergence behavior of the Euler--Maruyama scheme, before reaching a plateau due to the low-rank approximation error. This plateau decreases as the rank increases and disappears for $r=15$, as the integrators recover the Euler--Maruyama scheme.

\begin{figure}[H]
\begin{subfigure}[c]{\linewidth}
\center
\includegraphics[width=0.9\textwidth]{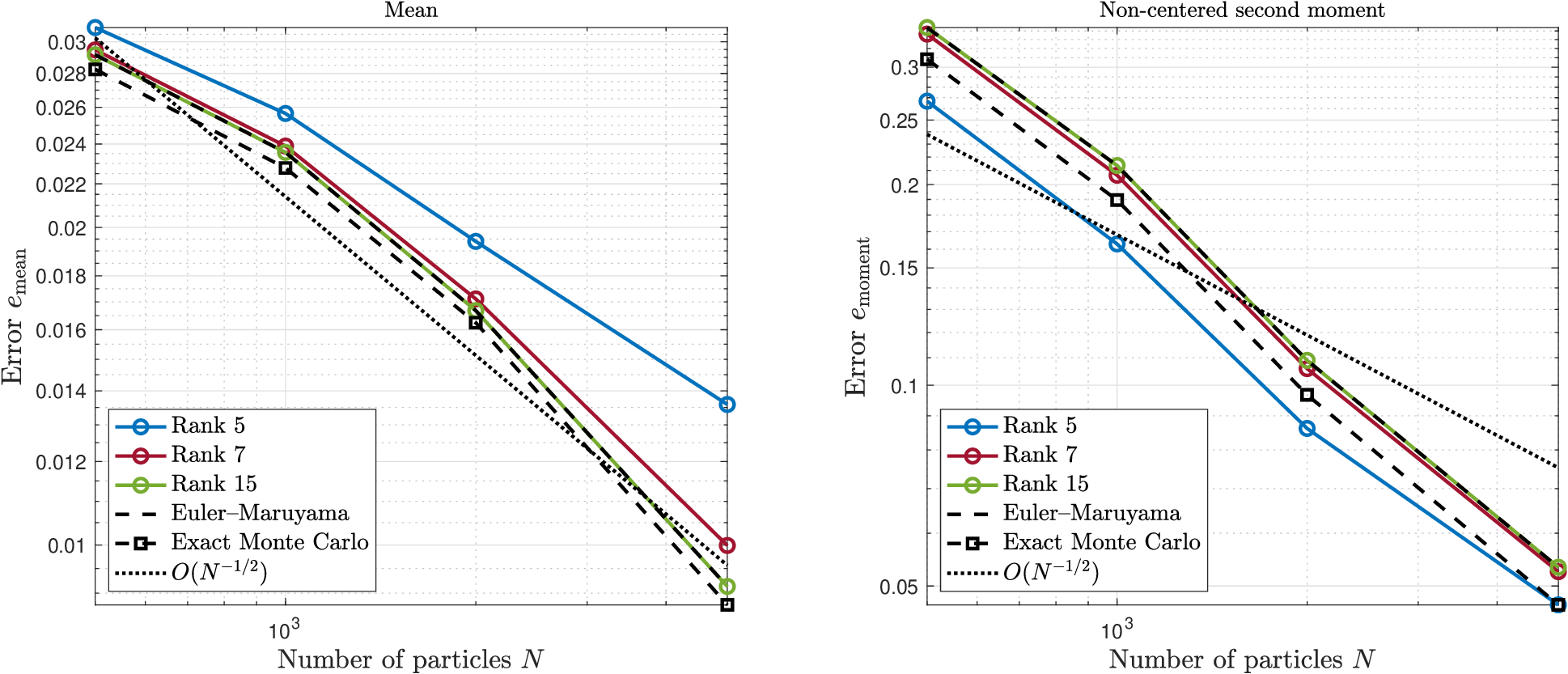}
\caption{Standard BUG}
\end{subfigure}\hfill
\begin{subfigure}[c]{\linewidth}
\center
\includegraphics[width=0.9\textwidth]{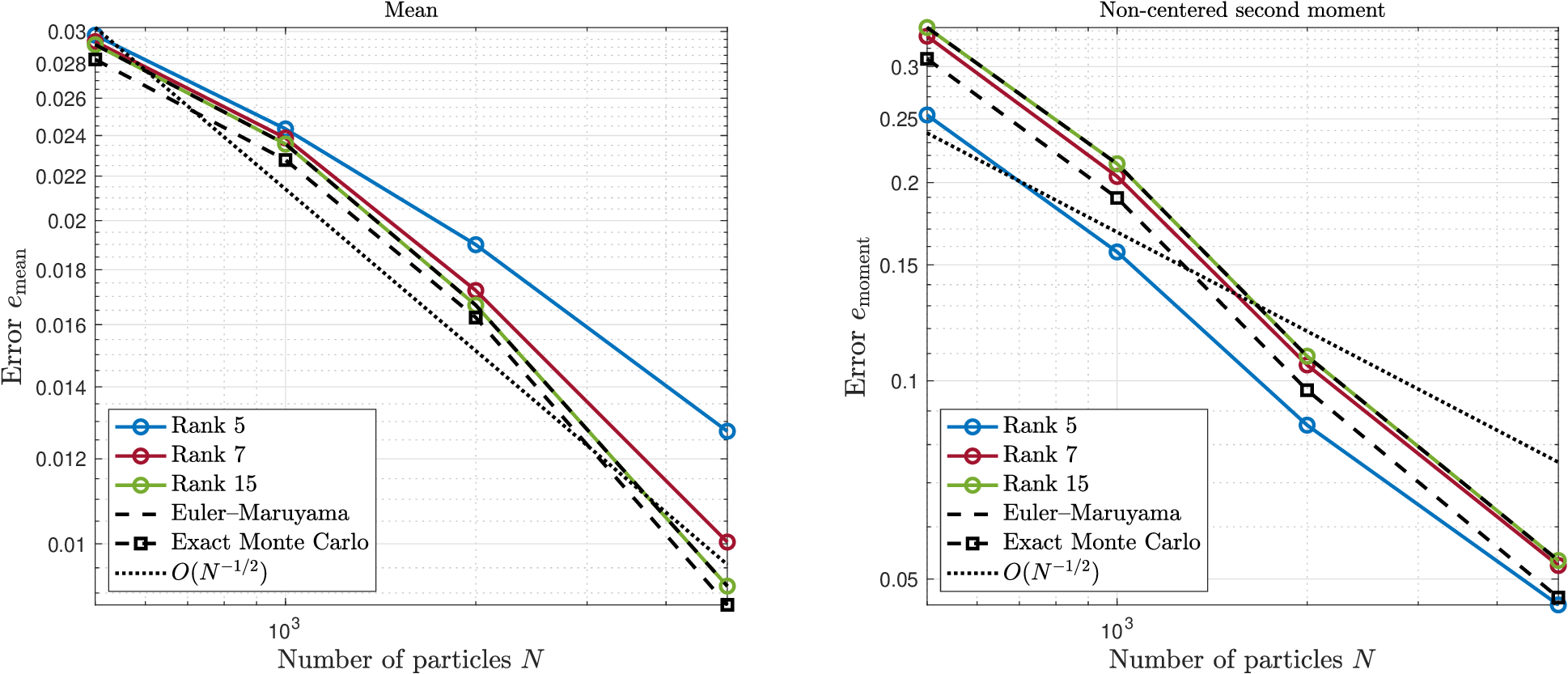}
\caption{Closed-moment EM--BUG}
\end{subfigure}\hfill
\caption{\label{fig:multiplicative-particles}Convergence with respect to the number of particles $N$ for the commuting multiplicative-noise problem, with $\Delta t=1/128$ and averaged over 5 runs.}
\end{figure}

Figure~\ref{fig:multiplicative-particles} reports the error in the mean and non-centered second moment as a function of the number of particles. The two BUG integrators yield similar results. For $r=15$, the error coincides with that of the Euler--Maruyama scheme. However, it does not necessarily decrease monotonically with the rank, particularly for the second moment. As $N$ increases, the error tends to decrease, exhibiting a behavior consistent with the Monte Carlo convergence.

\section{Conclusion}
\label{sec:conclusion}

In this paper, we have introduced BUG integrators based on the Euler--Maruyama method for SDEs. The proposed EM--BUG integrators extend the standard BUG formulation by enriching the physical basis with directions associated with the diffusion and exploiting the second moment to perform the Galerkin projection and rank truncation independently of the new Brownian increments.

For the closed-moment EM--BUG integrator, we have established several theoretical properties, including preservation of positive semidefiniteness, exactness with respect to Euler--Maruyama at full rank, and second-moment stability. We have also derived strong and weak convergence estimates showing that the time-discretization accuracy of the underlying Euler--Maruyama scheme is preserved up to the low-rank approximation error, with bounds that are robust with respect to small eigenvalues of the reduced second-moment matrix. In addition, Monte Carlo error estimates recover the standard sampling rate. Numerical experiments on linear SDEs with additive and multiplicative noise support these theoretical results and show that the proposed integrator is robust and achieves accuracy comparable to the standard BUG integrator.

Future work will focus on the analysis of the empirical conditional-moment EM--BUG integrator for general nonlinear SDEs, which will be addressed in Part II. Other possible directions include extending the proposed approach to higher-order stochastic time-integration schemes and exploiting the rank adaptivity of the EM--BUG integrators for applications such as ensemble-based data assimilation.

\begin{appendices}
\section{Auxiliary estimates for the analysis}

\begin{lemma}[Auxiliary estimates for the strong comparison]
\label{lem:analysis-strong-closed-aux}
Consider the closed-moment EM--BUG integrator with the reduced second moment initialized as $\bD_0=\bY_0^\top\bY_0/N$. Assume that the closed-moment EM--BUG integrator and the Euler–Maruyama scheme use the same Wiener increments. Let $\bX_n^{\mathrm{aux}}$ be the auxiliary ensemble initialized by $\bX_0^{\mathrm{aux}}=\bX_0^{\mathrm{BUG}}$ and defined by
\begin{equation*}
\label{eq:analysis-strong-closed-aux-update}
\bX_{n+1}^{\mathrm{aux}}=\bM\bX_{n}^{\mathrm{aux}}+\bB_{0,n+1}^{\parallel}\Delta\bW_{0,n}+\sum_{\ell=1}^q\bB_\ell\bX_{n}^{\mathrm{aux}}\Delta\bW_{\ell,n},
\end{equation*}
where $\bB_{0,n+1}^{\parallel}=\widehat{\bP}_{n+1}\bB_0$. Define also the spectral truncation defect
\begin{equation}
\label{eq:analysis-strong-closed-tau}
\tau_{n+1}:=\operatorname{tr}\left(\widehat{\bP}_{n+1}\bS_{n+1}\widehat{\bP}_{n+1}-\bP_{n+1}\bS_{n+1}\bP_{n+1}\right).
\end{equation}
Then,
\begin{equation}
\label{eq:analysis-strong-closed-aux-lr}
\max_{0\le n\le N_t}\E\left[\|\bX_n^{\mathrm{BUG}}-\bX_n^{\mathrm{aux}}\|_N^2\right]\leq C_T^{\mathrm{BUG,aux}}\sum_{n=0}^{N_t-1}\E[\tau_{n+1}],
\end{equation}
and
\begin{equation}
\label{eq:analysis-strong-closed-aux-em}
\E\left[\max_{0\le n\le N_t}\|\bX_n^{\mathrm{aux}}-\bX_n^{\mathrm{EM}}\|_N^2\right]\leq C_T^{\mathrm{EM,aux}}\left(\delta_{r,N}^2+\Delta t\sum_{n=0}^{N_t-1}\E\left[\|(\bI-\widehat{\bP}_{n+1})\bB_0\|_{\mathrm F}^2\right]\right),
\end{equation}
where
\begin{equation*}
\begin{aligned}
C_T^{\mathrm{EM,aux}}&=12\exp\left(T\left[3T\|\bA\|_2^2+12\sum_{\ell=1}^q\|\bB_\ell\|_2^2\right]\right),\\
C_T^{\mathrm{BUG,aux}}&=\exp\left(T\left[2\|\bA\|_2+T\|\bA\|_2^2+\sum_{\ell=1}^q\|\bB_\ell\|_2^2\right]\right).
\end{aligned}
\end{equation*}
\end{lemma}

\begin{proof}
Since the augmented basis $\widehat{\bU}_{n+1}$ contains the ranges of $\bM\bU_n$ and $\{\bB_\ell\bU_n\}_{\ell=1}^q$, the projected second-moment prediction satisfies
\begin{equation}
\label{eq:analysis-strong-closed-aux-prediction}
\begin{aligned}
\bS_{n+1}^{\parallel}&=\widehat{\bP}_{n+1}\bS_{n+1}\widehat{\bP}_{n+1}\\
&=\bM\bC_n^{\mathrm{BUG}}\bM^\top+\Delta t\,\bB_{0,n+1}^{\parallel}(\bB_{0,n+1}^{\parallel})^\top+\Delta t\sum_{\ell=1}^q\bB_\ell\bC_n^{\mathrm{BUG}}\bB_\ell^\top.
\end{aligned}
\end{equation}
Let $\mathcal G=\sigma(\bX_0^{\mathrm{EM}})$ denote the sigma-algebra generated by the initial Euler--Maruyama ensemble. Since the initial low-rank ensemble $\bX_0^{\mathrm{BUG}}$ is obtained from $\bX_0^{\mathrm{EM}}$ and the closed second-moment recursion is deterministic, the projectors $\widehat{\bP}_{n+1}$ and $\bP_{n+1}$ are $\mathcal G$-measurable. In the following, we therefore condition on $\mathcal G$, so that these projectors are deterministic conditionally on $\mathcal G$.

\emph{Step 1: Second-moment identities.}
Let $\bC_n^{\mathrm{aux}}=\E[\bX_n^{\mathrm{aux}}(\bX_n^{\mathrm{aux}})^\top\mid\mathcal G]/N$ and $\bC^{\mathrm{BUG,aux}}_n=\E[\bX_n^{\mathrm{aux}}(\bX_n^{\mathrm{BUG}})^\top\mid\mathcal G]/N$. Since the projectors are $\mathcal G$-measurable and the auxiliary and BUG ensembles use the same Wiener increments, their conditional second and cross moments satisfy
\begin{align}
\bC_{n+1}^{\mathrm{aux}}&=\bM\bC_n^{\mathrm{aux}}\bM^\top+\Delta t\,\bB_{0,n+1}^{\parallel}(\bB_{0,n+1}^{\parallel})^\top+\Delta t\sum_{\ell=1}^q\bB_\ell\bC_n^{\mathrm{aux}}\bB_\ell^\top \label{eq:analysis-strong-closed-aux-moment}\\
\bC^{\mathrm{BUG,aux}}_{n+1}&=\left(\bM\bC^{\mathrm{BUG,aux}}_n\bM^\top+\Delta t\,\bB_{0,n+1}^{\parallel}(\bB_{0,n+1}^{\parallel})^\top+\Delta t\sum_{\ell=1}^q\bB_\ell\bC^{\mathrm{BUG,aux}}_n\bB_\ell^\top\right)\bP_{n+1}. \label{eq:analysis-strong-closed-aux-cross}
\end{align}
Since $\bX_0^{\mathrm{aux}}=\bX_0^{\mathrm{BUG}}$ and $\bD_0=\bY_0^\top\bY_0/N$, we have $\bC^{\mathrm{BUG,aux}}_0=\bC_0^{\mathrm{BUG}}$. Moreover, since the columns of $\bU_{n+1}$ are eigenvectors of $\bS_{n+1}^{\parallel}$, the orthogonal projector $\bP_{n+1}=\bU_{n+1}\bU_{n+1}^\top$ satisfies
\begin{equation}
\label{eq:analysis-strong-closed-commute}
\bP_{n+1}\bS_{n+1}^{\parallel}=\bP_{n+1}\bS_{n+1}^{\parallel}\bP_{n+1}=\bS_{n+1}^{\parallel}\bP_{n+1}=\bC_{n+1}^{\mathrm{BUG}}.
\end{equation}
Consequently, using \eqref{eq:analysis-strong-closed-aux-cross}, \eqref{eq:analysis-strong-closed-aux-prediction}, and \eqref{eq:analysis-strong-closed-commute}, we obtain by induction $\bC^{\mathrm{BUG,aux}}_n=\bC_n^{\mathrm{BUG}}$ for every $n$.

Similarly, let $\bC_n^{\mathrm{BUG,cond}}=\E\left[\bX_n^{\mathrm{BUG}}(\bX_n^{\mathrm{BUG}})^\top\mid\mathcal G\right]/N$. Since $\bD_0=\bY_0^\top\bY_0/N$, we have $\bC_0^{\mathrm{BUG,cond}}=\bC_0^{\mathrm{BUG}}$. Moreover, since the projectors are $\mathcal G$-measurable,
\begin{equation*}
\bC_{n+1}^{\mathrm{BUG,cond}}=\bP_{n+1}\left(\bM\bC_n^{\mathrm{BUG,cond}}\bM^\top+\Delta t\,\bB_0\bB_0^\top+\Delta t\sum_{\ell=1}^q\bB_\ell\bC_n^{\mathrm{BUG,cond}}\bB_\ell^\top\right)\bP_{n+1}.
\end{equation*}
Consequently, we obtain by induction $\bC_n^{\mathrm{BUG,cond}}=\bC_n^{\mathrm{BUG}}$ for every $n$. Finally, since $\bC_n^{\mathrm{BUG,cond}}=\bC^{\mathrm{BUG,aux}}_n=\bC_n^{\mathrm{BUG}}$ and $\bC_n^{\mathrm{BUG}}$ is symmetric, it follows that
\begin{equation}
\label{eq:analysis-strong-closed-aux-error-moment}
\begin{aligned}
&\frac1N\E\left[(\bX_n^{\mathrm{aux}}-\bX_n^{\mathrm{BUG}})(\bX_n^{\mathrm{aux}}-\bX_n^{\mathrm{BUG}})^\top\mid\mathcal G\right]\\
&\qquad=\bC_n^{\mathrm{aux}}+\bC_n^{\mathrm{BUG,cond}}-\bC^{\mathrm{BUG,aux}}_n-(\bC^{\mathrm{BUG,aux}}_n)^\top\\
&\qquad=\bC_n^{\mathrm{aux}}-\bC_n^{\mathrm{BUG}}.
\end{aligned}
\end{equation}

\emph{Step 2: Control by the spectral truncation defect.}
Let $\bTheta_n=\bC_n^{\mathrm{aux}}-\bC_n^{\mathrm{BUG}}$. Using \eqref{eq:analysis-strong-closed-aux-prediction} and \eqref{eq:analysis-strong-closed-aux-moment}, we obtain
\begin{equation}
\label{eq:analysis-strong-closed-aux-H}
\bTheta_{n+1}=\bM\bTheta_n\bM^\top+\Delta t\sum_{\ell=1}^q\bB_\ell\bTheta_n\bB_\ell^\top+\bS_{n+1}^{\parallel}-\bP_{n+1}\bS_{n+1}^{\parallel}\bP_{n+1}.
\end{equation}
Since $\bC_n^{\mathrm{BUG}}\succeq0$, equation~\eqref{eq:analysis-strong-closed-aux-prediction} implies that $\bS_{n+1}^{\parallel}\succeq0$. Moreover, since $\bP_{n+1}$ and $\bS_{n+1}^{\parallel}$ commute according to \eqref{eq:analysis-strong-closed-commute}, we have
\begin{equation*}
\bS_{n+1}^{\parallel}-\bP_{n+1}\bS_{n+1}^{\parallel}\bP_{n+1}=(\bI-\bP_{n+1})\bS_{n+1}^{\parallel}(\bI-\bP_{n+1})\succeq0.
\end{equation*}
Consequently, since $\bTheta_0=0$ and the terms on the right-hand side of \eqref{eq:analysis-strong-closed-aux-H} are positive semidefinite when $\bTheta_n\succeq0$, we obtain by induction $\bTheta_n\succeq0$ for all $n$. It follows that, since $\bTheta_n\succeq0$ and the Löwner order is preserved under congruence transformations, 
\begin{equation*}
\begin{aligned}
\operatorname{tr}\left(\bM\bTheta_n\bM^\top+\Delta t\sum_{\ell=1}^q\bB_\ell\bTheta_n\bB_\ell^\top\right)&\leq\left(\|\bM\|_2^2+\Delta t\sum_{\ell=1}^q\|\bB_\ell\|_2^2\right)\operatorname{tr}(\bTheta_n)\\
&\leq(1+c_1\Delta t)\operatorname{tr}(\bTheta_n),
\end{aligned}
\end{equation*}
where $c_1=2\|\bA\|_2+T\|\bA\|_2^2+\sum_{\ell=1}^q\|\bB_\ell\|_2^2$. Taking the trace of \eqref{eq:analysis-strong-closed-aux-H} gives
\begin{equation*}
\operatorname{tr}(\bTheta_{n+1})\leq(1+c_1\Delta t)\operatorname{tr}(\bTheta_n)+\tau_{n+1},
\end{equation*}
and the discrete Gr\"onwall inequality yields
\begin{equation*}
\max_{0\le n\le N_t}\operatorname{tr}(\bTheta_n)\leq e^{c_1T}\sum_{n=0}^{N_t-1}\tau_{n+1}.
\end{equation*}
Taking the trace of \eqref{eq:analysis-strong-closed-aux-error-moment} and then taking expectations proves \eqref{eq:analysis-strong-closed-aux-lr}.

\emph{Step 3: Comparison with Euler--Maruyama.}
Let $\bX_n^{\mathrm{diff}}=\bX_n^{\mathrm{aux}}-\bX_n^{\mathrm{EM}}$. Since the auxiliary and EM ensembles use the same Wiener increments, their difference satisfies
\begin{equation*}
\bX_{n+1}^{\mathrm{diff}}=\bM\bX_n^{\mathrm{diff}}+\boldsymbol\xi_{n+1},
\end{equation*}
where
\begin{equation*}
\boldsymbol\xi_{n+1}=\sum_{\ell=1}^q\bB_\ell \bX_n^{\mathrm{diff}}\Delta \bW_{\ell,n}-(\bI-\widehat{\bP}_{n+1})\bB_0\Delta\bW_{0,n}.
\end{equation*}
Since $\bX_n^{\mathrm{diff}}$ and $\widehat{\bP}_{n+1}$ are $\mathcal F_{t_n}$-measurable, and the Wiener increments are independent of $\mathcal F_{t_n}$ and centered, we have $\E[\boldsymbol\xi_{n+1}\mid\mathcal F_{t_n}]=0$. Moreover, the conditional second moments of the Wiener increments satisfy
\begin{equation*}
\E[\Delta\bW_{\ell,n}\Delta\bW_{\ell,n}^\top\mid\mathcal F_{t_n}]=\Delta t\,\bI_N,\qquad \E[\Delta\bW_{0,n}\Delta\bW_{0,n}^\top\mid\mathcal F_{t_n}]=N\Delta t\,\bI_{m_0}.
\end{equation*}
Consequently, the independence of the Wiener components implies that the cross terms vanish, and we obtain
\begin{equation}
\label{eq:analysis-strong-closed-aux-noise}
\begin{aligned}
\E\left[\|\boldsymbol\xi_{n+1}\|_N^2\mid\mathcal F_{t_n}\right]&=\Delta t\sum_{\ell=1}^q\|\bB_\ell\bX_n^{\mathrm{diff}}\|_N^2+\Delta t\|(\bI-\widehat{\bP}_{n+1})\bB_0\|_{\mathrm F}^2\\
&\leq c_2\Delta t\|\bX_n^{\mathrm{diff}}\|_N^2+\Delta t\|(\bI-\widehat{\bP}_{n+1})\bB_0\|_{\mathrm F}^2,
\end{aligned}
\end{equation}
where $c_2=\sum_{\ell=1}^q\|\bB_\ell\|_2^2$. Summing the recursion of $\bX_n^{\mathrm{diff}}$ from $k=0$ to $n-1$, for any $n\leq N_t$, yields
\begin{equation*}
\bX_n^{\mathrm{diff}}=\bX_0^{\mathrm{diff}}+\Delta t\sum_{k=0}^{n-1}\bA\bX_k^{\mathrm{diff}}+\sum_{k=0}^{n-1}\boldsymbol\xi_{k+1}.
\end{equation*}
By Cauchy--Schwarz and since $n\Delta t\leq T$, we have
\begin{equation*}
\max_{0\leq j\leq n}\left\|\Delta t\sum_{k=0}^{j-1}\bA\bX_k^{\mathrm{diff}}\right\|_N^2\leq T\|\bA\|_2^2\Delta t\sum_{k=0}^{n-1}\|\bX_k^{\mathrm{diff}}\|_N^2.
\end{equation*}
Moreover, since $\E[\boldsymbol\xi_{k+1}\mid\mathcal F_{t_k}]=0$, the process $\sum_{k=0}^{j-1}\boldsymbol\xi_{k+1}$ is a martingale. Doob's $L^2$ maximal inequality and the orthogonality of martingale increments give
\begin{equation*}
\E\left[\max_{0\leq j\leq n}\left\|\sum_{k=0}^{j-1}\boldsymbol\xi_{k+1}\right\|_N^2\right]\leq4\sum_{k=0}^{n-1}\E[\|\boldsymbol\xi_{k+1}\|_N^2].
\end{equation*}
Consequently, using $\|a+b+c\|_N^2\leq3(\|a\|_N^2+\|b\|_N^2+\|c\|_N^2)$, we obtain
\begin{equation*}
\E\left[\max_{0\leq j\leq n}\|\bX_j^{\mathrm{diff}}\|_N^2\right]\leq3\E[\|\bX_0^{\mathrm{diff}}\|_N^2]+3T\|\bA\|_2^2\Delta t\sum_{k=0}^{n-1}\E[\|\bX_k^{\mathrm{diff}}\|_N^2]+12\sum_{k=0}^{n-1}\E[\|\boldsymbol\xi_{k+1}\|_N^2].
\end{equation*}
Using \eqref{eq:analysis-strong-closed-aux-noise}, the discrete Gr\"onwall inequality yields
\begin{equation*}
\E\left[\max_{0\le n\le N_t}\|\bX_n^{\mathrm{diff}}\|_N^2\right]\leq C_T^{\mathrm{EM,aux}}\left(\E[\|\bX_0^{\mathrm{diff}}\|_N^2]+\Delta t\sum_{n=0}^{N_t-1}\E\left[\|(\bI-\widehat{\bP}_{n+1})\bB_0\|_{\mathrm F}^2\right]\right).
\end{equation*}
Finally, since $\E[\|\bX_0^{\mathrm{diff}}\|_N^2]=\delta_{r,N}^2$, this proves \eqref{eq:analysis-strong-closed-aux-em}.
\end{proof}

\begin{lemma}[Maximal estimate for the strong comparison with additive noise]
\label{lem:analysis-strong-closed-max}
Suppose that the assumptions of Lemma~\ref{lem:analysis-strong-closed-aux} hold. Then, for additive noise ($q=0$),
\begin{equation}
\label{eq:analysis-strong-closed-max}
\E\left[\max_{0\le n\le N_t}\|\bX_n^{\mathrm{BUG}}-\bX_n^{\mathrm{aux}}\|_N^2\right]\leq C_T^{\mathrm{BUG,max}}\sum_{n=0}^{N_t-1}\E[\tau_{n+1}],
\end{equation}
where
\begin{equation*}
C_T^{\mathrm{BUG,max}}=32\exp\left(2T^2\|\bA\|_2^2\right).
\end{equation*}
\end{lemma}

\begin{proof}
As in the proof of Lemma~\ref{lem:analysis-strong-closed-aux}, we condition on the initial ensemble $\bX_0^{\mathrm{EM}}$, whose sigma-algebra is denoted by $\mathcal G$, so that the projectors are deterministic conditionally on $\mathcal G$. 

\emph{Step 1: Gaussian initial condition.}
Let $g\sim\mathcal N(0,\bI_d)$ be independent of the Wiener increments and of $\mathcal G$. Consider the two trajectories, $X_n^{\mathrm{aux},g}$ and $X_n^{\mathrm{BUG},g}$, initialized by $X_0^{\mathrm{aux},g}=X_0^{\mathrm{BUG},g}=X_0^g:=(\bC_0^{\mathrm{BUG}})^{1/2}g$ and defined by
\begin{equation*}
\begin{aligned}
X_{n+1}^{\mathrm{aux},g}&=\bM X_n^{\mathrm{aux},g}+\bB_{0,n+1}^{\parallel}\Delta W_{0,n},\\
X_{n+1}^{\mathrm{BUG},g}&=\bP_{n+1}\left(\bM X_n^{\mathrm{BUG},g}+\bB_0\Delta W_{0,n}\right),
\end{aligned}
\end{equation*}
where the Wiener increments and the projectors are those of the closed-moment EM--BUG integrator starting from the initial condition $X_0$, and are not recomputed from the Gaussian initial condition $X_0^g$.

Let $e_n^g=X_n^{\mathrm{aux},g}-X_n^{\mathrm{BUG},g}$. Since $\E[X_0^g(X_0^g)^\top\mid\mathcal G]=\bC_0^{\mathrm{BUG}}$, the conditional second-moment and cross-moment identities established in the proof of Lemma~\ref{lem:analysis-strong-closed-aux} remain valid, and we obtain
\begin{equation*}
\begin{aligned}
\E\left[e_n^g(X_n^{\mathrm{BUG},g})^\top\mid\mathcal G\right]&=\E\left[X_n^{\mathrm{aux},g}(X_n^{\mathrm{BUG},g})^\top\mid\mathcal G\right]-\E\left[X_n^{\mathrm{BUG},g}(X_n^{\mathrm{BUG},g})^\top\mid\mathcal G\right]\\
&=\bC_n^{\mathrm{BUG}}-\bC_n^{\mathrm{BUG}}=0.
\end{aligned}
\end{equation*}
Moreover, since the projectors are $\mathcal G$-measurable and the new Wiener increments are centered and independent of the past, we have, for $k\leq n$,
\begin{equation*}
\E\left[e_k^g(X_{n+1}^{\mathrm{BUG},g})^\top\mid\mathcal G\right]=\E\left[e_k^g(X_n^{\mathrm{BUG},g})^\top\mid\mathcal G\right]\bM^\top\bP_{n+1}.
\end{equation*}
Since $\E\left[e_k^g(X_k^{\mathrm{BUG},g})^\top\mid\mathcal G\right]=0$, it follows by induction that
\begin{equation}
\label{eq:analysis-strong-closed-max-orthogonality}
\E\left[e_k^g(X_n^{\mathrm{BUG},g})^\top\mid\mathcal G\right]=0,\qquad 0\leq k\leq n.
\end{equation}
Since $q=0$, the error satisfies
\begin{equation}
\label{eq:analysis-strong-closed-max-recursion}
e_{n+1}^g=\bM e_n^g+\zeta_{n+1}^g,
\end{equation}
where
\begin{equation*}
\zeta_{n+1}^g=(\bI-\bP_{n+1})\left(\bM X_n^{\mathrm{BUG},g}+\bB_{0,n+1}^{\parallel}\Delta W_{0,n}\right).
\end{equation*}
Using \eqref{eq:analysis-strong-closed-max-orthogonality} and the independence of the new Wiener increments, we obtain
\begin{equation}
\label{eq:analysis-strong-closed-max-independent}
\E\left[\zeta_{n+1}^g(e_k^g)^\top\mid\mathcal G\right]=0,\qquad 0\leq k\leq n.
\end{equation}
Since $X_0^g$ and the Wiener increments are Gaussian, and the projectors are deterministic conditionally on $\mathcal G$, the random vectors $(e_0^g,\ldots,e_n^g,\zeta_{n+1}^g)$ are jointly centered Gaussian conditionally on $\mathcal G$. Hence, the orthogonality relation \eqref{eq:analysis-strong-closed-max-independent} implies that $\zeta_{n+1}^g$ is independent of $(e_0^g,\ldots,e_n^g)$ conditionally on $\mathcal G$. Let $\mathcal E_n=\sigma(\mathcal G,e_0^g,\ldots,e_n^g)$ be the sigma-algebra generated by $\mathcal G$ and the errors up to time $t_n$. Then,
\begin{equation*}
\E[\zeta_{n+1}^g\mid\mathcal E_n]=0.
\end{equation*}
Consequently, $(\zeta_{n+1}^g)_{n\geq0}$ is a martingale difference sequence with respect to the filtration $(\mathcal E_n)_{n\geq0}$. Furthermore, by the second-moment recurrence established in Lemma~\ref{lem:analysis-strong-closed-aux}, we have $\E[X_n^{\mathrm{BUG},g}(X_n^{\mathrm{BUG},g})^\top\mid\mathcal G]=\bC_n^{\mathrm{BUG}}$. Since the Wiener increments satisfy $\E[\Delta W_{0,n}\Delta W_{0,n}^\top\mid\mathcal G]=\Delta t\,\bI_{m_0}$ and these increments are centered and independent of the past, the cross terms vanish, and we obtain
\begin{equation*}
\begin{aligned}
\E\left[\zeta_{n+1}^g(\zeta_{n+1}^g)^\top\mid\mathcal G\right]&=(\bI-\bP_{n+1})\left(\bM\bC_n^{\mathrm{BUG}}\bM^\top+\Delta t\bB_{0,n+1}^{\parallel}(\bB_{0,n+1}^{\parallel})^\top\right)(\bI-\bP_{n+1})\\
&=(\bI-\bP_{n+1})\bS_{n+1}^{\parallel}(\bI-\bP_{n+1})\\
&=\bS_{n+1}^{\parallel}-\bP_{n+1}\bS_{n+1}^{\parallel}\bP_{n+1},
\end{aligned}
\end{equation*}
where the second equality follows from \eqref{eq:analysis-strong-closed-aux-prediction} with $q=0$, and the last equality follows from the commutation property \eqref{eq:analysis-strong-closed-commute}. Then, taking the trace gives
\begin{equation}
\label{eq:analysis-strong-closed-max-innovation}
\E\left[\|\zeta_{n+1}^g\|_2^2\mid\mathcal G\right]=\operatorname{tr}\left(\bS_{n+1}^{\parallel}-\bP_{n+1}\bS_{n+1}^{\parallel}\bP_{n+1}\right)=\tau_{n+1}.
\end{equation}

Since $e_0^g=0$, summing the recursion \eqref{eq:analysis-strong-closed-max-recursion} gives
\begin{equation*}
e_n^g=\Delta t\sum_{k=0}^{n-1}\bA e_k^g+\sum_{k=0}^{n-1}\zeta_{k+1}^g.
\end{equation*}
By Cauchy--Schwarz and since $n\Delta t\leq T$, we have
\begin{equation*}
\max_{0\leq j\leq n}\left\|\Delta t\sum_{k=0}^{j-1}\bA e_k^g\right\|_2^2\leq T\|\bA\|_2^2\Delta t\sum_{k=0}^{n-1}\|e_k^g\|_2^2.
\end{equation*}
Moreover, since $(\zeta_{n+1}^g)_{n\geq0}$ is a martingale difference sequence, Doob's $L^2$ maximal inequality and the orthogonality of martingale increments yield
\begin{equation*}
\E\left[\max_{0\leq j\leq n}\left\|\sum_{k=0}^{j-1}\zeta_{k+1}^g\right\|_2^2\mid\mathcal G\right]\leq4\sum_{k=0}^{n-1}\E\left[\|\zeta_{k+1}^g\|_2^2\mid\mathcal G\right]=4\sum_{k=0}^{n-1}\tau_{k+1},
\end{equation*}
where the last equality follows from \eqref{eq:analysis-strong-closed-max-innovation}. Consequently, using $\|a+b\|_2^2\leq2(\|a\|_2^2+\|b\|_2^2)$, we obtain
\begin{equation*}
\begin{aligned}
\E\left[\max_{0\leq j\leq n}\|e_j^g\|_2^2\mid\mathcal G\right]&\leq2T\|\bA\|_2^2\Delta t\sum_{k=0}^{n-1}\E\left[\|e_k^g\|_2^2\mid\mathcal G\right]+8\sum_{k=0}^{n-1}\tau_{k+1}.
\end{aligned}
\end{equation*}
The discrete Gr\"onwall inequality then gives
\begin{equation}
\label{eq:analysis-strong-closed-max-gaussian}
\E\left[\max_{0\leq n\leq N_t}\|e_n^g\|_2^2\mid\mathcal G\right]\leq8\exp\left(2T^2\|\bA\|_2^2\right)\sum_{n=0}^{N_t-1}\tau_{n+1}.
\end{equation}

\emph{Step 2: Extension to the initial ensemble $\bX_0^{\mathrm{BUG}}$.}
We now consider the ensembles $\bX_n^{\mathrm{aux}}$ and $\bX_n^{\mathrm{BUG}}$ initialized by $\bX_0^{\mathrm{BUG}}$, instead of the Gaussian initial condition. Let $e_n^{(j)}=X_n^{\mathrm{aux},(j)}-X_n^{\mathrm{BUG},(j)}$. Since the recurrences for $X_n^{\mathrm{aux},(j)}$ and $X_n^{\mathrm{BUG},(j)}$ are linear with additive noise, the projectors are $\mathcal G$-measurable, and $\bP_{n+1}\bB_{0,n+1}^{\parallel}=\bP_{n+1}\bB_0$, we obtain by induction
\begin{equation*}
\begin{aligned}
X_n^{\mathrm{aux},(j)}&=\bM^n X_0^{\mathrm{BUG},(j)}+\sum_{k=0}^{n-1}\bM^{n-1-k}\bB_{0,k+1}^{\parallel}\Delta W_{0,k}^{(j)},\\
X_n^{\mathrm{BUG},(j)}&=\bT_{n,0}X_0^{\mathrm{BUG},(j)}+\sum_{k=0}^{n-1}\bT_{n,k+1}\bP_{k+1}\bB_{0,k+1}^{\parallel}\Delta W_{0,k}^{(j)},
\end{aligned}
\end{equation*}
where $\bT_{n,k}=\bP_n\bM\bP_{n-1}\bM\cdots\bP_{k+1}\bM$ for $0\leq k<n$, and $\bT_{n,n}=\bI$. Consequently,
\begin{equation*}
e_n^{(j)}=\bPi_n X_0^{\mathrm{BUG},(j)}+\chi_n^{(j)},
\end{equation*}
where $\bPi_n=\bM^n-\bT_{n,0}$ and $\chi_n^{(j)}=\sum_{k=0}^{n-1}(\bM^{n-1-k}-\bT_{n,k+1}\bP_{k+1})\bB_{0,k+1}^{\parallel}\Delta W_{0,k}^{(j)}$. Similarly, using the recursion \eqref{eq:analysis-strong-closed-max-recursion}, we obtain by induction
\begin{equation}
\label{eq:analysis-strong-closed-max-representation}
e_n^g=\bPi_n\left(\bC_0^{\mathrm{BUG}}\right)^{1/2}g+\chi_n,
\end{equation}
where
\begin{equation*}
\chi_n=\sum_{k=0}^{n-1}\left(\bM^{n-1-k}-\bT_{n,k+1}\bP_{k+1}\right)\bB_{0,k+1}^{\parallel}\Delta W_{0,k}.
\end{equation*}
In particular, conditionally on $\mathcal G$, the processes $(\chi_n^{(j)})_{n=0}^{N_t}$ and $(\chi_n)_{n=0}^{N_t}$ have the same law, and $(\chi_n)_{n=0}^{N_t}$ is independent of $g$.

Since $\bC_0^{\mathrm{BUG}}=\bX_0^{\mathrm{BUG}}(\bX_0^{\mathrm{BUG}})^\top/N$, we have
\begin{equation*}
\begin{aligned}
\frac1N\sum_{j=1}^N\|\boldsymbol{\Pi}_n X_0^{\mathrm{BUG},(j)}\|_2^2&=\frac1N\operatorname{tr}\left(\boldsymbol{\Pi}_n\bX_0^{\mathrm{BUG}}(\bX_0^{\mathrm{BUG}})^\top\boldsymbol{\Pi}_n^\top\right)\\
&=\|\boldsymbol{\Pi}_n(\bC_0^{\mathrm{BUG}})^{1/2}\|_{\mathrm F}^2.
\end{aligned}
\end{equation*}
Moreover, since $g$ is centered and independent of $\mathcal G$ and $(\chi_0,\ldots,\chi_{N_t})$, the identity \eqref{eq:analysis-strong-closed-max-representation} gives
\begin{equation*}
\E\left[\|e_n^g\|_2^2\mid\mathcal G\right]=\|\boldsymbol{\Pi}_n(\bC_0^{\mathrm{BUG}})^{1/2}\|_{\mathrm F}^2+\E\left[\|\chi_n\|_2^2\mid\mathcal G\right].
\end{equation*}
Consequently,
\begin{equation*}
\max_{0\le n\le N_t}\frac1N\sum_{j=1}^N\|\boldsymbol{\Pi}_n X_0^{\mathrm{BUG},(j)}\|_2^2\leq\E\left[\max_{0\le n\le N_t}\|e_n^g\|_2^2\mid\mathcal G\right].
\end{equation*}
Furthermore, since $g$ is centered and independent of $\mathcal G$ and $(\chi_0,\ldots,\chi_{N_t})$, the identity \eqref{eq:analysis-strong-closed-max-representation} gives
\begin{equation*}
\E[e_n^g\mid\mathcal G,\chi_0,\ldots,\chi_{N_t}]=\chi_n.
\end{equation*}
Since the function $(\chi_0,\ldots,\chi_{N_t})\mapsto\max_{0\le n\le N_t}\|\chi_n\|_2^2$ is convex, Jensen's inequality yields
\begin{equation*}
\max_{0\le n\le N_t}\|\chi_n\|_2^2\leq\E\left[\max_{0\le n\le N_t}\|e_n^g\|_2^2\mid\mathcal G,\chi_0,\ldots,\chi_{N_t}\right].
\end{equation*}
Taking conditional expectations with respect to $\mathcal G$, we obtain
\begin{equation*}
\E\left[\max_{0\le n\le N_t}\|\chi_n\|_2^2\mid\mathcal G\right]\leq\E\left[\max_{0\le n\le N_t}\|e_n^g\|_2^2\mid\mathcal G\right].
\end{equation*}
Then, combining $e_n^{(j)}=\bPi_nX_0^{\mathrm{BUG},(j)}+\chi_n^{(j)}$ with $\|a+b\|_2^2\leq2(\|a\|_2^2+\|b\|_2^2)$ gives
\begin{equation*}
\begin{aligned}
&\E\left[\max_{0\le n\le N_t}\sum_{j=1}^N\frac1N\|e_n^{(j)}\|_2^2\mid\mathcal G\right]\\
&\qquad\leq2\max_{0\le n\le N_t}\frac1N\sum_{j=1}^N\|\bPi_nX_0^{\mathrm{BUG},(j)}\|_2^2+\frac2N\sum_{j=1}^N\E\left[\max_{0\le n\le N_t}\|\chi_n^{(j)}\|_2^2\mid\mathcal G\right]\\
&\qquad=2\max_{0\le n\le N_t}\frac1N\sum_{j=1}^N\|\bPi_nX_0^{\mathrm{BUG},(j)}\|_2^2+2\E\left[\max_{0\le n\le N_t}\|\chi_n\|_2^2\mid\mathcal G\right]\\
&\qquad\leq4\E\left[\max_{0\le n\le N_t}\|e_n^g\|_2^2\mid\mathcal G\right],
\end{aligned}
\end{equation*}
where the equality follows from the identical conditional laws of $(\chi_n^{(j)})_{n=0}^{N_t}$ and $(\chi_n)_{n=0}^{N_t}$. Using \eqref{eq:analysis-strong-closed-max-gaussian}, we obtain
\begin{equation*}
\begin{aligned}
\E\left[\max_{0\le n\le N_t}\|\bX_n^{\mathrm{aux}}-\bX_n^{\mathrm{BUG}}\|_N^2\mid\mathcal G\right]&=\E\left[\max_{0\le n\le N_t}\sum_{j=1}^N\frac1N\|e_n^{(j)}\|_2^2\mid\mathcal G\right]\\
&\leq
32\exp(2T^2\|\bA\|_2^2)
\sum_{n=0}^{N_t-1}\tau_{n+1}.
\end{aligned}
\end{equation*}
Finally, taking expectations proves \eqref{eq:analysis-strong-closed-max}.
\end{proof}

\begin{lemma}[Auxiliary estimates for the quadratic Monte Carlo error]
\label{lem:analysis-mc-quadratic}
Suppose that the assumptions of Proposition~\ref{prop:analysis-mc-quadratic} hold. Then, there exist symmetric random matrices $\bXi_n$ and $\bGamma_n$ such that
\begin{equation}
\label{eq:analysis-mc-quadratic-decomposition}
\bC_{N,n}^{\mathrm{BUG}}=\bC_n^{\mathrm{EM}}+\bXi_n-\bGamma_n,
\end{equation}
where $\E[\bXi_n]=0$ and $\bGamma_n\succeq0$ satisfies
\begin{equation}
\label{eq:analysis-mc-quadratic-lemma-loss}
\begin{aligned}
\bGamma_0&=\frac1N\left(\bX_0^{\mathrm{EM}}(\bX_0^{\mathrm{EM}})^\top-\bX_0^{\mathrm{BUG}}(\bX_0^{\mathrm{BUG}})^\top\right),\\
\bGamma_{n+1}&=\bM\bGamma_n\bM^\top+\Delta t\sum_{\ell=1}^q\bB_\ell\bGamma_n\bB_\ell^\top+\bL_{n+1},
\end{aligned}
\end{equation}
with $\bL_{n+1}=\bS_{n+1}-\bP_{n+1}\bS_{n+1}\bP_{n+1}\succeq0$. Moreover, for every deterministic symmetric matrix $\bK\in\R^{d\times d}$,
\begin{equation}
\label{eq:analysis-mc-quadratic-fluctuation}
\max_{0\le n\le N_t}\left(\E\left[|\operatorname{tr}(\bK\bXi_n)|^2\right]\right)^{1/2}\leq\frac{\sigma_{T,K}}{\sqrt N},
\end{equation}
where $\sigma_{T,K}$ is independent of the rank $r$, the number of particles $N$, the time-step size $\Delta t$, the smallest eigenvalue of the reduced second moment, and spectral gaps.
\end{lemma}

\begin{proof}
\emph{Step 1: Uniform fourth-moment estimates.}
Since $\bP_{n+1}$ is $\mathcal{F}_{t_n}$-measurable, the closed-moment EM--BUG update can be written, for each particle, as
\begin{equation*}
X_{n+1}^{\mathrm{BUG},(j)}=\bP_{n+1}\left(\bM X_n^{\mathrm{BUG},(j)}+\boldsymbol\xi_n^{(j)}\right),
\end{equation*}
where $\boldsymbol\xi_n^{(j)}$ is a centered Gaussian vector, conditionally on $\mathcal{F}_{t_n}$, with covariance
\begin{equation*}
\bSigma_{n,j}=\Delta t\left(\bB_0\bB_0^\top+\sum_{\ell=1}^q\bB_\ell X_n^{\mathrm{BUG},(j)}(X_n^{\mathrm{BUG},(j)})^\top\bB_\ell^\top\right).
\end{equation*}
Since $\bP_{n+1}$ is an orthogonal projector, we have
\begin{equation*}
\begin{aligned}
\E\left[\|X_{n+1}^{\mathrm{BUG},(j)}\|_2^4\mid\mathcal{F}_{t_n}\right]&\leq\E\left[\left(\|\bM X_n^{\mathrm{BUG},(j)}\|_2^2+2(\bM X_n^{\mathrm{BUG},(j)})^\top\boldsymbol\xi_n^{(j)}+\|\boldsymbol\xi_n^{(j)}\|_2^2\right)^2\mid\mathcal{F}_{t_n}\right]\\
&=\|\bM X_n^{\mathrm{BUG},(j)}\|_2^4+2\|\bM X_n^{\mathrm{BUG},(j)}\|_2^2\operatorname{tr}(\bSigma_{n,j})\\
&\quad+4(\bM X_n^{\mathrm{BUG},(j)})^\top\bSigma_{n,j}(\bM X_n^{\mathrm{BUG},(j)})+(\operatorname{tr}(\bSigma_{n,j}))^2+2\operatorname{tr}(\bSigma_{n,j}^2)\\
&\leq\|\bM X_n^{\mathrm{BUG},(j)}\|_2^4+6\|\bM X_n^{\mathrm{BUG},(j)}\|_2^2\operatorname{tr}(\bSigma_{n,j})+3(\operatorname{tr}(\bSigma_{n,j}))^2,
\end{aligned}
\end{equation*}
where the equality follows from the vanishing of odd Gaussian moments and the Gaussian fourth-moment formula for $\boldsymbol\xi_n^{(j)}$, and the last inequality follows from $(\bM X_n^{\mathrm{BUG},(j)})^\top\bSigma_{n,j}(\bM X_n^{\mathrm{BUG},(j)})\leq\|\bM X_n^{\mathrm{BUG},(j)}\|_2^2\operatorname{tr}(\bSigma_{n,j})$ and $\operatorname{tr}(\bSigma_{n,j}^2)\leq(\operatorname{tr}(\bSigma_{n,j}))^2$, since $\bSigma_{n,j}\succeq0$. Moreover, we have
\begin{equation*}
\begin{aligned}
\operatorname{tr}(\bSigma_{n,j})&=\Delta t\left(\|\bB_0\|_{\mathrm F}^2+\sum_{\ell=1}^q\|\bB_\ell X_n^{\mathrm{BUG},(j)}\|_2^2\right)\\
&\leq\Delta t\left(\|\bB_0\|_{\mathrm F}^2+\sum_{\ell=1}^q\|\bB_\ell\|_2^2\|X_n^{\mathrm{BUG},(j)}\|_2^2\right)\\
&\leq c_1\Delta t(1+\|X_n^{\mathrm{BUG},(j)}\|_2^2),
\end{aligned}
\end{equation*}
with $c_1=\max\left(\|\bB_0\|_{\mathrm F}^2,\sum_{\ell=1}^q\|\bB_\ell\|_2^2\right)$. Since $s^2(1+s^2)\leq2(1+s^4)$ and $(1+s^2)^2\leq2(1+s^4)$ for $s\geq0$, using $\bM=\bI+\Delta t\bA$ and $\Delta t\leq T$ gives
\begin{equation*}
\begin{aligned}
\|\bM X_n^{\mathrm{BUG},(j)}\|_2^2\operatorname{tr}(\bSigma_{n,j})&\leq (1+T\|\bA\|_2)^2\|X_n^{\mathrm{BUG},(j)}\|_2^2c_1\Delta t\left(1+\|X_n^{\mathrm{BUG},(j)}\|_2^2\right)\\
&\quad\leq2c_1(1+T\|\bA\|_2)^2\Delta t\left(1+\|X_n^{\mathrm{BUG},(j)}\|_2^4\right),
\end{aligned}
\end{equation*}
and
\begin{equation*}
(\operatorname{tr}(\bSigma_{n,j}))^2\leq2Tc_1^2\Delta t(1+\|X_n^{\mathrm{BUG},(j)}\|_2^4).
\end{equation*}
Moreover, by the mean value theorem, we have
\begin{equation*}
\begin{aligned}
(1+\Delta t\|\bA\|_2)^4-1
&=4\|\bA\|_2\int_0^{\Delta t}(1+s\|\bA\|_2)^3\,ds\\
&\leq4\|\bA\|_2(1+\Delta t\|\bA\|_2)^3\Delta t.
\end{aligned}
\end{equation*}
Consequently,
\begin{equation*}
\begin{aligned}
\|\bM X_n^{\mathrm{BUG},(j)}\|_2^4&\leq(1+\Delta t\|\bA\|_2)^4
\|X_n^{\mathrm{BUG},(j)}\|_2^4\\
&\leq(1+c_2\Delta t)\|X_n^{\mathrm{BUG},(j)}\|_2^4,
\end{aligned}
\end{equation*}
where $c_2=4\|\bA\|_2(1+T\|\bA\|_2)^3$. Combining these estimates yields
\begin{equation*}
\E\left[\|X_{n+1}^{\mathrm{BUG},(j)}\|_2^4\mid\mathcal{F}_{t_n}\right]\leq(1+c_3\Delta t)\|X_n^{\mathrm{BUG},(j)}\|_2^4+c_3\Delta t,
\end{equation*}
where $c_3=c_2+12c_1(1+T\|\bA\|_2)^2+6Tc_1^2$. Taking expectations, averaging over the particles, and applying the discrete Gr\"onwall inequality give
\begin{equation}
\label{eq:analysis-mc-quadratic-fourth-moment}
\begin{aligned}
\max_{0\le n\le N_t}\frac1N\sum_{j=1}^N\E\left[\|X_n^{\mathrm{BUG},(j)}\|_2^4\right]&\leq e^{c_3T}\left(\frac1N\sum_{j=1}^N\E\left[\|X_0^{\mathrm{BUG},(j)}\|_2^4\right]+c_3T\right)\\
&\leq c_4,
\end{aligned}
\end{equation}
where $c_4=e^{c_3T}(\E[\|X(0)\|_2^4]+c_3T)$, and the last inequality follows from the contractivity of the orthogonal projection giving the initial low-rank ensemble $\bX_0^{\mathrm{BUG}}$. Consequently, Jensen's inequality yields
\begin{equation*}
\max_{0\le n\le N_t}\E\left[\|\bX_n^{\mathrm{BUG}}\|_N^4\right]\leq c_4.
\end{equation*}
The same fourth-moment estimate also holds for $\bX_n^{\mathrm{EM}}$ by taking $\bP_{n+1}=\bI$.

\emph{Step 2: Conditional second-moment fluctuations.}
Let $\mathcal G=\sigma(\bX_0^{\mathrm{EM}})$ denote the sigma-algebra generated by the initial Euler--Maruyama ensemble. Since the initial low-rank ensemble $\bX_0^{\mathrm{BUG}}$ is obtained from $\bX_0^{\mathrm{EM}}$ and the closed second-moment recursion is deterministic, the projectors $\widehat{\bP}_{n+1}$ and $\bP_{n+1}$ are $\mathcal G$-measurable. Consequently, the conditional second-moment recurrence gives
\begin{equation*}
\begin{aligned}
&\E\left[\bC_{N,n+1}^{\mathrm{BUG}}\mid\mathcal G\right]\\
&\qquad=\bP_{n+1}\left(\bM\E\left[\bC_{N,n}^{\mathrm{BUG}}\mid\mathcal G\right]\bM^\top+\Delta t\bB_0\bB_0^\top+\Delta t\sum_{\ell=1}^q\bB_\ell\E\left[\bC_{N,n}^{\mathrm{BUG}}\mid\mathcal G\right]\bB_\ell^\top\right)\bP_{n+1}.
\end{aligned}
\end{equation*}
In particular, this corresponds to the recurrence defining the closed second moment $\bC_n^{\mathrm{BUG}}=\bU_n\bD_n\bU_n^\top$. Since $\bD_0=\bY_0^\top\bY_0/N$, we have $\bC_0^{\mathrm{BUG}}=\bC_{N,0}^{\mathrm{BUG}}$, and it follows by induction that
\begin{equation}
\label{eq:analysis-mc-quadratic-closed-conditional}
\E\left[\bC_{N,n}^{\mathrm{BUG}}\mid\mathcal G\right]=\bC_n^{\mathrm{BUG}}.
\end{equation}
Moreover, since the Brownian motions are independent across particles and the projectors are $\mathcal G$-measurable, the closed-moment EM--BUG particles are independent conditionally on $\mathcal G$. Let $\bJ_{n,j}=X_n^{\mathrm{BUG},(j)}(X_n^{\mathrm{BUG},(j)})^\top$ and $\widetilde{\bJ}_{n,j}=\bJ_{n,j}-\E\left[\bJ_{n,j}\mid\mathcal G\right]$. Using \eqref{eq:analysis-mc-quadratic-closed-conditional}, we obtain
\begin{equation}
\label{eq:analysis-mc-quadratic-closed-variance}
\begin{aligned}
\E\left[\|\bC_{N,n}^{\mathrm{BUG}}-\bC_n^{\mathrm{BUG}}\|_{\mathrm F}^2\mid\mathcal G\right]&=\E\left[\bigg\|\frac1{N}\sum_{j=1}^N\widetilde{\bJ}_{n,j}\bigg\|_{\mathrm F}^2\mid\mathcal G\right]\\
&=\frac1{N^2}\sum_{i,j=1}^N\E\left[\langle\widetilde{\bJ}_{n,i},\widetilde{\bJ}_{n,j}\rangle_{\mathrm F}\mid\mathcal G\right]\\
&=\frac1{N^2}\sum_{j=1}^N\E\left[\|\widetilde{\bJ}_{n,j}\|_{\mathrm F}^2\mid\mathcal G\right]\\
&=\frac1{N^2}\sum_{j=1}^N\left(\E\left[\|\bJ_{n,j}\|_{\mathrm F}^2\mid\mathcal G\right]-\|\E[\bJ_{n,j}\mid\mathcal G]\|_{\mathrm F}^2\right)\\
&\leq\frac1{N^2}\sum_{j=1}^N\E\left[\|X_n^{\mathrm{BUG},(j)}\|_2^4\mid\mathcal G\right],
\end{aligned}
\end{equation}
where the third equality follows from the conditional independence and centering of the matrices $\widetilde{\bJ}_{n,j}$, and the last inequality follows from $\|\bJ_{n,j}\|_{\mathrm F}^2=\|X_n^{\mathrm{BUG},(j)}\|_2^4$. Taking expectations and using \eqref{eq:analysis-mc-quadratic-fourth-moment} yields
\begin{equation}
\label{eq:analysis-mc-quadratic-closed-fluctuation}
\max_{0\le n\le N_t}\E\left[\|\bC_{N,n}^{\mathrm{BUG}}-\bC_n^{\mathrm{BUG}}\|_{\mathrm F}^2\right]\leq\frac{c_4}{N}.
\end{equation}

\emph{Step 3: Decomposition and fluctuation estimate.}
Let
\begin{equation*}
\bE_0^{\mathrm{init}}=\bC_{N,0}^{\mathrm{EM}}-\bC_0^{\mathrm{EM}},\qquad \bE_{n+1}^{\mathrm{init}}=\bM\bE_n^{\mathrm{init}}\bM^\top+\Delta t\sum_{\ell=1}^q\bB_\ell\bE_n^{\mathrm{init}}\bB_\ell^\top.
\end{equation*}
Comparing the closed BUG second-moment recurrence with the Euler--Maruyama second-moment recurrence and using \eqref{eq:analysis-mc-quadratic-lemma-loss}, we obtain by induction
\begin{equation*}
\bC_n^{\mathrm{BUG}}=\bC_n^{\mathrm{EM}}+\bE_n^{\mathrm{init}}-\bGamma_n.
\end{equation*}
Consequently, setting
\begin{equation*}
\bXi_n=\bE_n^{\mathrm{init}}+\bC_{N,n}^{\mathrm{BUG}}-\bC_n^{\mathrm{BUG}}
\end{equation*}
gives \eqref{eq:analysis-mc-quadratic-decomposition}. Since $\bGamma_0\succeq0$ by the truncated SVD initialization and $\bL_{n+1}\succeq0$ under the complete-augmentation assumption, the recurrence \eqref{eq:analysis-mc-quadratic-lemma-loss} implies by induction that $\bGamma_n\succeq0$ for all $n$.

Moreover, since the initial Euler--Maruyama particles are independent and identically distributed, we have
\begin{equation*}
\E[\bE_0^{\mathrm{init}}]=\E\left[\bC_{N,0}^{\mathrm{EM}}\right]-\bC_0^{\mathrm{EM}}=0.
\end{equation*}
Since the recurrence of $\bE_n^{\mathrm{init}}$ is linear and deterministic, we obtain by induction $\E[\bE_n^{\mathrm{init}}]=0$ for all $n$. It follows from \eqref{eq:analysis-mc-quadratic-closed-conditional} that $\E[\bXi_n]=0$. Furthermore, since $\bE_n^{\mathrm{init}}$ is $\mathcal G$-measurable, we have
\begin{equation*}
\begin{aligned}
\E\left[\left\langle\bE_n^{\mathrm{init}},\bC_{N,n}^{\mathrm{BUG}}-\bC_n^{\mathrm{BUG}}\right\rangle_{\mathrm F}\right]&=\E\left[\left\langle\bE_n^{\mathrm{init}},\E\left[\bC_{N,n}^{\mathrm{BUG}}-\bC_n^{\mathrm{BUG}}\mid\mathcal G\right]\right\rangle_{\mathrm F}\right]\\
&=0,
\end{aligned}
\end{equation*}
where the last equality follows from \eqref{eq:analysis-mc-quadratic-closed-conditional}. Consequently, the cross term vanishes, and we obtain
\begin{equation}
\label{eq:analysis-mc-quadratic-fluctuation-decomposition}
\E\left[\|\bXi_n\|_{\mathrm F}^2\right]=\E\left[\|\bE_n^{\mathrm{init}}\|_{\mathrm F}^2\right]+\E\left[\|\bC_{N,n}^{\mathrm{BUG}}-\bC_n^{\mathrm{BUG}}\|_{\mathrm F}^2\right].
\end{equation}

For $\bM=\bI+\Delta t\bA$, the estimates $\|\bM\|_2\le 1+\Delta t\|\bA\|_2$ and $\Delta t\le T$ give
\begin{equation*}
\|\bM\|_2^2\leq(1+\Delta t\|\bA\|_2)^2\leq1+\Delta t(2\|\bA\|_2+T\|\bA\|_2^2).
\end{equation*}
Consequently,
\begin{equation*}
\begin{aligned}
\left\|\bM\bZ\bM^\top+\Delta t\sum_{\ell=1}^q\bB_\ell\bZ\bB_\ell^\top\right\|_{\mathrm F}&\leq\left(\|\bM\|_2^2+\Delta t\sum_{\ell=1}^q\|\bB_\ell\|_2^2\right)\|\bZ\|_{\mathrm F}\\
&\quad\leq(1+c_5\Delta t)\|\bZ\|_{\mathrm F},
\end{aligned}
\end{equation*}
where $c_5=2\|\bA\|_2+T\|\bA\|_2^2+\sum_{\ell=1}^q\|\bB_\ell\|_2^2$. Moreover, since $\Delta t\leq T$, we have
\begin{equation*}
(1+c_5\Delta t)^2\leq1+(2c_5+Tc_5^2)\Delta t.
\end{equation*}
Consequently,
\begin{equation*}
\left\|\bM\bZ\bM^\top+\Delta t\sum_{\ell=1}^q\bB_\ell\bZ\bB_\ell^\top\right\|_{\mathrm F}^2
\leq(1+c_6\Delta t)\|\bZ\|_{\mathrm F}^2,
\end{equation*}
where $c_6=2c_5+Tc_5^2$. Moreover, since the initial Euler--Maruyama particles are independent and identically distributed, the same calculation as in \eqref{eq:analysis-mc-quadratic-closed-variance} gives
\begin{equation*}
\begin{aligned}
\E\left[\|\bE_0^{\mathrm{init}}\|_{\mathrm F}^2\right]&=\frac1N\E\left[\|X(0)X(0)^\top-\bC_0^{\mathrm{EM}}\|_{\mathrm F}^2\right]\\
&=\frac1N\left(\E\left[\|X(0)\|_2^4\right]-\|\bC_0^{\mathrm{EM}}\|_{\mathrm F}^2\right)\\
&\le\frac1N\E\left[\|X(0)\|_2^4\right].
\end{aligned}
\end{equation*}
The discrete Gr\"onwall inequality then yields
\begin{equation*}
\max_{0\le n\le N_t}\E\left[\|\bE_n^{\mathrm{init}}\|_{\mathrm F}^2\right]\leq e^{c_6T}\E\left[\|\bE_0^{\mathrm{init}}\|_{\mathrm F}^2\right]\leq\frac{e^{c_6T}}{N}\E\left[\|X(0)\|_2^4\right].
\end{equation*}

By \eqref{eq:analysis-mc-quadratic-fluctuation-decomposition}, combining this estimate with \eqref{eq:analysis-mc-quadratic-closed-fluctuation} gives
\begin{equation*}
\label{eq:analysis-mc-quadratic-closed-Z}
\max_{0\le n\le N_t}\E\left[\|\bXi_n\|_{\mathrm F}^2\right]\leq\frac{c_7}{N},
\end{equation*}
where $c_7=e^{c_6T}\E[\|X(0)\|_2^4]+c_4$. For every deterministic symmetric matrix $\bK$, the Cauchy--Schwarz inequality gives
\begin{equation*}
\E\left[|\operatorname{tr}(\bK\bXi_n)|^2\right]\leq\|\bK\|_{\mathrm F}^2\E\left[\|\bXi_n\|_{\mathrm F}^2\right]\leq\frac{c_7\|\bK\|_{\mathrm F}^2}{N}.
\end{equation*}
Finally, setting $\sigma_{T,K}=\sqrt{c_7}\|\bK\|_{\mathrm F}$ proves \eqref{eq:analysis-mc-quadratic-fluctuation}.
\end{proof}

\end{appendices}

\backmatter

\section*{Declarations}

\bmhead{Data availability}

The source code used for the numerical experiments is publicly available at \url{https://github.com/sebastienriffaud/em-bug}. The datasets generated during and/or analyzed during the current study are available from the corresponding author on reasonable request.

\bmhead{Conflict of interest}

The authors declare that they have no conflict of interest.

\bmhead{Acknowledgements}

The author would like to thank Dr. T. Trigo Trindade and Dr. F. Zoccolan for fruitful discussions. This work has been supported by the Swiss National Science Foundation under the Project n°200518 "Dynamical low rank methods for uncertainty quantification and data assimilation".

\bibliography{sn-bibliography}

\end{document}